\documentclass[11pt, a4paper]{amsart}
\title{The Morrison cone conjecture under deformation}
\author{Wendelin Lutz}
\address{Leibniz University Hannover, Institute of Algebraic Geometry, Welfengarten 1,
	30167 Hannover, Germany.}
\email{lutz@math.uni-hannover.de}
\usepackage{amsfonts, amssymb}
\usepackage{tikz}
\usepackage{tikz-cd}
\usepackage{booktabs}
\usepackage[margin=2.5cm]{geometry}
\usepackage{hyperref}
\usepackage{xfrac}
\usepackage{epsfig}
\usepackage{setspace}
\usepackage[shortlabels]{enumitem}
\usepackage{xcolor}
\usetikzlibrary{matrix}

\usepackage{amsmath}

\newcommand{\Aut}{\mathrm{Aut}}
\newcommand{\Bbig}{\mathrm{Big}}
\newcommand{\Bir}{\mathrm{Bir}}

\newcommand{\Def}{\mathrm{Def}}

\newcommand{\Eff}{\mathrm{Eff}}

\newcommand{\GL}{\mathrm{GL}}

\newcommand{\id}{\mathrm{id}}

\newcommand{\Mov}{\mathrm{Mov}}
\newcommand{\NE}{\mathrm{NE}}
\newcommand{\Nef}{\mathrm{Nef}}
\newcommand{\NS}{\mathrm{NS}}

\newcommand{\oAut}{\mathrm{\overline{\Aut}}}
\newcommand{\oPsAut}{\mathrm{\overline{\PsAut}}}

\newcommand{\Pic}{\mathrm{Pic}}

\newcommand{\PsAut}{\mathrm{PsAut}}

\newcommand{\Span}{\mathrm{Span}}

\newcommand{\Supp}{\mathrm{Supp}}

\newcommand{\cC}{\mathcal{C}}
\newcommand{\cD}{\mathcal{D}}
\newcommand{\cE}{\mathcal{E}}

\newcommand{\cH}{\mathcal{H}}

\newcommand{\cL}{\mathcal{L}}

\newcommand{\cO}{\mathcal{O}}

\newcommand{\cU}{\mathcal{U}}
\newcommand{\cV}{\mathcal{V}}

\newcommand{\cX}{\mathcal{X}}
\newcommand{\cY}{\mathcal{Y}}
\newcommand{\cZ}{\mathcal{Z}}

\newcommand{\CC}{\mathbb{C}}

\newcommand{\PP}{\mathbb{P}}
\newcommand{\QQ}{\mathbb{Q}}
\newcommand{\RR}{\mathbb{R}}

\newcommand{\ZZ}{\mathbb{Z}}

\newcommand{\ip}{\raise1pt\hbox{\large $\lrcorner$}}

\theoremstyle{plain}
\newtheorem{theorem}{Theorem}[section]
\newtheorem{corollary}[theorem]{Corollary}
\newtheorem{lemma}[theorem]{Lemma}
\newtheorem{proposition}[theorem]{Proposition}
\newtheorem*{theorem*}{Theorem}
\newtheorem*{conjecture}{Conjecture}

\theoremstyle{definition}
\newtheorem{definition}[theorem]{Definition}
\newtheorem{example}[theorem]{Example}
\newtheorem{remark}[theorem]{Remark}
\newtheorem*{remark*}{Remark}

\begin{document}
	\begin{abstract}
We prove that the Morrison--Kawamata cone conjecture is invariant under deformation for smooth Calabi--Yau threefolds. The proof compares the nef and movable cones of $Y$ to those of a very general deformation of $Y$ via Weyl groups naturally associated to $Y$. We use this result to prove new cases of the cone conjecture.
	\end{abstract}
\maketitle
\section{Introduction}
A fundamental invariant of a projective variety $Y$ is its cone of curves, the convex cone spanned by numerical equivalence classes of curves on $Y$. Its dual is the cone $\Nef(Y)$, the cone of numerical classes of nef line bundles. In general, the geometry of $\Nef(Y)$ is not well-understood. The simplest behaviour occurs for Fano manifolds: by the cone theorem, $\Nef(Y)$ is rational polyhedral, i.e. it is the cone over a rational polytope. 
If $Y$ is a Calabi--Yau manifold, by which we mean that $Y$ has trivial canonical bundle and $h^i(\cO_Y)=0$ for $0 <i<\dim(Y)$, the nef cone can already exhibit complicated geometric behaviour, such as having infinitely many faces. 
Nevertheless, the Morrison cone conjecture, originally formulated in 1993, predicts:
\begin{conjecture}[\cite{Morrison} nef cone conjecture]
	Let $Y$ be a Calabi--Yau manifold.
The action of the automorphism group $\Aut(Y)$ on $\Nef^+(Y)$ admits a rational polyhedral fundamental domain (RPFD).
\end{conjecture} 
Here, $\Nef^+(Y)$ denotes the convex hull of the rational points of $\Nef(Y)$.
One very surprising consequence of the cone conjecture is that a Calabi--Yau manifold whose nef cone is not rational polyhedral must have infinite automorphism group. \\
The cone conjecture was generalized by Kawamata to the closure of the movable cone $\overline{\Mov}(Y)$, which can roughly be thought of as an enlargement of the nef cone preserved under small modifications. We recall that a small modification is a birational map that is an isomorphism in codimension one, such as a flop. 
A birational self-map $Y \dashrightarrow Y$ which is an isomorphism in codimension one is called a pseudo-automorphism. We denote by $\PsAut(Y)$ the group of pseudo-automorphisms of $Y$. 
\begin{conjecture}[\cite{KawamataCY} movable cone conjecture]
	Let $Y$ be a Calabi--Yau manifold.
The action of $\PsAut(Y)$ on ${\Mov}^+(Y)$ admits a rational polyhedral fundamental domain.
\end{conjecture}
The cone conjecture has been verified for surfaces by Sterk~\cite{Sterk}, following suggestions of Looijenga. In dimension $n>2$, the conjecture has been verified in several nontrivial examples (see~\cite{Fryers, GrassiMorrison, Gachet, CantatOguiso, Yanez, Lazic, Ito}). There has recently also been progress on the closely related cone conjecture for klt log Calabi--Yau pairs (see \cite{Totaro, Li, GachetStenger, Xu}). Nevertheless, the cone conjecture for Calabi--Yau manifolds of dimension $n>2$ remains wide open.
Our main result is that in dimension three, the cone conjecture for smooth Calabi--Yau threefolds is invariant under deformation:
\begin{theorem}\label{thm:mainA} (see Theorem~\ref{thm:MCClocal})
	If the movable cone conjecture, respectively the nef cone conjecture, holds for a smooth Calabi--Yau threefold $Y$, then it also holds for any smooth Calabi--Yau threefold deformation-equivalent to $Y$.
\end{theorem}
Our theorem suggests a path to proving the cone conjecture in a wider range of examples, by deforming $Y$ to a smooth threefold for which the cone conjecture is easier to prove, for instance because the movable and/or the nef cone are rational polyhedral. Using this strategy, we give a new, very short proof of the cone conjecture for any smooth $(2,2,2,2)$-hypersurface in $(\PP^1)^4$ (for a very general hypersurface, this was first proved by Cantat--Oguiso~\cite{CantatOguiso}).
We also use our result to show that the cone conjecture holds for the deformation families of some of the Calabi--Yau threefolds studied by Borcea and Voisin~\cite{VoisinMiroirs, Borcea}. 
\subsection*{Outline of proof}
To prove Theorem~\ref{thm:mainA}, it suffices to prove the following local statement: the nef cone conjecture (resp. movable cone conjecture) holds for a smooth Calabi--Yau threefold $Y$ if and only if it holds for a very general fiber $Y^{\mathrm{gen}}$ of the Kuranishi family of $Y$. Applying this statement on overlapping analytic neighbourhoods then gives invariance on the whole deformation class. 
We achieve this by relating the nef and movable cones, as well as the automorphism and pseudo-automorphism groups of $Y$ and $Y^{\mathrm{gen}}$.\\
Our first result (see Section~\ref{sec:deformation}) concerns the behaviour of $\overline{\Mov}(Y)$ under deformation, which is the analogue of a theorem of Wilson~\cite{WilsonCone} on $\Nef(Y)$. To state it precisely, we make the following definitions:
\begin{definition}\label{def:qruled}
	Let $E$ be a projective Cohen--Macaulay irreducible surface.
	\begin{itemize}
		\item $E$ is a conic fibration surface if there is a flat surjective morphism $p \colon E \rightarrow C$ with connected fibers to a smooth curve $C$, and the general fiber of $p$ is either $\PP^1$ or the union of two $\PP^1$'s meeting in a node. We say that $E$ is quasi-ruled if in addition $p$ is locally trivial in the \'etale topology.
		\item Suppose that $E$ is contained in a smooth projective threefold $Y$. We say that $E$ is an SQM-conic fibration surface if there is a small $\QQ$-factorial modification (SQM) 
		\[
		\alpha \colon Y \dashrightarrow Z
		\]
		 such that the strict transform of $E$ admits the structure of a conic fibration surface
		 \[
		 p \colon E_Z \rightarrow C.
		 \]
		 In particular, if $E \subset Y$ is a conic fibration surface, then it is an SQM-conic fibration surface by taking $\alpha$ the identity. We call the pair $(\alpha, p)$ an SQM-model of $E$. 
		 We say that an SQM-conic fibration surface $E\subset Y$ is over a curve of positive genus if it admits an SQM-model
		 \[
		 (\alpha\colon Y\dashrightarrow Z,  p\colon E_Z\to C)
		 \]
		 with $g(C)>0$.
		\end{itemize}
\end{definition}
For a smooth Calabi--Yau threefold $Y$, the vanishing of $h^i(\cO_Y)$ induces an isomorphism $\Pic(Y) \cong H^2(Y, \ZZ)$. If $\pi \colon \cY \rightarrow B$ is a smooth family of Calabi--Yau threefolds over a polydisk $B$, we may identify $H^2(\cY_b) \cong H^2(\cY)$, where $\cY_b$ denotes the fiber of $\pi$ over $b \in B$. Moreover, the nef cones of fibers (resp. movable cones) are constant over a dense subset of $B$. In view of this, we may compare the nef and movable cones of the fibers $\cY_b$ in the fixed vector space $H^2(Y, \RR)$, and may unambiguously speak about the nef and movable cones of $Y^{\mathrm{gen}}$.
\begin{theorem}\label{thm:main1}
		Let $\pi \colon \cY \rightarrow \Def(Y)$ be the Kuranishi family of a smooth Calabi--Yau threefold. 
\begin{itemize}
	\item 
The nef cone is constant in fibers of $\pi$ if and only if no fiber $Y_b, b \in \Def(Y)$ contains a surface $E$ quasi-ruled over a curve of genus $1$ \cite{WilsonCone, WilsonErratum}.
	\item 
The movable cone is constant in fibers of $\pi$ if and only if no fiber $Y_b, b \in \Def(Y)$ contains an SQM-conic fibration surface over a curve of positive genus.
\end{itemize}	
\end{theorem}
Theorem~\ref{thm:main1} shows that the nef and movable cones can only jump in the presence of certain special surfaces $E$. The jumping is in fact measured by a Coxeter group:
inspired by Wilson's work, we introduce a reflection group $W_Y^{big}$, which we call the big Wilson Weyl group. Let $E \subset Y$ be a SQM-conic fibration surface over a curve of positive genus, and let $(\alpha \colon Y \dashrightarrow Z, p \colon E_Z \rightarrow C)$ be an SQM-model with $g(C)>0$. Let $\ell_Z \in H_2(Z, \ZZ)$ be the class of a fiber of $p$ and let $\ell_E=\alpha_*^{-1}(\ell_Z) \in H_2(Y, \ZZ)$. We will show in Lemma~\ref{lem:propertiesofSQMconicfibration} that the isomorphism class of $C$ and the class $\ell_E$ are intrinsic to $E$. We define
\[
\sigma_{E} \colon H^2(Y, \ZZ) \rightarrow H^2(Y, \ZZ), \quad x \mapsto x+(x\cdot \ell_E)[E].
\] 
Since $E \cdot \ell_E=-2$, $\sigma_E$ fixes the hyperplane $\ell_E^\perp$ pointwise, sends $[E]$ to $-[E]$, and is an involution. We call $\sigma_E$ the reflection associated to $E$.  However, unlike in the case of surfaces, there is no canonical symmetric bilinear form on $H^2(Y, \ZZ)$.
The group $W_Y^{big}$ is generated by the reflections $\sigma_E$ as $E$ ranges over SQM-conic fibration surfaces $E$ over curves of positive genus. The small Weyl group $W_Y^{sm}$ is the subgroup generated by $\sigma_E$ with $E$ quasi-ruled over a curve of genus $1$ (cf. \cite{WilsonElliptic}).\\
We will see in Section~\ref{sec:roots} that we may view the pairs $([E], \ell_E)$ as simple roots in a generalized sense. 
Even though the classes $[E]$ of SQM-conic fibration surfaces are usually not linearly independent in $H^2(Y, \RR)$, classical results of Bourbaki and Looijenga~\cite{Bourbaki, Looijenga} remain valid, and we combine them with the minimal model program to show:
\begin{theorem}\label{thm:main2}
	\hfill
	\begin{enumerate}
		\item $\Nef(Y) \cap \Bbig(Y)$ is a fundamental domain for the action of $W_Y^{sm}$ on $\Nef(Y^{\mathrm{gen}}) \cap \Bbig(Y^{\mathrm{gen}})$. 
		\item 	
		$\overline{\Mov}(Y) \cap \Bbig(Y)$ is a fundamental domain for the action of $W_Y^{big}$ on $\overline{\Mov}(Y^{\mathrm{gen}}) \cap \Bbig(Y^{\mathrm{gen}})$. 
	\end{enumerate}
	Here, $\Bbig(Y)$ is the cone of big divisors of $Y$. 
\end{theorem}
Under the assumption that $W_Y^{sm}$ is finite, Theorem~\ref{thm:main2}(1) was proved in \cite[Corollary]{WilsonElliptic}.
Although Theorem~\ref{thm:main2} only compares nef and movable cones after intersecting with the big cone, this is sufficient for the proof of Theorem~\ref{thm:mainA}.
The group $W_Y^{sm}$ is the three-dimensional analogue of the group generated by reflections in $(-2)$-curves on a $K3$ surface. Indeed, we show that elements of $W_Y^{sm}$ are realized by monodromy transformations, in analogy with the Picard-Lefschetz reflection in dimension $2$. The group $W_Y^{big}$ has no direct analogue in dimension $2$, and its elements are realized geometrically as composites of parallel transport in a family of Calabi--Yau threefolds, followed by flops. See Section~\ref{sec:monodromy}. \\
The group $W_Y^{big}$ gives a geometric explanation for the appearance of abstract Coxeter groups in the work of Cantat--Oguiso~\cite{CantatOguiso} on Calabi--Yau hypersurfaces. Indeed, we show in Section~\ref{sec:CO} that in dimension three, the group appearing in their work can be interpreted as the big Wilson Weyl group of a suitable degeneration, and their Theorem~4.1 can be seen as an instance of Theorem~\ref{thm:main2}.\\
In Section~\ref{sec:H3} and \ref{sec:MCC}, we use $W_Y^{sm}$ (resp. $W_Y^{big}$) to compare the automorphism groups (resp. pseudo-automorphism groups) of $Y$ and $Y^{\mathrm{gen}}$. If $\pi \colon \cY \rightarrow S$ denotes a very general deformation of $Y$, and $\cY_\eta$ is the generic fiber, we show that there is a natural inclusion 
\[
\Aut(\cY_\eta) \rightarrow W_Y^{sm} \rtimes \Aut(Y).
\]
The group $W_Y^{sm} \rtimes \Aut(Y)$ naturally acts on the deformation space $\Def(Y)$, and Hodge-theoretic considerations (Proposition~\ref{prop:actsfinitegroup}) show that it must act via a finite group. In particular, a finite index subgroup of $W_Y^{sm} \rtimes \Aut(Y)$ lifts to automorphisms of any small deformation $Y'$ of $Y$, a result which may be of independent interest. \\
Using the foundational work of Looijenga~\cite{LooijengaCones}, we are then able to show that the nef cone conjecture for $Y$ is equivalent to the nef cone conjecture for $Y^{\mathrm{gen}}$. \\
In a similar fashion, we use $W_Y^{big}$ to relate the pseudo-automorphism groups of $Y$ and $Y^{\mathrm{gen}}$ to show that the movable cone conjecture for $Y$ is equivalent to the movable cone conjecture for $Y^{\mathrm{gen}}$. \\
We remark that the cone conjecture is often stated using the effective nef cone $\Nef^e(Y)$ (resp. the effective movable cone $\Mov^e(Y)$). The equivalence of the two formulations is conjectural, and reduces to the question whether every rational class on the boundary of $\Nef^+(Y)$ or $\Mov^+(Y)$ on a Calabi--Yau threefold is effective. The generalized abundance conjecture (see for example \cite{LazicAbundance}) would imply that this is true for the nef cone, but as far as we know, a similar statement for the movable cone is not even conjectured. Both of these questions seem very hard, and this paper offers no new results in that direction.
\subsection*{Acknowledgements}
I am very grateful to Paul Hacking for his help with technical aspects of this paper. I would also like to thank Eyal Markman for helpful comments on a prior version of this paper, and Stefano Filipazzi for help with Proposition~\ref{prop:flopofroot}(2).
During the revision process, ChatGPT was used for proofreading, to suggest improvements to the notation in Sections~4--6, and to assist in improving the exposition of some proofs. It also identified omitted hypotheses in Definition~\ref{def:genroots} and Proposition~\ref{prop:doublecover} and several minor textual inaccuracies. All AI-generated suggestions were verified by the author, and these revisions did not affect the results of the paper. 
\section{Preliminaries}
Throughout this paper, $Y$ will denote a Calabi--Yau threefold in the strict sense. By this, we mean that $Y$ is a smooth projective threefold over the field of complex numbers with trivial canonical bundle satisfying $h^1(Y, \cO_Y)=h^2(Y, \cO_Y)=0$.
$S$ will denote either a smooth variety, or a smooth germ of a complex space. A family $\pi \colon \cY \rightarrow S$ of Calabi--Yau threefolds is a smooth projective morphism of varieties (resp. complex spaces). \\
Since $H^2(Y, \CC)=H^{1,1}(Y)$, we have isomorphisms $\NS(Y)_\RR \cong \Pic(Y)_\RR \cong H^2(Y, \RR)$, and the space of numerical equivalence classes of $1$-cycles on $Y$ is identified with $H_2(Y, \RR)$. 
In view of this, we will view all the cones commonly used in algebraic geometry as cones in $H^2(Y, \RR)$ or $H_2(Y, \RR)$.
Given any cone $C$, we denote by $\bar{C}$ its closure, and $C^{\circ}$ its interior. 
The Mori cone $\NE(Y)$ is the convex cone in $H_2(Y, \RR)$ spanned by the set of effective curves on $Y$. 
The nef cone $\Nef(Y) \subset H^2(Y, \RR)$ is the dual cone of $\overline{\NE}(Y)$, and by the Kleiman criterion, the ample cone $A(Y)$ is identified with the interior of $\Nef(Y)$. 
Recall that an effective divisor $D$ is movable if the stable base locus of $D$ has codimension at least $2$ in $Y$. We denote by $\Mov(Y) \subset H^2(Y, \RR)$ the cone spanned by movable divisors. In general, $\Mov(Y)$ is neither open nor closed, and we denote $\overline{\Mov}(Y)$ its closure. 
The effective cone $\Eff(Y)$ is the cone spanned by effective divisor classes (which is also not open or closed in general), and its interior is the cone of big divisors $\Bbig(Y)$.
We have inclusions
\begin{align*}
	A(Y)  &\subset \Mov(Y)^{\circ} \subset \Bbig(Y)\\
	{\Nef}(Y) & \subset \overline{\Mov}(Y) \subset \overline{\Eff}(Y).
\end{align*}
Integral points on the boundary of $\Nef(Y)$ (resp. $\overline{\Mov}(Y)$) are not always represented by effective divisor classes, so we define the effective nef and movable cones
\[
{\Nef}^e(Y)=\Nef(Y) \cap \Eff(Y), \qquad {\Mov}^e(Y)=\overline\Mov(Y) \cap \Eff(Y)
\]
A small $\QQ$-factorial modification (SQM) of $Y$ is a birational map $\alpha \colon Y \dashrightarrow Z$ with $Z$ projective and $\QQ$-factorial which is an isomorphism in codimension one. It is well-known (see for example \cite[Theorem 4.9]{KollarFlops}) that any SQM between smooth Calabi--Yau threefolds can be decomposed into a composite of flops.
A SQM $\alpha \colon Y \dashrightarrow Z$ induces an isomorphism $\alpha^* \colon H^2(Z, \ZZ) \rightarrow H^2(Y, \ZZ)$ by taking strict transform of divisors. The strict transform of a movable divisor is again movable, so that $\alpha^*(\Mov(Z))=\Mov(Y)$. Similar equalities hold for ${\Mov}^e(Y), \overline{\Mov}(Y), \Bbig(Y), \Eff(Y), \overline{\Eff}(Y)$. 
 The nef cone is not preserved under flops: If $\alpha \colon Y \dashrightarrow Z$ flops a curve $C$, then the intersection of $\alpha^*(\Nef(Z))$ and $\Nef(Y)$ is contained in the hyperplane $[C]^\perp \subset H^2(Y, \RR)$. By \cite[Theorem 2.3]{KawamataCY}, we have a decomposition:
\begin{equation}\label{MovingConeDecomposition}
	{\Mov}^e(Y)= \bigcup_{\alpha \colon Y \dashrightarrow Z} \alpha^*\Nef^e(Z)
\end{equation}
where the union is over all SQMs of $Y$. We emphasize that it might well happen that several or even all of the varieties $Z$ are abstractly isomorphic to $Y$, so $\alpha$ cannot be omitted from notation.
We define $\Nef^+(Y)$ to be the convex hull of the rational points of $\Nef(Y)$, and ${\Mov}^+(Y)$ to be the convex hull of the rational points of $\overline{\Mov}(Y)$. 
As a consequence of the log abundance theorem for threefolds we have the following inclusions of cones:
\begin{lemma}\cite[Proposition 2.4]{KawamataCY}
	Suppose that $Y$ is a smooth Calabi--Yau threefold. Then
	\begin{itemize}
		\item $\Nef^e(Y) \subset \Nef^+(Y) $
		\item ${\Mov}^e(Y) \subset {\Mov}^+(Y)$
	\end{itemize}
\end{lemma}
The nef cone and the movable cone of a Calabi--Yau threefold can exhibit complicated behaviour, for instance, they can fail to be even locally rational polyhedral. However, after intersecting with the big cone, this ``bad" behaviour disappears:
\begin{proposition}\cite[Corollary 2.7]{KawamataCY}\label{prop:big-rational}
	Let $Y$ be a smooth Calabi--Yau variety. The intersection $\overline{\Mov}(Y) \cap \Bbig(Y)$ is locally rational polyhedral inside the open cone $\Bbig(Y)$.\\
	Every codimension $k$ face $F$ of $\Bbig(Y) \cap \overline{\Mov}(Y)$ corresponds to a contraction morphism $p_F \colon Z \rightarrow  \bar{Z}$ of relative Picard rank $k$ of some SQM $\alpha\colon Y \dashrightarrow Z$ of $Y$.
\end{proposition}
We will call a contraction morphism of relative Picard rank $1$ primitive.
Here, the contraction morphism $p_F$ is the map defined by the linear system $|mD|$, for $D$ a $\QQ$-Cartier divisor in the relative interior of $F$, and $m$ sufficiently large and divisible. Note that by \eqref{MovingConeDecomposition}, $D$ is big and nef on some SQM $\alpha \colon Y \dashrightarrow Z$, so that $p_F \colon Z \rightarrow \bar{Z}$ is a morphism by the base-point free theorem. 
It is easy to show that $p_F$ as above is always a divisorial contraction:
\begin{lemma}\label{lem:MovingConeFacesAreDivisorial}
	Let $Y$ be a smooth Calabi--Yau threefold, let $F$ be a codimension one face of $\overline{\Mov}(Y) \cap \Bbig(Y)$, and let $\alpha \colon Y \dashrightarrow Z$ be a SQM such that $F$ is a face of $\alpha^*\Nef(Z)$. Then the associated contraction morphism 
	$p_F \colon Z \rightarrow \bar{Z}$ is a divisorial contraction.  
\end{lemma}
\begin{proof}
	Since $F$ is contained in $\Bbig(Y)$, $p_F$ is birational by definition, and we only need to show that $p_F$ is not small. But by \cite[Lemma 2.3]{KawamataCrepant}, $p_F$ is small if and only if $F^{\circ}$ is contained in the interior of $\overline{\Mov}(Y)$. 
\end{proof}
In view of Lemma~\ref{lem:MovingConeFacesAreDivisorial}, every codimension one face $F$ of $\overline{\Mov}(Y) \cap \Bbig(Y)$ has an associated exceptional divisor, which we denote $E_F$. 
\section{The movable cone under deformation}\label{sec:deformation}
Throughout this section, $S$ denotes a smooth germ of a complex analytic manifold.
\subsection{The behaviour of the movable cone in families}
Recall that deformations of a smooth projective Calabi--Yau threefold $Y$ are unobstructed by the Bogomolov--Tian--Todorov theorem, so the deformation space (Kuranishi space) $\Def(Y)$ is a smooth germ. Moreover, by the vanishing of $h^i(\cO_Y)$ for $0<i<3$, line bundles lift uniquely, so that the Kuranishi family $\pi \colon \cY \rightarrow \Def(Y)$ is projective. 
We write $Y^{\mathrm{gen}}$ for a very general deformation of $Y$.
	Since $h^1(\cO_Y)=h^2(\cO_Y)=0$, the inclusion $i \colon Y_s \rightarrow \cY$ of any fiber induces an isomorphism $i^* \colon \Pic(\cY) \cong \Pic(Y_s)$, which we will refer to as parallel transport. Moreover, we have isomorphisms $H^2(Y, \RR) \cong \Pic(Y)_\RR \cong H^2(Y_s, \RR)$, so we may view $\Nef(Y_s)$ and $\Mov(Y_s)$ as cones in the fixed vector space $H^2(Y ,\RR)$. 
	A standard argument shows that these cones are constant away from a countable union of closed analytic subsets in $\Def(Y)$:
	\begin{proposition}\cite[Proposition 2.12]{KawamataCY}\label{prop:ConesConstant}
		Let $Y$ be a Calabi--Yau threefold, and $\cY \rightarrow S$ a deformation of $Y$ over a germ $S$.
		Then the cone $\Nef(Y_s)$ (respectively $\overline{\Mov}(Y_s)$) ($s \in S$) is constant on the complement $S^0$ of a countable union of closed analytic subsets of $S$.
	\end{proposition}
	In view of this, we define:
	\begin{definition}
		Let $\pi \colon \cY \rightarrow S$ be a very general deformation of a Calabi--Yau threefold $Y$ over a smooth germ. We define $\Nef(Y^{\mathrm{gen}})$ to be the nef cone of a very general fiber of $\pi$. We define $\overline\Mov(Y^{\mathrm{gen}})$ analogously.
	\end{definition}
	We first investigate how the various cones we consider change under deformation. 
	\begin{lemma}\label{lem:ConeInclusions}
		We have the following inclusions:
		\begin{enumerate}
			\item $\Nef(Y) \subset \Nef(Y^{\mathrm{gen}})$
			\item $\Eff(Y^{\mathrm{gen}}) \subset \Eff(Y)$
			\item $\Bbig(Y^{\mathrm{gen}}) \subset \Bbig(Y)$
			\item $\Nef(Y) \cap \Bbig(Y) = \Nef(Y) \cap \Bbig(Y^{\mathrm{gen}})$ 
			\item $\overline{\Mov}(Y) \cap \Bbig(Y) = \overline{\Mov}(Y) \cap \Bbig(Y^{\mathrm{gen}})$
			\item $\overline{\Mov}(Y) \subset \overline{\Mov}(Y^{\mathrm{gen}})$. 
		\end{enumerate}
	\end{lemma}
	\begin{proof}
		We note first that the dual statement $\NE(Y^{\mathrm{gen}}) \subset \NE(Y)$ of $(1)$ is a consequence of the properness of the relative Hilbert scheme of $\pi \colon \cY \rightarrow S$. Taking closures and dualizing, we obtain (1).
		Similarly, we have $\Eff(Y^{\mathrm{gen}}) \subset \Eff(Y)$, yielding (2), and taking interiors, we obtain (3). \\
		For (4), one inclusion is immediate from (3). For the other inclusion, note that if $D$ is big and nef on $Y$, then $D^3 > 0$. By (1), $D \in \Nef(Y^{\mathrm{gen}})$, so $D \in \Bbig(Y^{\mathrm{gen}})$ as well. \\
		For (5), if $D$ is contained in the left hand side, $D \in {\Mov}^e(Y)$, so there is a SQM $\alpha \colon Z \dashrightarrow Y$ such that $\alpha^*D \in \Nef^e(Z)$. We then argue as in (4).\\
		For (6), let $x \in \Mov(Y)$. We again have $\alpha^*(x) \in \Nef(Z)$ for some SQM $\alpha \colon Z \dashrightarrow Y$. Since $Y$ is smooth Calabi--Yau, $\alpha$ is a composition of flops so by \cite[Theorem 11.10]{KollarMoriFlips}, $\alpha$ deforms to a SQM $\tilde{\alpha} \colon \cZ \dashrightarrow \cY$ over $S$, and therefore induces a SQM $\alpha' \colon Z^{\mathrm{gen}} \dashrightarrow Y^{\mathrm{gen}}$. 
		In particular we have that $\alpha^*(\Mov(Y^{\mathrm{gen}}))=\Mov(Z^{\mathrm{gen}})$. By $(1)$ we have $\alpha^*(x)\in \Nef(Z^{\mathrm{gen}}) \subset \overline{\Mov}(Z^{\mathrm{gen}})$, and therefore $x \in \overline{\Mov}(Y^{\mathrm{gen}})$. This gives $\Mov(Y) \subset \overline{\Mov}(Y^{\mathrm{gen}})$, and taking closures we obtain $\overline{\Mov}(Y) \subset \overline{\Mov}(Y^{\mathrm{gen}})$. 
	\end{proof}
\begin{lemma}\label{lem:preservegen}
	Let $\alpha \colon Y \dashrightarrow Z$ be an SQM between smooth Calabi--Yau threefolds. Then \[
	\alpha^*\overline{\Mov}(Z^{\mathrm{gen}})=\overline{\Mov}(Y^{\mathrm{gen}})
	\]
\end{lemma}
\begin{proof}
Let $\cY \rightarrow \Def(Y)$ and $\cZ \rightarrow \Def(Z)$ denote the Kuranishi families of $Y$ and $Z$. By \cite[12.6.2]{KollarMoriFlips}, $\alpha$ deforms to an SQM $\tilde{\alpha} \colon \cY \dashrightarrow \cZ$ over an {isomorphism} $\alpha^* \colon \Def(Y) \rightarrow \Def(Z)$, and therefore induces an SQM between a very general deformation of $Y$ and a very general deformation of $Z$. This gives the result.
\end{proof}
	\subsection{Deformation of divisorial contractions}
	We have seen in Lemma~\ref{lem:ConeInclusions} that the movable cone of a smooth Calabi--Yau threefold can only become larger under deformation. We now determine along which faces $F$ this jumping occurs. This turns out to be closely related to the deformation theory of the corresponding exceptional divisor $E_F$. \\
	We start by observing that the birational contraction morphism $p_F \colon Y \rightarrow \bar{Y}$ associated to a face $F$ of $\Nef(Y) \cap \Bbig(Y)$ always deforms: since $h^1(\cO_Y)=h^2(\cO_Y)=0$, we see using cohomology and base change that there is a unique lift of $L$ to a base-point-free line bundle $\cL$ on $\cY$. We obtain a birational morphism 
	$P_F \colon \cY \rightarrow \bar{\cY}$ over $\Def(Y)$ which restricts to $p_F$ on $Y$. A curve $C \subset \cY$ is contracted by $P_F$ if and only if $[C] \in F^\perp$. \\
	We will say that $E \subset Y$ disappears under deformation if $E$ does not extend to an effective divisor $\cE \subset \cY$, flat over $\Def(Y)$.
	We now give a summary of the results by Wilson~\cite{WilsonCone}, who classified the exceptional divisors of primitive contractions of $Y$, and investigated their behaviour under deformation.
	\begin{proposition}\cite{WilsonCone,WilsonErratum}; see also
		\cite[Proposition~3.1]{WilsonSymplectic}. \label{prop:movablechanges}
		Let $Y$ be a smooth Calabi--Yau threefold, and let $p \colon Y \rightarrow \bar{Y}$ be a primitive divisorial contraction with exceptional divisor $E$. 
		\begin{enumerate}
			\item We have that either
			\begin{itemize}
			\item $p$ is of type {\rm II}, i.e. $p$ contracts $E$ to a point, and $E$ is a generalized del Pezzo surface. 
			\item $p$ is of type {\rm III}, i.e. $p$ contracts $E$ to a curve $C \subset \bar{Y}$ of singularities. $C$ is automatically smooth, and $p|_E \colon E \rightarrow C$ exhibits $E$ as a conic fibration surface. Moreover, the classes of components of fibers of $p|_E$ are all proportional in $H_2(Y, \ZZ)$.
			\end{itemize}
		\item $E$ disappears under a very general deformation if and only if $E$ is a conic fibration surface over a curve of positive genus. 
		\item Suppose that $E$ is a conic fibration surface with fiber class $\ell_E$. Then no curve $Z$ with $[Z]$ proportional to $\ell_E$ deforms to $Y^{\mathrm{gen}}$ if and only if $E$ is quasi-ruled over a curve of genus $1$. 
		\end{enumerate}
	\end{proposition}
For the definition of a quasi-ruled surface see Definition~\ref{def:qruled}.
Keeping the notation of Proposition~\ref{prop:movablechanges}, we note that if $E \rightarrow C$ is a conic fibration surface with general fiber $\PP^1$, then $\bar{Y}$ has transverse $A_1$-singularities at a general point of $C$. If the general fiber is the union of two $\PP^1$'s, the singularity at a general point of $C$ is transverse $A_2$. However, $C$ might have more complicated singularities at points corresponding to singular fibers of $E$, often called ``dissident" points.\\
Let now $F$ be a codimension one face of $\overline{\Mov}(Y) \cap \Bbig(Y)$ with supporting hyperplane $H$. If $H$ supports a codimension one face $F'$ of $\overline{\Mov}(Y^{\mathrm{gen}})$ 
we say that $F$ persists under deformation. Otherwise, we say that $F$ disappears under deformation.
	\begin{lemma}\label{prop:deformationcriterion}
		Let $\pi \colon \cY \rightarrow S$ be a deformation of a Calabi--Yau threefold $Y$ over a smooth germ, and let $F$ be a codimension one face of ${\Nef}(Y) \cap \Bbig(Y)$ such that the associated contraction morphism $p_F \colon Y \rightarrow \bar{Y}$ is divisorial, with exceptional divisor $E_F$. Let $L$ be a line bundle giving the contraction.
		The following are equivalent:
		\begin{enumerate}
			\item The extension of $p_F$ to a birational contraction $P_F \colon \cY \rightarrow \bar{\cY}$ is divisorial.
			\item The corresponding facet of $\overline{\Mov}(Y)$ persists under deformation in the family $\cY \rightarrow S$.
			\item $E_F$ deforms to $\cY$.
		\end{enumerate}
	\end{lemma}
	\begin{proof}
		$1) \iff 2)$. We have that $L \in F^\circ \subset \partial \overline{\Mov}(Y)$, so the corresponding facet of $\overline{\Mov}(Y)$ persists if and only if $L \in \partial \overline{\Mov}(Y^{\mathrm{gen}})$. By \cite[Lemma 2.3]{KawamataCrepant}, this happens if and only if the restriction of $P_F$ to a general fiber is divisorial, which happens if and only if $P_F \colon \cY \rightarrow \bar{\cY}$ is a divisorial contraction.\\
		$3) \implies 1)$.
		Suppose first that $E_F \subset Y$ deforms, that is, there is a flat family $\cE \rightarrow S$ with $\cE_0 \cong E_F$. We claim that $P_F$ contracts $\cE$. Let $\cL$ be the lift of $L$ to $\cY$, write $\cL=P_F^*\cH$ for relatively ample $\cH$, and write $H$ for the restriction of $\cH$ to the central fiber. Since the morphism $p_F$ on central fibers contracts $E_F$, we have $L^2 \cdot E_F=H^2 \cdot p_F(E_F)= 0$. By continuity of intersection numbers, we have $\cH^2 \cdot P_F(\cE)=\cL^2 \cdot \cE=L^2 \cdot E_F = 0$ as well, so that $\cE$ is contracted.\\
		$1) \implies 3)$
		Suppose now that $E_F \subset Y$ does not deform. By Proposition~\ref{prop:movablechanges}, $E_F$ is a conic fibration surface over a curve $C$. Let $\ell$ be the class of a fiber of $E_F$, we have $E_F \cdot \ell=-2$.  
		Suppose for a contradiction that $P_F \colon \cY \rightarrow \bar{\cY}$ is divisorial, and denote $\cE$ the exceptional divisor. $L$ is nef on every fiber of $\pi$ and big on the central fiber. Possibly shrinking $S$, we may assume $L$ is big on every fiber and so the restriction of $\cE \rightarrow S$ to every fiber is a Cartier divisor, and the induced morphism $\cE \rightarrow S$ is flat. $\cE_0$ is supported on $E_F$, so that $[\cE_0]=k[E_F]$.
		The restriction of $P_F$ to a very general fiber must also be a contraction of type III, and we denote its fiber class by $\ell^{gen}$. We have $\cE \cdot \ell^{gen}=-2$ and we may write $\ell^{gen}=k'\ell$. This shows that $kk'=1$, so that $[\cE_0]=[E_F]$, contradicting our assumption that $E_F$ doesn't deform.
	\end{proof}
	\begin{corollary}\label{cor:conicfibrationsmall}
Keeping the notation of Lemma~\ref{prop:deformationcriterion}, suppose that $E_F$ is a conic fibration surface over a curve of positive genus, and that $\cY \rightarrow S$ is a very general deformation (for example, $S=\Def(Y)$). Then $P_F$ is a small contraction. 
	\end{corollary}
	\begin{proof}
		Since $E_F$ is a conic fibration surface over a curve of positive genus, Proposition~\ref{prop:movablechanges} shows that $E_F$ disappears under deformation. By Lemma~\ref{prop:deformationcriterion}, the birational contraction morphism $P_F \colon \cY \rightarrow \bar{\cY}$ must be small.
		\end{proof}
	\subsection{Proof of Theorem~\ref{thm:main1}}
	Before our next definition, we establish a few basic properties of SQM-conic fibration surfaces over curves of positive genus.
Let $E \subset Y$ be an SQM-conic fibration surface over a curve of positive genus, and let 
\[
(\alpha \colon Y \dashrightarrow Z, p \colon E_Z \rightarrow C)
\]
be an SQM-model of $E$. 
Let $\ell_{Z} \in H_2(Z, \ZZ)$ be the class of a fiber of $p$, and let $\ell_{E}={\alpha_*}^{-1}\ell_{Z}$
be the conic fiber class of $E$. The following lemma shows that $\ell_E$ and $C$ are independent of the chosen SQM-model.
\begin{lemma}\label{lem:propertiesofSQMconicfibration}
Suppose that $E \subset Y$ is an SQM-conic fibration surface with an SQM-model 
\[
(\alpha \colon Y \dashrightarrow Z, \;p \colon E_Z \rightarrow C)
\]
with $g(C)>0$. Then every SQM-model has base curve isomorphic to $C$, and the fiber class $\ell_E$ is intrinsic to $E$. In particular, $\ell_E$ does not depend on the choice of SQM-model. Moreover:
 	\begin{itemize}
		\item $E \cdot \ell_E=-2$.
		\item The class $[E] \in H^2(Y, \ZZ)$ has a unique effective representative.
	\end{itemize}
	\end{lemma}
	\begin{proof}
	Let 
	\[
	(\alpha_i \colon Y \dashrightarrow Z_i, p_i \colon E_{Z_i} \rightarrow C_i), \quad i=1,2
	\]
	be SQM-models for $E$, and suppose that $g(C_1)>0$. Let $\ell_{Z_i} \in H_2(Z_i, \ZZ)$ be the fiber class of $p_i$, and let 
	Let $\Phi=\alpha_1 \circ \alpha_2^{-1}$, and let $V_i \subset Z_i$ be the common isomorphism loci:
	\[
	\Phi|_{V_2} \colon V_2 \rightarrow V_1
	\]
	and let $U_i=E_{Z_i} \cap V_i$. Since $\Phi$ is an SQM, $Z_i \setminus V_i$ is a union of rational curves, hence $E_{Z_i} \setminus U_i$ is a union of points and rational curves.
	Let $\phi \colon E_{Z_2} \dashrightarrow E_{Z_1}$ be the induced birational map of surfaces. Since $g(C_1)>0$, we must have that $E_{Z_1} \setminus U_1$ is disjoint from a general fiber of $p_1$. Let $B$ be a general fiber of $E_{Z_2}$. By definition, $B$ is either a smooth rational curve or a reducible conic. Let $B_i$ be an irreducible component of $B$. Its strict transform is a general rational curve $L_i$ on $E_{Z_1}$, so that $L_i \subset U_{1}$. Consequently, $\phi^{-1}(L_i) \subset U_2$. Since $\phi^{-1}(L_i)$ is proper and agrees with $B_i$ on an open subset, it follows that $\phi^{-1}(L_i)=B_i$ and hence $B_i \subset U_2$. Therefore, if $B=B_1 \cup B_2$, then $L_1$ and $L_2$ meet and are contained in the same $p_1$-fiber. We thus obtain a dominant rational map $\bar{\phi} \colon C_2 \dashrightarrow C_1$. In particular, $g(C_2)>0$ as well. Reversing the roles of $E_{Z_1}$ and $E_{Z_2}$ we see that $\bar{\phi}$ has an inverse on an open subset. It follows that $\bar{\phi}$ is birational, and therefore an isomorphism since the $C_i$ are smooth. By taking $\alpha_1=\alpha_2$, this shows in particular that the fiber class $\ell_Z$ of $p \colon E_Z \rightarrow C$ is uniquely determined by $E_Z$. 
		Since $\alpha_1\circ \alpha_2^{-1}\colon Z_2 \dashrightarrow Z_1$ maps $\ell_{Z_2}$ to $\ell_{Z_1}$, the 
		conic fiber class $\ell_E$ of $E$ is independent of the choice of SQM. \\
		Next, 
%
		we have that $E \cdot \ell_E=E_Z \cdot \ell_Z$, so it is enough to prove $E_Z \cdot \ell_Z=-2$.
		Recall that a general fiber $B$ of $p \colon E_Z \rightarrow C$ is either $\PP^1$ or a union of two $\PP^1$'s meeting in a node, in particular, $\omega_B$ has degree $-2$. $E_Z$ is Cohen--Macaulay, and $B$ is an effective Cartier divisor on $E_Z$, so the adjunction formula \cite[5.73]{KollarMori} shows that 
		\begin{align*}
			\omega_{E_Z}& \cong \cO_{E_Z}(K_Z+E_Z)\\
			\omega_B&\cong \omega_{E_Z}(B)|_B
		\end{align*}
		We therefore have  
		$
		\deg \omega_{E_Z}|_B=K_Z \cdot B+E_Z \cdot B=E_Z \cdot B$ since $K_Z=0$.
		Since $\deg \cO_B(B)=0$, we conclude that 
		\[
		E_Z \cdot \ell_Z=E_Z \cdot B=\deg(\omega_B)=-2.
		\]	
		For the last part, suppose that $F$ is an effective divisor with $[E]=[F]$. Then $F \cdot \ell_E=-2$, so every curve of class $\ell_E$ must be contained in $F$. Since such curves sweep out $E$, $E$ must be an irreducible component of $F$. We may therefore write 
		\[
		F=mE+D
		\]
		for some effective divisor $D$. Intersecting with $H^2$ for an ample divisor $H$
		and using $[F]=[E]$, we see that $m=1$ and $D=0$.
	\end{proof}
	We make the following definition:
	\begin{definition}\label{def:roots}
	Let $Y$ be a Calabi--Yau threefold. 
	\begin{itemize}
		\item We define $\Delta_Y^{big}$ as the set of classes $\alpha \in H^2(Y, \ZZ)$ such that $\alpha=[E]$ for some SQM-conic fibration surface $E$ over a curve of positive genus. 
		\item We define $\Delta_Y^{sm} \subset \Delta_Y^{big}$ as the subset consisting of classes of quasi-ruled surfaces over a curve of genus $1$.
		\end{itemize}
		Writing $\ell_E$ for the fiber class of $E$, we define 
		\begin{itemize}
		\item 
		$\bar{\cC}^{big}$ as the intersection of the closed half-spaces 
	\[
	\left\{x \in H^2(Y, \RR) \mid x \cdot \ell_E \geq 0\right\}, \quad  [E] \in \Delta_Y^{big}. 
	\]
	\item $\bar{\cC}^{sm}$ as the intersection of the closed half-spaces 
		\[
	\left\{x \in H^2(Y, \RR) \mid x \cdot \ell_E \geq 0\right\}, \quad [E] \in \Delta_Y^{sm}. 
	\]
	\end{itemize}
	\end{definition}
	In view of Lemma~\ref{lem:propertiesofSQMconicfibration}, we will sometimes abuse notation and write $E \in \Delta_Y^{big}$.
We will see in Section~\ref{sec:roots} that the pairs $([E], \ell_E)$ can be interpreted as simple roots in a generalized sense. 
	We have:
	\begin{proposition}\label{prop:movableisintersection}
		We have that
		\begin{align}
				\Nef(Y)=\bar{\cC}^{sm} \cap \Nef(Y^{\mathrm{gen}}).	\label{eq:nef}\\
		\overline{\Mov}(Y)=\bar{\cC}^{big} \cap  \overline{\Mov}(Y^{\mathrm{gen}}). \label{eq:mov}
		\end{align}
	\end{proposition}
	\begin{proof}
		The first equality	is proved in~\cite[Lemma 1.1]{WilsonElliptic}.\\
		We now prove the second equality.
		Since curves of class $\ell_{E_i}$ sweep out the divisor $E_i$, any movable divisor $D$ satisfies $D \cdot \ell_{E_i} \geq 0$, so 
		\[
		\overline{\Mov}(Y) \subset \bar{\cC}^{big} \cap \overline{\Mov}(Y^{\mathrm{gen}})
		\]
		For the reverse inclusion, suppose $D \in \bar{\cC}^{big} \cap \overline{\Mov}(Y^{\mathrm{gen}})$. Let $H$ be a generic ample divisor on $Y$, and consider the line $L: \{D+ tH \mid t \in \RR \}$ in the vector space $H^2(Y, \RR)$. $L$ must meet $\partial \overline{\Mov}(Y)$ transversely, and we set 
		\[
		t_0=\inf\{t \in \RR \mid D+tH \in \overline{\Mov}(Y) \}
		\]
		and $x=D+t_0H \in \partial  \overline{\Mov}(Y) $.
		We claim that $t_0 \leq 0$, which will show that $D \in \overline{\Mov}(Y)$ and complete the proof.\\ 
		Suppose first that $x \notin \Bbig(Y)$. By Lemma~\ref{lem:ConeInclusions}, $x \notin \Bbig(Y^{\mathrm{gen}})$, so $x$ must be a point on $\partial \overline{\Mov}(Y^{\mathrm{gen}})$. Since $D=D+0 \cdot H \in \overline{\Mov}(Y^{\mathrm{gen}})$, we have $t_0 \leq 0$ in this case.\\
		Suppose now that $x \in \Bbig(Y)$. By Proposition~\ref{prop:big-rational}, $\overline{\Mov}(Y)$ is rational polyhedral in a neighbourhood of $x$, and we let $F$ be the face of $\overline{\Mov}(Y)$ whose interior contains $x$. Since $H$ was chosen generically, we may assume that $F$ has codimension one. 
		If $F$ persists under deformation, then $x \in \partial \overline{\Mov}(Y^{\mathrm{gen}})$, and we have $t_0 \leq 0$ as before.
		If $F$ disappears under deformation, let $\alpha \colon Z \dashrightarrow Y$ be a SQM such that $\alpha^*(F)$ is a facet of $\Nef(Z)$, and let 
		$\ell_F$ be the primitive generator of the corresponding extremal ray of $\NE(Z)$, so that $\alpha^*(F)$ is contained in the hyperplane $\ell_F^\perp$. In particular we have 
	$\alpha^*(x) \cdot \ell_F=0$.
		
		The corresponding exceptional divisor $E_F \subset Z$ disappears under deformation, by Lemma~\ref{prop:deformationcriterion}, so by Proposition~\ref{prop:movablechanges}, $E_F$ is a conic fibration surface over a curve of positive genus. The strict transform $E=\alpha^*(E_F)$ is then (by definition) an SQM-conic fibration surface over a curve of positive genus and so $\alpha_*(\ell_F)=\ell_E$. We have
		\[
		0=\alpha^*(x) \cdot \ell_F=x \cdot \ell_E=(D \cdot \ell_E) +t_0(H \cdot \ell_E)
		\]
		Since $H$ is ample, $H \cdot \ell_E>0$, and $D \cdot \ell_E \geq 0$ since $D \in \bar{\cC}^{big}$. It follows that 
		\[
		t_0=-\frac{D \cdot \ell_E}{H \cdot \ell_E} \leq 0
		\]
		 in this case as well.
	\end{proof}
Having shown that the movable cone can only jump along faces corresponding to SQM-conic fibration surfaces, we now show that each SQM-conic fibration surface has a corresponding face of $\overline{\Mov}(Y) \cap \Bbig(Y)$:
	\begin{proposition}\label{prop:everysurfacecontracts}
		Let $Y$ be a Calabi--Yau threefold and let $E \in \Delta_Y^{big}$. Then there exists a sequence of flops $Y \dashrightarrow Z$, and a primitive divisorial contraction $Z \rightarrow \bar{Z}$ which contracts the strict transform $E_Z$ of $E$ to $C$ along its ruling. In particular, there is a codimension one face of $\overline{\Mov}(Y) \cap \Bbig(Y)$ corresponding to $E$.
	\end{proposition}
	\begin{proof}
		Let 
		\[
		(\alpha \colon Y \dashrightarrow Y', p \colon E_{Y'} \rightarrow C)
		\]
		be an SQM-model for $E$.
		We run an MMP for the pair $(Y', K_{Y'}+\epsilon E_{Y'})$, where $0< \epsilon  \ll 1$. (Note $(Y',\epsilon E_{Y'})$ is klt and $K_{Y'}+ \epsilon E_{Y'} \sim \epsilon E_{Y'}$.) This will terminate with a minimal model $(\bar{Z}, B)$ with $K_{\bar{Z}}+B$ nef (we can't get a Mori fiber space because $Y'$ is not uniruled). If $\Gamma$ is a rational curve with $E_{Y'} \cdot \Gamma<0$, then $\Gamma$ must be an irreducible component of a fiber of the ruling (as $C$ has positive genus), so that the MMP must be a sequence of flops $Y' \dashrightarrow Z$ followed by the divisorial contraction of $E_Z$, so that the minimal model is $(\bar{Z},B)=(\bar{Z},0)$ where $\bar{Z}$ has a curve $C$ of cDV singularities. Note that the MMP has to contract $E_{Y'}$ because it is negative on the irreducible components of the general fiber of the ruling. The contraction has relative Picard rank 1 by the cone theorem \cite[Theorem 3.7]{KollarMori}, and therefore corresponds to a codimension one face of ${\Nef}(Z) \cap \Bbig(Z)$ and $\overline{\Mov}(Z) \cap \Bbig(Z)$, and therefore a codimension one face of $\overline{\Mov}(Y) \cap \Bbig(Y)$.
	\end{proof}
If $E \in \Delta_Y^{big}$ can already be contracted via a primitive contraction (i.e. no flops are necessary in Proposition~\ref{prop:everysurfacecontracts}), then we say that $E$ is contractible on $Y$. If $E$ is quasi-ruled, then $E$ is always contractible on $Y$. In fact, it is clear from the proof of Proposition~\ref{prop:everysurfacecontracts} that $E$ is contractible on $Y$ precisely if $E$ is a conic fibration surface and the subspace of $H_2(Y, \QQ)$  spanned by irreducible components of fibers of $p \colon E \rightarrow C$ is $1$-dimensional.
We now give the proof of Theorem~\ref{thm:main1}.
	\begin{proof}[Proof of Theorem~\ref{thm:main1}] \label{pf:1}
		It is immediate from Proposition~\ref{prop:movableisintersection} that $\overline{\Mov}(Y)=\overline{\Mov}(Y^{\mathrm{gen}})$ if $Y$ does not contain any such surface $E$. If none of the threefolds $Y_b$ contains such a surface, we have $\overline{\Mov}(Y_b)=\overline{\Mov}(Y^{\mathrm{gen}})$ for all $b \in B$, so that the movable cone is constant in the family. Suppose on the other hand that one threefold $Y_b$  in the family contains at least one such surface $E$. Suppose without loss of generality that $Y_b=Y$. By Proposition~\ref{prop:everysurfacecontracts}, there is a SQM $\alpha \colon Y \dashrightarrow Z$ such that the strict transform $E_Z$ is the exceptional divisor of a primitive contraction $p \colon Z \rightarrow \bar{Z}$. Let $F$ be the corresponding facet of $\Nef(Z) \cap \Bbig(Z)$, so that $p=p_F$. Since $E_Z$ is a conic fibration surface over a curve of positive genus, Corollary~\ref{cor:conicfibrationsmall} shows that the induced contraction morphism $P_F \colon \cZ \rightarrow \bar{\cZ}$ is small.
		So if $D \in F^\circ \subset \partial \overline{\Mov}(Z)$ is a divisor giving the contraction, then $D$ is in the interior of $\overline{\Mov}(Z^{gen})$ so that $\overline{\Mov}(Z)$ jumps under deformation. Using Lemma~\ref{lem:preservegen}, we conclude that $\overline{\Mov}(Y)$ also jumps under deformation, i.e. $\overline{\Mov}(Y) \subsetneq \overline{\Mov}(Y^{\mathrm{gen}})$. It follows that the movable cone is not constant in fibers of the family.
	\end{proof}
\section{Roots and the Wilson Weyl group}\label{sec:roots}
\subsection{Generalized Simple Roots}
\begin{definition}\label{def:genroots}
Let $V$ be a vector space, and $V^*$ its dual. 
A subset $\Delta \subset V$ together with an injection $\vee \colon \Delta \rightarrow V^*$ is called a set of generalized roots if 
\[
 \alpha^\vee(\alpha)=-2
 \]
for all $\alpha \in \Delta$, and for all distinct $\alpha, \beta \in \Delta$, we have
\begin{enumerate}
	\item $\alpha^\vee(\beta) \in \ZZ_{\geq 0}$ 
	\item $\alpha^\vee(\beta) =0$ if and only if $\beta^\vee(\alpha)=0$
	\item $\{\alpha, \beta\}$ (respectively $\{\alpha^\vee, \beta^\vee\}$) is a linearly independent set in $V$ (respectively in $V^*$) 
\end{enumerate}
\end{definition}
\begin{remark}
	Our definition is more general than the one usually considered in the literature \cite{Bourbaki, Looijenga}. We do not require that the set of generalized roots is linearly independent, nor do we require that $\Delta$ is finite. It is for this reason that we avoid calling $\Delta$ a ``root basis". Our definition is closest in spirit to the ``root data" considered in \cite[Section 2]{MoodyPianzola}. There, the set of generalized roots is also not assumed to be finite or linearly independent, but has to satisfy an independence condition which is stronger than (3) above. This stronger condition is necessary to show that the fundamental Weyl chamber is nonempty. With our weaker condition (3), this might well fail in general (see Example~\ref{ex:rootsfail}). However, in all geometric examples we consider, it will be obvious that the fundamental Weyl chamber is in fact nonempty, see Lemma~\ref{lem:nonempty}.
\end{remark}
If the injection $\vee$ is clear from context, we simply say that $(V, \Delta)$ is a set of generalized roots. We will often find it convenient to index the elements of $\Delta$.
Our task is to show that in this setup, the classical results of \cite[Chap.~V, §4, nos.~4--6]{Bourbaki} on Weyl groups and fundamental chambers continue to hold. 
For each $\alpha, \beta\in \Delta$, define $m(\alpha, \beta) \in \ZZ \cup \infty$ by
\begin{equation}\label{eq:defnofm}
m(\alpha, \beta)=
\begin{cases}
	1 \quad \text{if} \quad \alpha=\beta\\
	2,3,4,6\quad  \text{if} \quad  \alpha^\vee(\beta)\cdot\beta^\vee(\alpha) =0,1,2,3\\
	\infty \quad  \text{otherwise} 
\end{cases}
\end{equation}
We remark that $M=(m(\alpha, \beta))_{\alpha, \beta\in \Delta}$ is a Coxeter matrix in the sense of \cite[Chap.~VI, §1, no.~9]{Bourbaki}. Let $ W(\Delta)$ be the group defined by a family of generators $(g_\alpha)_{\alpha\in \Delta}$ and relations $(g_\alpha g_{\beta})^{m(\alpha, \beta)}=1$ if $m(\alpha, \beta) \neq \infty$. 
We call $ W(\Delta)$ the Weyl group associated to the set of generalized roots.

There is a natural representation of the Weyl group of a set of generalized roots as a reflection group:
\begin{definition}
	For each $\alpha\in \Delta$ we define the reflection
	\[
	\sigma_{\alpha} \colon V \rightarrow V, \quad x \mapsto x+\alpha^\vee(x) \alpha
	\]
\end{definition}
\begin{proposition}\label{prop:WeylgroupactsonH2}
The assignment $g_\alpha \mapsto \sigma_\alpha$ defines a representation of $W(\Delta)$.
\end{proposition}
 \begin{proof}
 	It is enough to show that $(\sigma_\alpha\sigma_{\beta})^{m(\alpha, \beta)}=1$. If $\alpha=\beta$ this follows from the fact that $\sigma_\alpha$ has order $2$. If $\alpha^\vee(\beta)\,\beta^\vee(\alpha)=0$ then in fact $\alpha^\vee(\beta)=\beta^\vee(\alpha) =0$ (by our definition of generalized roots). A calculation shows that $\sigma_\alpha\sigma_{\beta}$ has order $2$ in this case. For the other cases, note that since the sets $\{\alpha, \beta\} $ and $ \{\alpha^\vee, \beta^\vee\}$ are linearly independent, they form a root basis in $V$ in the sense of Looijenga~\cite{Looijenga}. By \cite[1.2]{Looijenga}, we have that 
$\sigma_\alpha\sigma_{\beta}$ has order $3,4,6$, depending on whether $\alpha^\vee(\beta) \cdot \beta^\vee(\alpha)=1,2,3$. 
\end{proof}
This argument in fact shows that $(W(\Delta), \Delta)$ is a Coxeter system, in the sense of \cite[IV.1.3 Definition 3]{Bourbaki}.
Given $w \in W$, we denote $w(x)$ the transform of $x$ by $\sigma(w)$. 
For each $\alpha\in \Delta$, the inequality $\langle x, \alpha^\vee \rangle>0 $ defines an open half-space $A_\alpha$. It follows from the definition of $\sigma_\alpha$ that $\sigma_\alpha(A_\alpha) \cap A_\alpha =\emptyset$:
We define the open cone (fundamental chamber)
\[
\cC_\Delta=\bigcap_{\alpha\in \Delta}A_\alpha \subset V
\]
In our generalized setup, it can easily happen that the cone $\cC_\Delta$ is empty:
\begin{example}\label{ex:rootsfail}
Let $V=\RR^2$, and define
\[
\alpha_1=\begin{pmatrix}1\\0 \end{pmatrix}, \alpha_2=\begin{pmatrix}0\\1 \end{pmatrix}, \alpha_3=\begin{pmatrix}-1\\-1 \end{pmatrix}, \qquad \alpha_1^\vee=\begin{pmatrix}-2\\1 \end{pmatrix}^T, \alpha_2^\vee=\begin{pmatrix}1\\-2 \end{pmatrix}^T, \alpha_3^\vee=\begin{pmatrix}1\\1 \end{pmatrix}^T.
\]
The sets $\Delta=\{\alpha_1, \alpha_2, \alpha_3\}$ and $\{\alpha_1^\vee, \alpha_2^\vee, \alpha_3^\vee\}$ satisfy the properties of Definition~\ref{def:genroots}, but the cone $\cC_\Delta$ is empty. 
\end{example}
	 In the geometric setups we consider, it will however be clear that $\cC_\Delta$ is nonempty, and even full-dimensional, see Lemma~\ref{lem:nonempty}. We have the following result: 
\begin{lemma}\label{thm:weylchambersdontmeet}
	Let $(V, \Delta)$ be a set of generalized roots, and suppose that the cone $\cC_\Delta$ is nonempty. 
	Then $w\cC_\Delta \cap \cC_\Delta \neq \emptyset$ implies that $w=1$. 
\end{lemma}
\begin{proof}
	The proof of \cite[V.4.4, Theorem 1]{Bourbaki} only uses abstract properties of Coxeter systems, and goes through mostly unchanged. Note that the condition that $\cC_\Delta \neq \emptyset$ is necessary for the base case of the induction in the proof of $P_n$ and $Q_n$. The only necessary change is the proof of Chap. {\rm V}, §5, Lemma 1, for $\Delta=\{\alpha, \beta\}$ a set of two generalized simple roots. However, such a set is a generalized root basis in the sense of Looijenga, as any two distinct generalized simple roots are linearly independent. In this setting, the lemma in question is proved in \cite[(1.2)] {Looijenga}.
\end{proof}
The main result of this section is the following: 
\begin{theorem}\label{thm:Weylgroupactionisfaithful}
	Let $(V, \Delta)$ be a set of generalized roots, and suppose that the cone $\cC_\Delta$ is nonempty. 
	The action of $W= W(\Delta)$ on $V$ is faithful. In particular, the only relations in the subgroup of $\GL(V)$ generated by $(\sigma_\alpha)_{\alpha\in \Delta}$ are $\sigma_\alpha^2=1$ and if $\langle \alpha, \beta^\vee \rangle \langle \beta, \alpha^\vee \rangle=0,1,2,3$, then $\sigma_\alpha\sigma_\beta$ is of order $2,3,4,6$ respectively. 
\end{theorem}
\begin{proof}
Suppose that $w \in W$ acts trivially on $V$. Since $\cC_\Delta$ is nonempty, the set $w\cC_\Delta \cap \cC_\Delta=\cC_\Delta$ is also nonempty. By Lemma~\ref{thm:weylchambersdontmeet}, it follows that $w=1$. 
\end{proof}
We define the Tits cone	\[
	U=\bigcup_{w \in W} w(\bar{\cC}_\Delta)
	\]	
	Given any root $\alpha\in \Delta$, denote $H_\alpha$ the hyperplane in $V$ defined by $\alpha^\vee(x)=0$. We construct a partition (stratification) of $\bar{\cC}_\Delta$ as follows. 
	Given any subset $X \subset \Delta$, we define the cone 
		\[
	C_X=\left( \bigcap_{\alpha \in X} H_\alpha \right) \cap \left( \bigcap_{\alpha\in \Delta \setminus X} A_\alpha \right)
	\]
	Note that in our generalized setup, the various $C_X$ for $X \subset \Delta$ do not necessarily give a partition of $\bar{\cC}_\Delta$. For example, it is possible that $C_X=C_{X'}$ for distinct subsets $X$ and $X'$. 
	To get around this issue, we do the following:
	given $x \in \bar{\cC}_\Delta$, let $X \subset \Delta$ be the set of $\alpha\in \Delta$ such that $x \in H_\alpha$. As $x$ varies over $\bar{\cC}_\Delta$, we obtain in this way a collection of subsets of $\Delta$, which we denote $B(\Delta)$.
	By construction, the cones $C_X$, for $X \in B(\Delta)$ then give a partition of $\bar{\cC}_\Delta$. 
	Given $X \in B(\Delta)$, we denote $W_X$ the subgroup of $W$ generated by $X$. Then clearly $w(x)=x$ if $w \in W_X$ and $x \in C_X$. 
	The same proof as \cite[V. 4.6 Proposition 6]{Bourbaki} shows:
\begin{corollary}\label{cor:fundamentaldomain}
Suppose that $\cC_\Delta$ is nonempty. Then $\bar{\cC}_\Delta$ is a fundamental domain for the action of $ W(\Delta)$ on $U$ in the following strong sense: if $w\in W(\Delta)$ and $X \in B(\Delta)$ are such that $w(C_X) \cap \bar{\cC}_\Delta \neq \emptyset$, then $w \in W_X$ and (hence) $w$ leaves $C_X$ pointwise fixed.
\end{corollary}
In view of this, we call $\bar{\cC}_\Delta$ the fundamental chamber. 
\subsection{Generalized roots on Calabi--Yau threefolds}
We now show that the set $\Delta_Y^{big}$ (see Definition~\ref{def:roots}) forms a set of generalized roots in $H^2(Y, \RR)$. 
\begin{proposition}\label{prop:Eirootsystem}
	Let $Y$ be a Calabi--Yau threefold. Then $(H^2(Y, \RR), \Delta_Y^{big})$ is a set of generalized roots, with the injection 
	\[
	\vee \colon \Delta_Y^{big} \rightarrow H^2(Y, \RR)^* \cong H_2(Y, \RR)
	\]
	given by sending an SQM-conic fibration surface $E$ to its conic fiber class $\ell_E$. 
\end{proposition}
\begin{proof}
We have seen in Lemma~\ref{lem:propertiesofSQMconicfibration} that $E$ uniquely determines the fiber class $\ell_E \in H_2(Y, \ZZ)$ so that $\vee$ is well-defined. 
We have also seen that 
\[
E \cdot \ell_E=-2.
\]
	Let now $E_1,  E_2 \in \Delta_Y^{big}$ be distinct, and denote the conic fiber classes $\ell_1$ and $\ell_2$. Since the $E_i$ and $\ell_i$ are integral classes, we have that $E_i \cdot \ell_j \in \ZZ$. The intersection $E_1 \cap E_2$ has dimension at most $1$, and the class $\ell_1$ sweeps out $E_1$. Therefore $\ell_1$ can be represented by a curve with support not entirely contained in $E_2$, so that $E_2 \cdot \ell_1 \geq 0$. This proves (1).  
	This also proves that $\vee$ is indeed an injection: if $\ell_{1}=\ell_2$, then we have
	\[
	-2=E_1 \cdot \ell_1=E_1 \cdot \ell_2 \geq 0,
	\]
	contradiction.\\
	If $[E_1]=\lambda[E_2]$ in $H^2(Y, \RR)$ for some $\lambda \in \RR$, intersecting with $\ell_1$ shows that $\lambda$ must be negative, which is a contradiction, by intersecting with an ample class. The same argument shows that we cannot have $\ell_1=\lambda \ell_2$ in $H_2(Y, \RR)$, so $\{[E_1], [E_2]\}$ and $\{\ell_1, \ell_2\}$ are linearly independent sets. This proves (3). \\
	Finally, suppose for a contradiction that $E_1 \cdot \ell_2=0$ and $E_2 \cdot \ell_1 \neq 0$.
	Let 
	\[
	(\alpha_i \colon Y \dashrightarrow Z_i, p_i \colon E_{i,Z_i} \rightarrow C_i), \qquad i=1,2
	\]  
	be SQM-models for $E_1$ and $E_2$. By definition of $\Delta_Y^{big}$, we have $g(C_i)>0$.
	We write $E_{i, Z_j}$ for the strict transform of $E_i$ on $Z_j$ and $\ell_{i, Z_j}$ for the image of $\ell_i$ under the induced identification of curve classes. \\
	Taking strict transforms on $Z_1$, $E_{1, Z_1}$ is a conic fibration surface over $C_1$, with fiber class $\ell_{1, Z_1}$. We have 
	\[
	E_2 \cdot \ell_1=E_{2, Z_1} \cdot \ell_{1, Z_1}=E_{2, Z_1}|_{E_{1, Z_1}} \cdot \ell_{1, Z_1}>0,
	\]
	so there must be a component $B$ of the intersection $\Gamma_1=E_{1, Z_1} \cap E_{2, Z_1}$ horizontal over $C_1$ which is therefore not rational.\\
	Taking strict transforms on $Z_2$, we similarly have:
	\[
	E_1 \cdot \ell_2=E_{1, Z_2} \cdot \ell_{2, Z_2}=E_{1, Z_2}|_{E_{2, Z_2}} \cdot \ell_{2, Z_2}=0
	\]
	so every component of the intersection $\Gamma_2=E_{1, Z_2} \cap E_{2, Z_2}$ is vertical over $C_2$. In particular every component of $\Gamma_2$ is rational. 
	This is a contradiction: the exceptional locus of the composite SQM $\alpha_2 \circ \alpha_1^{-1} \colon Z_1 \dashrightarrow Z_2$ is a composite of flops, and therefore its exceptional locus is a union of rational curves. It follows that $B$ is not contracted by $\alpha_2 \circ \alpha_1^{-1} $, so that its strict transform is a component of $\Gamma_2$. But we showed that $B$ is not rational and that every component of $\Gamma_2$ is rational, contradiction.
	Condition (2) of Definition~\ref{def:genroots} now follows by symmetry, completing the proof.
%
%
\end{proof}
To declutter notation, we will write
\begin{equation*}
	W_Y^{big}=W(\Delta_Y^{big}), \qquad
	\cC^{big}=\cC(\Delta_Y^{big}),\qquad
	U^{big}=U(\Delta_Y^{big})
\end{equation*}
and make similar definitions for $\Delta_Y^{sm}$. This agrees, of course, with our notation in Section~\ref{sec:deformation}. The cone $\bar{\cC}^{big}$ defined there is the closure of $\cC^{big}$.
As in the previous section, we denote the reflection defined by $E \in \Delta_Y^{big}$ by $\sigma_E$. To reiterate, 
\[
\sigma_E \colon H^2(Y, \ZZ) \rightarrow H^2(Y, \ZZ) \qquad D \mapsto D+ (D \cdot \ell_E) E.
\]
By Theorem~\ref{thm:Weylgroupactionisfaithful}, the action of $W_Y^{big}$ on $H^2(Y, \RR)$ is faithful, so we can identify $W_Y^{big}$ with the subgroup of $\GL(H^2(Y, \RR))$ generated by the $\sigma_{E}$ for $E \in \Delta_Y^{big}$, and we will denote the generator $g_E$ of $W_Y^{big}$ by $\sigma_E$ from now on.
\begin{lemma}\label{lem:nonempty}
The fundamental chambers $\cC^{sm}$ and $\cC^{big}$ are nonempty.
\end{lemma}
\begin{proof}
	The defining equations of $\cC^{sm}$ and $\cC^{big}$ are all of the form $\langle x, C \rangle>0$ for $C$ a curve on $Y$. Therefore both cones contain the ample cone of $Y$.
\end{proof}
\subsection{Geometric realization of Weyl groups}
Throughout this section, $S$ denotes the germ of a smooth complex analytic manifold. Correspondingly, an open subset $U \subset S$ is understood to be an open set in the analytic topology. \\
The automorphism $\sigma_{E}$ of $H^2(Y, \ZZ)$ arises geometrically in two closely related ways:  In Section~\ref{sec:monodromy}, we will see that if $E \in \Delta_Y^{sm}$, then $\sigma_E$ is realized by a monodromy transformation. In this section, we first show that $\sigma_E$ is induced by a flop of the total space $\cY$ of a very general deformation of $Y$, under appropriate identifications of cohomology groups. 
We introduce the following convenient notation: 
let
\[
\Phi \colon \mathcal Y \dashrightarrow \mathcal Y'
\]
be a birational map between two deformation families of $Y$. We denote by
\[
\Phi^0 \colon H^2(Y,\mathbb R) \rightarrow H^2(Y,\mathbb R)
\]
the composite isomorphism
\[
H^2(Y,\mathbb R)
\cong 
H^2(\mathcal Y,\mathbb R)
\xrightarrow{\Phi_*}
H^2(\mathcal Y',\mathbb R)
\cong
H^2(Y',\mathbb R)
\]
obtained by parallel transport and 
strict transform $\Phi_*$ of divisor classes.\\
We recall the case of surfaces~\cite{Burns, Reid2}. Let $X$ be a smooth $K$-trivial surface containing a $(-2)$-curve $B$, and consider a $1$-parameter deformation $\cX$ in which $B$ does not deform. Let 
\[
\Phi_B \colon \cX \dashrightarrow \cX'
\]
be the flop of $B$. By~\cite[Theorem 6.3]{Reid2}, the induced birational map on central fibers extends to an isomorphism $\bar{\phi}_B \colon X \rightarrow X'$ and
\[
\Phi^0_B=\sigma_B,
\]
where $\sigma_B$ is the reflection in the $(-2)$-curve $B$, and we identify $H^2(X)$ and $H^2(X')$ using the isomorphism $\bar{\phi}_B$.
Of course, it is also well-known that $\sigma_B$ is realized by monodromy, it is the Picard-Lefschetz transformation associated to a degeneration of $X$ to an ordinary double point.\\
We recall that $E \in \Delta_Y^{big}$ is {contractible} if there is a primitive contraction morphism $p_E \colon Y \rightarrow \bar{Y}$ contracting $E$ along its ruling. We first show that we may flop $E$ in the total space of a very general deformation: 
\begin{proposition}\label{prop:flopofroot}
Let $\pi \colon \cY \rightarrow S$ be a very general deformation family of $Y$ over a smooth germ $(0 \in S)$. Suppose that $E \in \Delta_Y^{big}$ is contractible, and let $p_E \colon Y \rightarrow \bar{Y}$ be the associated primitive contraction. Denote $P_E \colon \cY \rightarrow \bar{\cY}$ its extension to a contraction over $S$. 
\begin{enumerate}
	\item 
$P_E \colon \cY \rightarrow \bar{\cY}$ is a flopping contraction, and the flop $P_E^+ \colon \cY^+ \rightarrow \bar{\cY}$ over $S$ exists.
\item 
Denoting 
\[
\Phi_E \colon \cY \dashrightarrow \cY^+
\]
the flop, the induced rational map 
\[
\phi_E \colon Y \dashrightarrow Y^+
\]
on central fibers extends to an isomorphism, denoted $\bar{\phi}_E \colon Y \xrightarrow{\sim} Y^+$.
In particular, $\cY^+ \rightarrow S$ is again a deformation of $Y$ (possibly after shrinking $S$).
\item  The induced rational map $\cY_t \dashrightarrow \cY_t^+$ on a very general fiber is biregular if $E \in \Delta_Y^{sm}$.
\end{enumerate}
\end{proposition}
We emphasize that we are not claiming that $\Phi_E \colon \cY \dashrightarrow \cY^+$ is defined everywhere on $Y$. Only after restricting the birational map of total spaces to central fibers does it extend to an isomorphism. 
 \begin{proof}
 	\begin{enumerate}
 		\item 
	Let $H \in \ell_E^\perp$ be a big and nef line bundle defining the contraction of $E$. After shrinking $S$, we may assume that $\Pic(\cY) \cong \Pic(Y)$ and we denote $\cH$ and $\cE$ the line bundles corresponding to $H$ and $E$. After shrinking $S$ further, we may also assume that $p_E \colon Y \rightarrow \bar{Y}$ extends to a contraction morphism $P_E \colon \cY \rightarrow \bar{\cY}$ over $S$ defined by $\cH \in \Pic(\cY)$.
Since $E \in \Delta_Y^{big}$ and $\cY \rightarrow S$ is a very general deformation, $P_E$ is a small contraction by Corollary~\ref{cor:conicfibrationsmall}. 
 	To show that the flop exists, we show that $P_E$ is a flipping contraction for a suitable klt pair. \\
 	Since the contraction $P_E$ defined by $\cH$ is birational on a very general fiber, $\cH$ is relatively big,
 	and therefore 
 	\[
 	\cL=\cH+\epsilon \cE \in \Bbig(\cY/S)
 	\]
 	as well for $\epsilon$ a small positive rational number. By the relative Kodaira lemma, we may write 
 	\[
 	\cL \sim_{\RR, S} \mathcal{A}+\mathcal{D}
 	\]
 	where $\mathcal{A}$ is relatively ample and $\mathcal{D}$ is effective. Shrinking $S$ and replacing $\mathcal{A}$ with an $\RR$-linearly equivalent effective divisor, it follows that there is $\Delta \geq 0$ which is $\RR$-linearly equivalent to $\cL$. Since $\cY/S$ is smooth, the pair $(\cY, \delta \Delta)$ is klt for $0 \leq \delta \ll 1$. 
 We have 
 	\[
 	\Delta \cdot \ell_E=(\cH+\epsilon \cE)\cdot \ell_E =0+\epsilon E \cdot \ell_E<0
 	\]
 	Since $\rho(\cY/\bar{\cY})=1$ and $K_{\cY/S}=0$, it follows that 
 	\[
 	-(K_{\cY/S}+\delta \Delta)
 	\]
 	is $P_E$-ample, so that $P_E \colon \cY \rightarrow \bar{\cY}$ is a flopping contraction for the pair $(\cY, \delta\Delta)$. By \cite[Theorem 17.9]{Fujino}, the flop $P_E^+ \colon \cY^+ \rightarrow \bar{\cY}$ exists. 
 	\item 
 	After restricting $S$ to a smooth very general curve germ $(0 \in T)$, the restriction
 	\[
 	{\Phi_E}|_T \colon \cY|_T \dashrightarrow {\cY}^+|_T
 	\]
 	is the flop of 
 	\[
 	P_E|_T \colon \cY|_T \rightarrow \bar{\cY}|_T
 	\]
 	by uniqueness of flips.
 	The central fiber $Y^+ \subset \cY^+|_T$ is a Cartier divisor. 
 	Since the exceptional locus of ${\Phi_E}|_T$ has codimension at least two, it identifies prime divisors on $\cY_T$ and $\cY_T^+$, so that $Y^+$ is the strict transform of $Y$. 
 	Since $Y$ is smooth, the pair $(\cY|_T, Y)$ is plt. The pair $(\cY^+|_T, Y^+)$ is then plt as well. It follows that $Y^+$ is normal by \cite[Proposition 5.51]{KollarMori}.
 	We now show that the restriction of $\Phi_E$ to central fibers, denoted
 	\[
 	\phi_E \colon Y \dashrightarrow Y^+
 	\]
 	extends to an isomorphism. 
 	We first show that $\phi_E$ is a small modification. By upper-semicontinuity, the exceptional locus of $P_E^+$ has to meet the central fiber, so $p_E^+$ is not an isomorphism.
 	The contraction $p_E \colon Y \rightarrow \bar{Y}$ is the divisorial contraction of an extremal ray associated to the klt pair $(Y, \epsilon E)$, and therefore $\bar{Y}$ is $\QQ$-factorial. Since $Y^+$ is normal, $p_E^+ \colon Y^+ \rightarrow \bar{Y}$ cannot be a small contraction, so $p_E^+$ is a divisorial contraction. Since $Y^+$ is Calabi--Yau, $p_E^+$ is crepant. Its exceptional divisor corresponds to an exceptional divisorial valuation $v$ on $\bar{Y}$. We claim that the only such valuation on $\bar{Y}$ corresponds to $E$, so that $(p_E^+)^{-1} \colon \bar{Y} \dashrightarrow Y^+$ must extract $E$ and $\phi_E$ is small. Indeed, if the center of $v$ on $Y$ is a divisor, then $v=v_E$. If the center of $v$ on $Y$ has codimension at least two, then its discrepancy $a(v, Y)>0$ since $Y$ is smooth. Since $p_E$ is crepant, we have
 	\[
 	a(v, \bar{Y})=a(v, Y),
 	\]
 	so that $v$ is not crepant in this case. \\
 	 It follows that $\phi_E \colon Y \dashrightarrow Y^+$ is defined at the generic point of $E$ and therefore a small map. By \cite[Lemma 3.38]{KollarMori}, $Y^+$ is terminal, and therefore $\phi_E$ is a composite of flops. The first flop in this composite must flop a rational curve contained in $E$. However, curves of such class sweep out $E$ and cannot be flopped. It follows that $\phi_E$ is an isomorphism.
Finally, we note that since $(\cY^+, \delta \Delta^+)$ is klt, the space $\cY^+$ has rational singularities by \cite[Theorem 5.22]{KollarMori}, and is therefore Cohen--Macaulay by \cite[Theorem 5.10]{KollarMori}, so that $\pi^+ \colon \cY^+ \rightarrow S$ has Cohen--Macaulay source, a smooth target, and central fiber of the expected dimension, so after shrinking $S$,
 	$\pi^+$ is flat. Since $Y^+ \cong Y$ is smooth, $\pi^+$ is also smooth after shrinking $S$ further.
 	\item 
	If $E \in \Delta_Y^{sm}$, then no curve $Z$ in $E$ deforms to $Y^{\textrm{gen}}$ by Proposition~\ref{prop:movablechanges}, so the exceptional locus of $P_E$ is disjoint from the very general fiber. 
  \end{enumerate}
 \end{proof}
The following result establishes the geometric realization of the reflection $\sigma_E$ associated to a contractible $E \in \Delta_Y^{big}$. 
\begin{proposition}\label{prop:geometricrealization}
Let $\Phi_E \colon \cY \dashrightarrow \cY^+$ be the flop of total spaces constructed in Proposition~\ref{prop:flopofroot}. Identifying the central fibers $Y$ and $Y^+$ using the isomorphism $\bar{\phi}_E$, we have an equality	$\Phi_E^0=\sigma_E$ as elements of $\Aut(H^2(Y, \ZZ))$.
\end{proposition}
\begin{proof}
	We consider a general hyperplane section $\bar{H}$ through a general point $p \in C \subset \bar{Y}$, and let $H=p_E^{-1}(\bar{H})$. We let $j \colon \bar{H} \rightarrow \bar{Y}$ and $i \colon H \rightarrow Y$ be the inclusions. 
	Note that $\bar{H}$ has an $A_1$ or $A_2$ singularity at $p$, depending on whether the general fiber of $E \rightarrow C$ is a $\PP^1$ or the union of two $\PP^1$'s. In the $A_1$ case, we write $R=B$ for the exceptional curve of $H \rightarrow \bar{H}$. In the $A_2$ case, we write $B_1$ and $B_2$ for the exceptional curves and set $R=B_1+B_2$. In either case, we have 
	\[
	E|_H=R\quad \textrm{and} \quad i_*[R]=\ell_E.
	\]
	Working locally over a disk $\Delta$, we extend $\bar{H}$ to a relative hyperplane section $\bar{\cH}$ over $(0 \in \Delta)$, and we write $\cH =P_E^{-1}(\bar{\cH})$ and $\cH^+=(P_E^+)^{-1}(\bar{\cH})$. The induced birational map $\Phi_R \colon \cH \dashrightarrow \cH^+$ is the corresponding threefold flop so that
	\[
	\Phi_R^0 \circ i^*=i^* \circ \Phi_E^0.
	\]
	In the $A_1$-case, \cite[Theorem 5.22]{Reid2} shows that $\Phi_R^0$ is the reflection 
	\[
	\sigma_{R} \in \Aut(H^2(H, \ZZ)), \qquad x \mapsto x+ (x \cdot R)R.
	\] 
	In the $A_2$-case, the birational map $\Phi_R$ is the simultaneous flop of $B_1$ and $B_2$. By the Weyl group description (see \cite[Section 4]{KatzMorrison}), this flop factors as the flop of $B_1$, followed by the flop of the strict transform of $B_2$, followed by the flop of the strict transform of $B_1$. With respect to the original lattice, the three roots are
	\[
	B_1,B_1+B_2, B_2,
	\]
	so that the induced action on $H^2(H, \ZZ)$ is given by 
	\[
	\sigma_{B_2}\sigma_{B_1+B_2}	\sigma_{B_1}=\sigma_R
	\] 
	We therefore have $\Phi_R^0=\sigma_R$ in both cases.
	We have that 
	\[
	\sigma_R(i^*x)=i^*x+(i^*x\cdot R)R=i^*x+(x\cdot i_*R)i^*E=i^*x+(x\cdot \ell_E)i^*E=i^*\sigma_E(x)
	\]
	so that $i^*\Phi_E^0(x)=i^*\sigma_E(x)$. It remains to show that $i^*$ is injective. 
	We have a commutative diagram of Picard groups
	\[
	\begin{tikzcd}
		0 \ar[r]& \Pic(\bar{Y}) \ar[r, "p_E^*"] \ar[d,"j^*"] & \Pic(Y)  \ar[d,"i^*" ] \\
		0\ar[r]&	\Pic(\bar{H}) \ar[r, "p_E^*|_H"] & \Pic({H})
	\end{tikzcd} 
	\]
	where the rows are exact and the top row continues as 
	\[
	\begin{tikzcd}
		0 \ar[r]& \Pic(\bar{Y}) \ar[r, "p_E^*"]& \Pic(Y) \ar[r, "(-)\cdot \ell_E"]  &\ZZ,
	\end{tikzcd}
	\]
	(see \cite[Corollary 3.17]{KollarMori}).
	If $i^*L=0$, then 
	\[
	L \cdot \ell_E=L \cdot i_*R=i^*L \cdot R=0.
	\]
	It follows that $L=p_E^*\bar{L}$ for some $\bar{L} \in \Pic(\bar{Y})$ and by commutativity, we have 
	\[
	p_E^*|_H \circ j^*\bar{L}=0
	\]
	Since $\bar{Y}$ has hypersurface singularities, $j^*$ is injective by the Lefschetz theorem (\cite[Theorem, p.24]{Goresky2}). It follows that $L=0$ and $i^*$ is injective.
\end{proof}
By definition, the strict transform of $E \in \Delta_Y^{big}$ under an SQM $\alpha \colon Y \dashrightarrow Z$ is an element of $\Delta_Z^{big}$. The respective reflections commute with $\alpha_*$:
\begin{lemma}\label{lem:Wunderflops}
Let $\alpha \colon Y \dashrightarrow Z$ be an SQM, let $E \in \Delta_Y^{big}$. Let $E_Z \in \Delta_Z^{big}$ be its strict transform. 
Then we have an equality 
	\begin{equation}
 \sigma_E \circ \alpha^*=\alpha^* \circ \sigma_{E_Z}
	\end{equation}
\end{lemma}
\begin{proof}
	Under the dual identification $\alpha_* \colon H_2(Y) \cong H_2(Z)$, the fiber class $\ell_E$ maps to the fiber class $\ell_{E_Z}$. The statement now follows from the duality between $\alpha^*$ and $\alpha_*$.
\end{proof}
Lemma~\ref{lem:Wunderflops} shows that conjugation by $\alpha^*$ induces an isomorphism $W_Y^{big} \cong W_Z^{big}$, so that $W_Y^{big}$ is a birational invariant of $Y$. 
\subsection{The Weyl group action on the movable and the nef cone}
We will now use the geometric description of the action of $W_Y^{big}$ and $W_Y^{sm}$ on $H^2(Y, \RR)$ to show that the action of the big Weyl group preserves the movable cone of $Y^{\mathrm{gen}}$, whereas the action of the small Weyl group preserves the nef cone of $Y^{\mathrm{gen}}$.
\begin{corollary}\label{thm:weylpreservescones}
The action of $W_Y^{big}$ on $H^2(Y, \RR)$ preserves $\overline{\Mov}(Y^{\mathrm{gen}})$, $\Eff(Y^{\mathrm{gen}})$ and  $\Bbig(Y^{\mathrm{gen}})$. 
The action of $W_Y^{sm}$ additionally preserves $\Nef(Y^{\mathrm{gen}})$ and $A(Y^{\mathrm{gen}})$.
\end{corollary}
\begin{proof}
	Suppose first that $E$ is contractible, and $E \in \Delta_Y^{sm}$. In this case, 
	$\Phi_E$ restricts to an isomorphism on a very general fiber by Proposition~\ref{prop:flopofroot}. It follows that $\sigma_E \in \Aut(H^2(Y, \RR))$ is the composite of parallel transport and maps induced by isomorphisms, and therefore preserves $\Nef(Y^{\mathrm{gen}})$, and its interior $A(Y^{\mathrm{gen}})$.\\
	Suppose now that $E$ is contractible, and $E \in \Delta_Y^{big}$. In this case, 
	 $\Phi_E$ restricts to a flop on a very general fiber. It follows that $\sigma_E$ is a composite of parallel transport, and maps induced by a flop and therefore preserves $\overline{\Mov}(Y^{\mathrm{gen}})$.\\
	 If $E$ is not contractible, then we choose again an SQM $\alpha \colon Y \dashrightarrow Z$ so that the strict transform $E_Z$ is contractible on $Z$. By
	 Lemma~\ref{lem:Wunderflops}, we have that 
	 \[
	 \sigma_E=(\alpha^{-1})^* \circ \sigma_{E_Z} \circ \alpha^*.
	 \]
	 $\alpha$ induces an SQM $Y^{\textrm{gen}} \dashrightarrow Z^{\textrm{gen}}$, and therefore preserves the effective, movable, and big cone of a very general fiber. Combined with the previous part, the composite $\sigma_E$ does too.
\end{proof}
We are now ready to prove Theorem~\ref{thm:main2}. Ideally, we would like to show that any class $D$ in $\Nef(Y^{\mathrm{gen}})$ is a $W_Y^{sm}$-translate of $\Nef(Y)$. 
We are unable to prove this in general. However, if we assume that $D$ is also big, we are able to use the minimal model program to show just that. We note that Fujino~\cite{Fujino} has proved that the results of \cite{BCHM} continue to hold for projective morphisms of complex analytic spaces. Recall the notation $\Phi^0$ introduced at the beginning of Section 4.3. 
\begin{theorem}\label{thm:MMP}
Let $\pi \colon \cY \rightarrow S$ be a very general deformation family of Y over a smooth germ $(0 \in S)$. 
	\begin{enumerate}
	\item Let $D \in \Nef(Y^{\mathrm{gen}}) \cap \Bbig(Y^{\mathrm{gen}})$. Possibly after shrinking $S$, there is a sequence of flops 
	\[
	\Phi_D \colon \cY \dashrightarrow \cY_D
	\]
	over $S$ such that the induced map $\phi_D \colon Y \dashrightarrow Y_D$ on central fibers extends to an isomorphism 
	$\bar{\phi}_D \colon Y \cong Y_D$, and 
		\begin{itemize}
			\item ${\Phi^0_{D}} \in W_Y^{sm}$
			\item ${\Phi^0_{D}}(D) \in \Nef(Y) \cap \Bbig(Y)$
		\end{itemize}
	In particular, $D$ is a $W_Y^{sm}$-translate of an element in $\Nef(Y) \cap \Bbig(Y)$. 
		\item 
	Let $D \in \overline{\Mov}(Y^{\mathrm{gen}}) \cap \Bbig(Y^{\mathrm{gen}})$. 
	Possibly after shrinking $S$, there is a sequence of flops 
	\[
	\Phi_D \colon \cY \dashrightarrow \cY_D
	\]
	over $S$ such that the induced map $\phi_D \colon Y \dashrightarrow Y_D$ on central fibers extends to an isomorphism 
	$\bar{\phi}_D \colon Y \cong Y_D$, and 
	\begin{itemize}
	\item ${\Phi^0_{D}} \in W_Y^{big}$
	\item ${\Phi^0_{D}}(D) \in \overline{\Mov}(Y) \cap \Bbig(Y)$
\end{itemize}
			In particular, $D$ is a $W_Y^{big}$-translate of an element in $\overline\Mov(Y) \cap \Bbig(Y)$. 
		\end{enumerate}
\end{theorem}
\begin{proof}
	We start with the first assertion. After shrinking $S$, we identify $H^2(\cY, \ZZ)$ with the cohomology of the central fiber.
	Let $D \in \Nef(Y^{\mathrm{gen}}) \cap \Bbig(Y^{\mathrm{gen}})$. Since $\cD$ is big on a very general fiber, we have that $D \in \Bbig(\cY/S)$.
	By the relative Kodaira lemma, we may write 
	\[
	\cD \sim_{\RR, S} \mathcal{A}+\mathcal{G}
	\]
	where $\mathcal{A}$ is relatively ample and $\mathcal{G}$ is an effective $\RR$-divisor. Shrinking $S$ and replacing $A$ with an $\RR$-linearly equivalent effective divisor, it follows that there is $\Delta' \geq 0$ which is $\RR$-linearly equivalent to $\cD$. Let $\Delta$ be the horizontal part of $\Delta'$, obtained from $\Delta'$ by removing all vertical components. Since vertical components are numerically trivial, we have 
	\[
	[\Delta]=[\Delta']=D \in H^2(Y, \RR).
	\]
	Since $\cY/S$ is smooth, the pair $(\cY, \epsilon \Delta)$ is klt for $0 \leq \epsilon \ll 1$. 
Since $\Delta$ is relatively big, we may run a relative MMP for $K_{\cY/S}+\epsilon \Delta$ by \cite[Theorem 1.7]{Fujino}. Since $K_{\cY/S}=0$, this is a $\Delta$-MMP. After shrinking $S$, we may assume that the exceptional locus of each step of the MMP meets the central fiber.\\
Consider the contraction $P$ of the first extremal ray $R$ in the $\Delta$-MMP, and let $\cE$ be the exceptional locus of $P$.
We have $\Delta \cdot R=D \cdot R<0$, so that $\cE \subset \Supp(\Delta)$. 
Since $\Delta$ is nef on a very general fiber, $\cE$ cannot dominate $S$, and since $\Delta$ has no vertical components, it follows that $\cE$ has codimension at least $2$ in $\cY$ and $P$ is a small contraction. Therefore, the first step in the MMP is a flop 
\[
\begin{tikzcd}
\cY \ar[rr, dashed, "\alpha"] \ar[dr]&& \cY^+ \ar[dl]\\
&S&
\end{tikzcd}
\]
We claim that $\alpha$ is the flop $\Phi_E$ of some $E \in \Delta_Y^{sm}$ constructed in Proposition~\ref{prop:flopofroot}. Indeed, let $Z \subset Y$ be a curve in the exceptional locus of $p=P|_Y$. Note that the corresponding face $F$ of $\Nef(Y)$ satisfies 
$F \subset [Z]^\perp$.
Since $[Z]\in R$, we have that 
\begin{equation}\label{eq:facenefdisappears}
D \cdot Z<0, \qquad D \in \Nef(Y^{\mathrm{gen}}),
\end{equation}
which shows that the face $F$ of $\Nef(Y)$ disappears under deformation. By Wilson's results (in particular, Equality~\eqref{eq:nef}), $R$ is the extremal ray spanned by $\ell_E$ for some $E \in \Delta_Y^{sm}$, so that $p=p_E$ is the contraction of $E$ considered in Proposition~\ref{prop:flopofroot}. It follows that $P=P_E$. Therefore,
$\cY^+ \rightarrow S$ is again a smooth family with central fiber isomorphic to $Y$, and by Proposition~\ref{prop:geometricrealization}
 \[
 [\Delta^+]=\Phi^0_E([\Delta])=\sigma_E([\Delta]),
 \]
where $\Delta^+$ is the strict transform of $\Delta$. 
By Proposition~\ref{prop:flopofroot}, the exceptional locus of $\alpha$ doesn't meet a very general fiber of $\cY$, so that 
\[
[\Delta^+] \in \Nef(Y^{\mathrm{gen}}) \cap \Bbig(Y^{\mathrm{gen}}).
\]
Repeating the same argument, the $\Delta$-MMP terminates after a finite number of flops 
\[
\Phi_D \colon \mathcal{Y}_0 \dashrightarrow \dots \dashrightarrow \mathcal{Y}_n=\mathcal{Y}_D 
\]
Let $Y_i$ be the central fiber of $\cY_i \rightarrow S$. By induction, each birational map $\cY_i \dashrightarrow \cY_{i+1}$ is the flop associated to some root $E_i \in \Delta_{Y_i}^{sm}$ and therefore $Y_i \cong Y_{i+1}$ for all $i$. It follows that all $Y_i$ are isomorphic to $Y$, and using these identifications, we view the $E_i$ as elements of $\Delta_{Y}^{sm}$. It follows that $\phi_D$ acts on $H^2(Y)$ by the product 
of the reflections $\sigma_{E_i}$. 
Therefore
\[
\Phi^0_{D}=\Phi^0_{E_{n-1}}\circ \dots \circ \Phi^0_{E_0}=\sigma_{E_{n-1}}\dots \sigma_{E_0} \in W_Y^{sm}.
\]
Let $\Delta_n$ denote the final strict transform of $\Delta$. Then $\Delta_n$ is nef since the MMP has terminated. Therefore
\[
\Phi^0_{D}(D)=[\Delta_n] \in \Nef(Y). 
\]
Since relative bigness is preserved by the MMP, we also have that 
\[
\Phi^0_{D}(D) \in \Nef(Y) \cap \Bbig(Y^{\textrm{gen}})=\Nef(Y) \cap \Bbig(Y)
\]
where the last equality is Lemma~\ref{lem:ConeInclusions}(4).\\
We now prove the second assertion. 
Let $D \in \overline{\Mov}(Y^{\mathrm{gen}}) \cap \Bbig(Y^{\mathrm{gen}})$. 
The same argument as in the first part yields an effective $\RR$-divisor $\Delta$ with no vertical components and 
\[
[\Delta]=D \in H^2(Y, \RR)
\]
and we run a relative MMP for the pair $(\cY, \epsilon \Delta)$. As before, this is just a $\Delta$-MMP.  
Consider the contraction $P$ of the first extremal ray $R$ in the MMP, and let $\cE$ be the exceptional locus of $P$.
We have $\Delta \cdot R=D \cdot R<0$, so that $\cE \subset \Supp(\Delta)$. 
 Since $D$ is in the closure of the movable cone on a very general fiber, the intersection $\cE \cap \cY_t$ has codimension at least $2$ for $t$ very general: otherwise, curves with class in $R$ would sweep out a divisor on $\cY_t$, and a class $D \in \overline{\Mov}(Y_t)$ has nonnegative intersection with curves sweeping out a divisor. Since $\Delta$ has no vertical components, it follows that $\cE$ has codimension at least $2$ in $\cY$ and $P$ is a small contraction. Therefore, the first step in the MMP is a flop 
\[
\begin{tikzcd}
	\cY \ar[rr, dashed, "\alpha"] \ar[dr]&& \cY^+ \ar[dl]\\
	&S&
\end{tikzcd}
\]
 let $Z \subset Y$ be a curve in the exceptional locus of $p=P|_Y$. Note that the corresponding face $F$ of $\Nef(Y)$ satisfies $F \subset [Z]^\perp$. Since $[Z] \in R$, we have that 
\begin{equation}\label{eq:facemovdisappears}
D \cdot Z<0, \qquad D \in \overline{\Mov}(Y^{\mathrm{gen}}).
\end{equation}
There are now two cases to consider. 
\begin{enumerate}[(a)]
	\item $F$ does not support a codimension one face of $\overline{\Mov}(Y) \cap \Bbig(Y)$.
	\item $F$ supports a codimension one face of $\overline{\Mov}(Y) \cap \Bbig(Y)$.
\end{enumerate}
In case (a), the contraction $p \colon Y \rightarrow \bar{Y}$ on central fibers is small, and the corresponding flop $\phi \colon \cY \dashrightarrow \cY'$ over $S$ restricts to a flop on every fiber. These flops do not correspond to elements of $\Delta_Y^{big}$, and we call them type A flops.\\
In case (b), \eqref{eq:facemovdisappears} shows that the face $F$ of the movable cone disappears under deformation. By Equality~\eqref{eq:mov}, $R$ is the extremal ray spanned by $\ell_E$ for some $E \in \Delta_Y^{big}$, so that $p=P|_Y$ is the contraction $p_E$ considered in Proposition~\ref{prop:flopofroot} and it follows that $P=P_E$. In particular, 
$\cY^+ \rightarrow S$ is again a smooth family with central fiber isomorphic to $Y$, and 
\[
[\Delta^+]=\sigma_E([\Delta])
\]
as before, where $\Delta^+$ is the strict transform of $\Delta$. 
We call these flops type B flops. \\
For both type A and B flops, we have that 
\[
[\Delta^+] \in \overline{\Mov}(Y^{\mathrm{gen}}) \cap \Bbig(Y^{\mathrm{gen}}).
\]
Repeating the same argument, the $\Delta$-MMP terminates after a finite number of type A and B flops 
\[
\Psi_D \colon \mathcal{Y}=\mathcal{Y}_0  \dashrightarrow \dots \dashrightarrow \mathcal{Y}_n=\mathcal{Y'}
\]
Let $Y_i$ be the central fiber of $\cY_i \rightarrow S$. 
The induced map $\Psi_D^0 \colon H^2(Y) \rightarrow H^2(Y')$
factors as a composite of reflections $\sigma_{E_i}$ corresponding to type B flops and isomorphisms $H^2(Y_i) \cong H^2(Y_{i+1})$ induced by type A flops.
The induced birational map $Y_i \dashrightarrow Y_{i+1}$ on central fibers is either a flop (in case (a)) or extends to an isomorphism $Y_i \dashrightarrow Y_{i+1}$ (in case (b)). 
Let 
\[
\psi_D \colon Y \dashrightarrow Y'
\]
be the induced composite SQM on central fibers. 
We may use Lemma~\ref{lem:Wunderflops} to factor 
\[
\Psi_D^0=(\psi_D)_* \circ w
\]
where $w \in W_Y^{big}$ is a product of reflections $\sigma_{E_i}$ for $E_i \in \Delta_{Y}^{big}$. 
By \cite[11.10]{KollarMoriFlips}, the SQM $\psi_D^{-1} \colon Y' \dashrightarrow Y$ deforms to a relative SQM $\Theta_D \colon \cY' \dashrightarrow \cY_D$
for some deformation $\cY_D \rightarrow S$ of $Y$. We set  
 \[
\Phi_D=\Theta_D \circ \Psi_D \colon \cY \dashrightarrow \cY_D.
\]
By construction, the induced birational map $\phi_D \colon Y \dashrightarrow Y_D$ on central fibers extends to an isomorphism and 
\[
\Phi_D^0=(\psi_D^{-1})_* \circ (\psi_D)_* \circ w=w \in W_Y^{big}.
\]
We have that
\[
\Psi^0_{D}(D)=[\Delta_n] \in \Nef(Y')
\]
since the MMP has terminated, and therefore
\[
\Phi^0_D(D)=\Theta_D^0([\Delta_n]) \in \overline{\Mov}(Y),
\]
since $\Theta_D$ restricts to an SQM on central fibers. 
Finally, relative bigness is preserved by the MMP and $\Theta_D$, so we also have that 
\[
\Phi_{D}^0(D) \in \overline{\Mov}(Y) \cap \Bbig(Y^{\textrm{gen}})=\overline{\Mov}(Y) \cap \Bbig(Y)
\]
where the last equality is Lemma~\ref{lem:ConeInclusions}(5).
\end{proof}
We are now ready to give the proof of Theorem~\ref{thm:main2}.
	\begin{proof}[Proof of Theorem~\ref{thm:main2}]
		Recall that the fundamental chamber $\bar{\cC}^{sm}$ is a fundamental domain for the action of $W_Y^{sm}$ on the Tits cone $U^{sm}$. By  Corollary~\ref{thm:weylpreservescones}, the action of $W_Y^{sm}$ preserves $\Nef(Y^{\mathrm{gen}})$, so that $W_Y^{sm}$ acts on $U^{sm} \cap \Nef(Y^{\mathrm{gen}})$ with fundamental domain 
		\[
		\bar{\cC}^{sm} \cap \Nef(Y^{\mathrm{gen}})=\Nef(Y)
		\]
		where the last equality is \eqref{eq:nef}. 
		By Theorem~\ref{thm:MMP}, every element of $\Nef(Y^{\mathrm{gen}}) \cap \Bbig(Y^{\mathrm{gen}})$ is a $W_Y^{sm}$-translate of $\Nef(Y)$, so that
		\[
		\Nef(Y^{\mathrm{gen}}) \cap \Bbig(Y^{\mathrm{gen}}) \subset U^{sm}
		\]
		This shows that 
		$W_Y^{sm}$ acts on $\Nef(Y^{\mathrm{gen}}) \cap \Bbig(Y^{\mathrm{gen}})$ with fundamental domain 
		\[
		\Nef(Y) \cap \Bbig(Y^{\mathrm{gen}})=\Nef(Y) \cap \Bbig(Y),
		\]
		as required (we have used Lemma~\ref{lem:ConeInclusions}(4) for the last equality). \\
		The proof of the second part is exactly analogous, using equality~\eqref{eq:mov} instead.
	\end{proof}
\begin{remark}
	Let $S=\Def(Y)$. By Kawamata~\cite{KawamataCY}, the relative movable cone $\Mov(\cY/S)$ has a decomposition into nef cones of its various relative SQMs. Since $\Mov(\cY/S)=\Mov(Y^{\mathrm{gen}})$, one can ask how this decomposition is related to the decomposition of $\Mov(Y^{\mathrm{gen}})$ into $W_Y^{big}$-chambers. 
	We have seen that $W_Y^{big}$ acts on the set of minimal models of $\cY/S$, so that the decomposition of $\Mov^e(Y^{\mathrm{gen}})$ into nef cones of minimal models is a refinement of the decomposition into Weyl chambers. 
	In general, the decomposition into Weyl chambers is strictly coarser: one Weyl chamber is given by $\Mov(Y)$, whereas one chamber of the decomposition into nef cones is given by $\Nef(\cY/S)=\Nef(Y)$. The decomposition of $\Mov(Y)$ into nef cones is nontrivial as soon as $\Nef(Y)$ has a codimension one face corresponding to a small contraction. The corresponding flop $\phi \colon Y \dashrightarrow Y'$ gives rise to a minimal model of $Y$ and $\cY/S$ that doesn't correspond to an element of $W_Y^{big}$.
\end{remark}
\subsection{The Weyl group and monodromy}\label{sec:monodromy}
We complete this section by explaining the relationship between the Wilson Weyl group and monodromy. 
The results of this section are not used elsewhere in this paper.\\
The automorphism $\sigma_E \in \Aut(H^2(Y, \ZZ))$ is naturally realized by monodromy in a family of Calabi--Yau threefolds.
Indeed, let $p \colon Y \rightarrow \bar{Y}$ denote the contraction of $E$. Since $E$ is an elliptic (quasi)-ruled surface, $\bar{Y}$ has a curve $C$ of uniform ADE singularities of either type $A_1$ or $A_2$. We only give the argument in the former case.
	Let $\pi \colon \cY \rightarrow \Def(Y)$ and $\bar{\pi} \colon \bar{\cY} \rightarrow \Def(\bar{Y})$ denote the Kuranishi families of $Y$ and $\bar{Y}$, respectively. The morphism $\pi$ is smooth, so we can identify the monodromy groups of $Y$ and $Y^{\mathrm{gen}}$ by parallel transport. By \cite[11.4]{KollarMoriFlips}, $\pi$ deforms to give a commutative diagram 
	\[
	\begin{tikzcd}
		\cY \ar[r, "p"] \ar[d, "\pi"] &\bar{\cY}\ar[d, "\bar{\pi}"]\\
		\Def(Y) \ar[r, "f"] & \Def(\bar{Y})
	\end{tikzcd}
	\]
	and since $E$ is ruled over an elliptic curve $C$, \cite[Theorem 2.4]{Namikawa} shows that $f$ is finite. We denote the branch locus $H$. It follows that ${Y}^{\textrm{gen}} \cong \bar{Y}^{\textrm{gen}}$
	and we will show that $\sigma_E$ is realized by monodromy around a loop $\gamma$ around $H$ in $\Def(\bar{Y})$.
	Let $0 \neq v \in T_0\Def(\bar{Y}) \setminus T_0H$. Locally around $C$, $v$ corresponds to a first order deformation $\cY'$ of $\bar{Y}$ over $\Delta_t$ of the form $xy+z^2+tg(w)=0 \subset \CC^4_{xyzw} \times \Delta_t$, for some function $g(w)$, with $g$ not identically zero. 
	To compute the monodromy on $H^2(\cY'_t)$, 
	we may restrict to a tubular neighbourhood $\cV$ of the singular locus $C$ in the total space $\cY'$ and compute $T \colon H^2(\cV_t) \rightarrow H^2(\cV_t)$. Even though $C$ does not deform to the general fiber of $\cY'$ as a complex manifold, it does deform as a smooth manifold, which is topologically a $2$-torus $T^2$. Topologically, $\cV \rightarrow \Delta_t$ has general fiber $\cV_t=S^2 \times T^2 \times B^2$ (where $S^2$ is the topological $2$-sphere, and $B^2$ is the topological $2$-ball), so is homotopy-equivalent to $S^2 \times T^2$. We therefore have an isomorphism
	\[
	H^*(\cV_t) \cong H^*(T^2) \otimes H^*(S^2)
	\]
	of graded rings.
	The manifold $T^2$ in every fiber extends to a (trivial) $T^2$-bundle over $\Delta_t$, so the monodromy acts trivially on $H^*(T^2)$. 
	On the other hand, the restriction of $\cV \rightarrow \Delta \times T^2$ to $\Delta \times p$ has a local equation of the form $xy+z^2+g(p)t+O(t^2)=0$. This is a degeneration to an ordinary double point with vanishing sphere $[S^2]$.
	We may assume that $p \in T^2$ is general, so the linear term $g(p)$ is nonzero, and the monodromy around $0$ is given by a Picard-Lefschetz transformation, so that $T([S^2] \times p)=-([S^2] \times p)$. 
	This completely determines the monodromy action on $H^*(\cV_t)$. 
	Using similar arguments as in \cite[Theorem 3.16]{Voisin2}, one then shows that the monodromy action on $H^*(\cY_t)$ is given by the reflection $\sigma_E$. 
\begin{remark}	The automorphism $\sigma_E$ is sometimes called a ``fibered Dehn twist", and the ``vanishing bundle" $\cV_t \rightarrow T^2$ is known to be a co-isotropic submanifold (see for example \cite{Wehrheim}).
\end{remark}
Elements $\sigma_E \in W_Y^{big}$ are usually not realized by monodromy. Rather, $\sigma_E$ arises as a composite of parallel transport operators $H^2(Y) \rightarrow H^2(Y')$ and morphisms $H^2(Y') \rightarrow H^2(Z)$ induced by flops $Y' \dashrightarrow Z$. This description follows from Proposition~\ref{prop:geometricrealization} and Lemma~\ref{lem:Wunderflops}.
\section{The Weyl group action on the deformation space}\label{sec:H3}
Let $Y$ be a smooth Calabi--Yau threefold, and denote $\cY \rightarrow \Def(Y)$ the Kuranishi family. 
We will construct actions of $W_Y^{big} \rtimes \PsAut(Y)$ on the Kuranishi family, on $\Def(Y)$, and on $H^3(Y, \ZZ)$. 
The action on $H^3(Y, \ZZ)$ is via Hodge isometries and compatible with the period map \[
P \colon \Def(Y) \rightarrow \PP(H^3(Y, \CC)), \qquad s \mapsto[\Omega_s],
\]
where $\Omega_s$ is the image in $H^3(Y, \CC)$ of the class of the holomorphic volume form of $\cY_s$ under parallel transport.
Recall that $P$ is an immersion by the infinitesimal Torelli theorem (see for example \cite[Theorem 10.27]{Voisin1}).


\subsection{Brief review of intersection cohomology}
We recall some results on intersection cohomology that we will use in this section. All references in brackets in what follows refer to \cite{Goresky}. We use $\ZZ$-coefficients throughout, but the discussion is valid for more general rings of coefficients. 
Given a pseudo-manifold $X$ of real dimension $n$, one defines a complex of sheaves $IC^\bullet_X \in D^b(X)$, called the intersection complex with middle perversity (3.1), and the intersection cohomology groups are its hypercohomology groups (we will use cohomological grading):
\[
IH^k(X)=\mathbb H^{k-n}(X,{IC}_X^\bullet).
\]
The construction of $IC^\bullet_X \in D^b(X)$ depends on a stratification of $X$ (1.1), but the resulting complex in the derived category is unique up to canonical isomorphism (4.1). 
If $X$ is nonsingular, then $X \supset \emptyset$ is a stratification of $X$, so that $IC^\bullet_X \cong \underline{\ZZ}_X[n]$, where $\underline{\ZZ}_X$ denotes the constant sheaf. 
\begin{equation}\label{eq:IHsmooth}
	IH^*(X)\cong H^*(X).
\end{equation}
If $j\colon U\hookrightarrow X$ is the inclusion of the nonsingular
locus, the canonical adjunction morphism ${IC}^\bullet_X \rightarrow Rj_*j^*{IC}^\bullet_X$ induces a homomorphism
\[
j^*\colon IH^*(X)\longrightarrow H^*(U).
\]
If $f\colon X\to\bar X$ is a proper small map, then by (6.2), there is an isomorphism
\[
Rf_*{IC}_X^\bullet\cong{IC}_{\bar X}^\bullet,
\]
Taking hypercohomology gives isomorphisms
\begin{equation}\label{eq:IHsmall}
	IH^*(X)\cong IH^*(\bar X).
\end{equation}
Consequently, if $\phi\colon X\dashrightarrow X'$ is the flop of a
flopping contraction $X\to\bar X$, then
\[
\phi_*\colon IH^*(X)\xrightarrow{\sim}IH^*(\bar X)
\xrightarrow{\sim}IH^*(X').
\]
If $X$ and $X'$ are smooth, this gives an isomorphism
\begin{equation}\label{eq:22}
	\phi_*\colon H^*(X,\ZZ)\xrightarrow{\sim}H^*(X',\ZZ).
\end{equation}
We do not assume here that $X$ and $X'$ are compact; in particular, we
will apply \eqref{eq:22} to total spaces of the Kuranishi families of a Calabi--Yau threefold.
Note that in \cite[4.12]{KollarFlops}, Equality~\eqref{eq:22} is proved for $H^*(X, \QQ)$ for threefolds with terminal $\QQ$-factorial singularities. Since we only apply this result to nonsingular complex spaces, the argument works with $\ZZ$-coefficients, by (6.2).
If $j \colon U\hookrightarrow X$ and
$j' \colon U'\hookrightarrow X'$ are the open subsets on which $\phi$ is an isomorphism, then the isomorphism \eqref{eq:22} is compatible with restriction to
$U\cong U'$, in the sense that there is a commutative diagram
\begin{equation}\label{eq:H3restriction}
	\begin{tikzcd}
		H^3(X',\ZZ) \ar[r,"\phi^{-1}_*"] \ar[d, "j'^*"] & H^3(X,\ZZ) \ar[d, "j^*"]\\
		H^3(U',\ZZ) \ar[r,"(\phi|_U)^*"] & H^3(U,\ZZ).
	\end{tikzcd}
\end{equation}
We introduce the following notation. Let
\[
\begin{tikzcd}
	\mathcal Y \ar[r,dashed,"\Phi"] \ar[d] & \mathcal Y' \ar[d]\\
	(S,0) \ar[r,"\bar\phi"] & (S',0)
\end{tikzcd}
\]
be a composite of flops between deformation families of Calabi--Yau threefolds with
central fibers $Y$ and $Y'$. Using \eqref{eq:22} and parallel transport, we obtain an isomorphism
\[
\Phi^{0,3}\colon H^3(Y,\ZZ)\rightarrow H^3(Y',\ZZ).
\]
Throughout this section, we write $\Phi^{0, 2}$, rather than $\Phi^0$ for the corresponding isomorphism $H^2(Y, \ZZ) \rightarrow H^2(Y', \ZZ)$ induced by parallel transport and strict transform of divisor classes defined in Section 4.3.
We have the following compatibility result:
\begin{lemma}\label{lem:compatibilityH3}
Let $\Phi \colon \cY \dashrightarrow \cY'$ be a composite of flops between Kuranishi families of Calabi--Yau threefolds. Let $\mathcal{E}$ be the exceptional locus, and let $s \in \Def(Y)$ be any point such that $\mathcal{E} \cap \cY_s$ has codimension at least two in $\cY_s$. Then the induced isomorphism 
	\[
	H^3(Y, \ZZ) \xrightarrow{PT} H^3(\cY_s, \ZZ) \xrightarrow{\Phi_{s,*}} H^3(\cY'_{s}, \ZZ) \xrightarrow{PT} H^3(Y', \ZZ)
	\]
	agrees with $\Phi^{0,3}$. 
\end{lemma} 
\begin{proof}
	Let $\mathcal{U} \subset \cY$ and $\mathcal{U'} \subset \cY'$ be open subsets on which $\Phi$ is an isomorphism, and write
	\[
	\cU_s=\mathcal{U} \cap \cY_s, \qquad \cU'_s=\mathcal{U'} \cap \cY'_s.
	\]
We have a diagram
\[
\begin{tikzcd}[column sep=large,row sep=large]
	H^3(\cY_s,\ZZ)
	\ar[r,"\sim"]
	\ar[d,"\Phi_{s,*}"]
	&
	H^3(\cY,\ZZ)
	\ar[r]
	\ar[d,"\Phi_*"]
	&
	H^3(\cU,\ZZ)
	\ar[r]
	\ar[d,"\sim"]
	&
	H^3(\cU_s,\ZZ)
	\ar[d,"\sim"]
	\\
	H^3(\cY'_{\bar{\phi}(s)},\ZZ)
	\ar[r,"\sim"']
	&
	H^3(\cY',\ZZ)
	\ar[r]
	&
	H^3(\cU',\ZZ)
	\ar[r]
	&
	H^3(\cU'_{\bar{\phi}(s)},\ZZ).
\end{tikzcd}
\]
The big outer square and the middle square commute by \eqref{eq:H3restriction}, and the right square commutes because the underlying maps of spaces commute. 
Since $\cY_s \setminus \cU_s$ has codimension at least $2$, the restriction map 
\[
j^* \colon H^3(\cY_s, \ZZ) \rightarrow H^3(\cU_s, \ZZ)
\]
is injective, and it follows that the left square commutes. Parallel transport to the central fiber gives the required result. 
\end{proof}
\subsection{The action of $\PsAut(Y)$}
We have the following well-known result:
\begin{proposition}\label{prop:PsAutH3}\cite[4.12]{KollarFlops}
	Let $\phi \colon Y \dashrightarrow Y'$ be a birational map between smooth Calabi--Yau threefolds. Then $\phi$ induces a Hodge isometry $H^3(Y, \ZZ) \xrightarrow{\sim} H^3(Y', \ZZ)$.
	\end{proposition}
\begin{proof}
	Factoring $\phi$ as a sequence of flops and
	using \eqref{eq:22}, we obtain an isomorphism
	\[
	\phi_*\colon H^3(Y,\ZZ)\xrightarrow{\sim}H^3(Y',\ZZ).
	\]
	This isomorphism is independent of the chosen factorization: by
	\eqref{eq:H3restriction}, any two such isomorphisms agree after
	restriction to the common open subset $U$ on which $\phi$ is an
	isomorphism, and $j^*$ and $j'^*$ are injective because $Y \setminus U$ has codimension at least $2$.
	By \cite[4.12, 4.13]{KollarFlops}, $\phi_*$ is an isomorphism of polarized Hodge structures. 
\end{proof}
We will frequently use the following lemma:
\begin{lemma}\label{lem:h2enough}
	Let $\cY \rightarrow \Def(Y)$ be the Kuranishi family of a smooth Calabi--Yau threefold $Y$, and suppose we are given a diagram
\[
\begin{tikzcd}
	\cY \ar[r, dashed, "{\Phi}"] \ar[d, ]& \cY \ar[d]\\
	\Def(Y) \ar[r, "\bar{\phi}"] &\Def(Y)
\end{tikzcd}
\]
where $\Phi$ is a small birational map and $\bar{\phi}$ an isomorphism of germs. Suppose that the restriction of $\Phi$ to the central
fiber extends to the identity on $Y$, and that $\Phi^{0,2}=\id$
on $H^2(Y,\mathbb Z)$. Then
\[
\Phi=\id_{\cY},
\qquad
\bar{\phi}=\id_{\Def(Y)}, \qquad \Phi^{0,3}=\id_{H^3(Y)}.
\]
\end{lemma}

\begin{proof}
Choose a relatively ample line bundle $\cH$ on $\cY$. 
The assumption $\Phi^{0,2}=\id$ implies that the strict transform of $\cH$
under $\Phi$ is isomorphic to $\cH$, so that $\Phi$ is a small birational map preserving a relatively ample line bundle. It follows that $\Phi$
is biregular.
Since the restriction of $\Phi$ to central fibers extends to the identity map, we have that $\bar{\phi}^*\cY \cong \cY$ as deformations. 
By the universality of the Kuranishi family, we have
\[
\bar{\phi}=\id_{\Def(Y)}.
\]
Therefore $\Phi$ is an automorphism of $\cY$ over $\Def(Y)$ whose restriction to
the central fiber is the identity. 
 By the argument in \cite[Theorem 1.11]{Deligne-Mumford}, the scheme $\Aut_H(\cY/\Def(Y))$ of polarized relative automorphisms is quasi-finite over $\Def(Y)$ and since $H^0(Y, T_Y)=0$, it is also unramified over $\Def(Y)$. We therefore must have $\Phi=\id_{\cY}$ (possibly after shrinking $\Def(Y)$). It follows that $\Phi^{0,3}=\id_{H^3(Y)}$ as well.
%
\end{proof}
Our next result shows that every (pseudo)-automorphism deforms to a (pseudo)-automorphism between the Kuranishi families, covering an automorphism of $\Def(Y)$.

\begin{proposition}\label{prop:PsAutbirational}
	Let $Y$ be a smooth Calabi--Yau threefold. 
\begin{itemize}
	\item Every $\theta \in \Aut(Y)$ deforms uniquely to an automorphism 
		\[
	\begin{tikzcd}
		\cY \ar[r, "{\Theta}"] \ar[d, ]& \cY \ar[d]\\
		(\Def(Y), 0) \ar[r, "\bar{\theta}"] &(\Def(Y), 0)
	\end{tikzcd}
	\]
	between Kuranishi families.
\item Every $\theta \in \PsAut(Y)$ deforms uniquely to a birational automorphism 
	\[
\begin{tikzcd}
	\cY \ar[r, "{\Theta}", dashed] \ar[d, ]& \cY \ar[d]\\
	(\Def(Y), 0) \ar[r, "\bar{\theta}"] &(\Def(Y), 0)
\end{tikzcd}
\]
between Kuranishi families. 
\end{itemize}
\end{proposition}
\begin{proof}
	The first part is standard: since $H^0(Y, T_Y)=0$, the versal deformation space $\Def(Y)$ is in fact universal, and therefore any automorphism lifts uniquely to an automorphism of the Kuranishi family. \\
	For the second part, note first that a flop $\phi \colon Y \dashrightarrow Y'$ uniquely deforms to a flop
	\[
	\begin{tikzcd}
		\cY \ar[r, dashed, "{\Phi}"] \ar[d, ]& \cY' \ar[d]\\
		(\Def(Y), 0) \ar[r, "\bar{\phi}"] &(\Def(Y'), 0)
	\end{tikzcd}
	\]
	between the Kuranishi families, and gives rise to an isomorphism $\bar{\phi} \colon \Def(Y) \rightarrow \Def(Y')$ (\cite[12.6.2]{KollarMoriFlips}). Any pseudo-automorphism $\theta$ of a Calabi--Yau threefold can be decomposed as a composite of flops and therefore induces a birational automorphism $(\Theta, \bar{\theta})$ of the Kuranishi family as above. The pair $(\Theta, \bar{\theta})$ is independent of the decomposition of $\theta$ into flops: if $(\Theta, \bar{\theta})$ and $(\Theta', \bar{\theta}')$ arise from two different decompositions of $\theta$ into flops, then $\Theta^{-1} \circ \Theta'$ satisfies the assumptions of Lemma~\ref{lem:h2enough}, and therefore $(\Theta, \bar{\theta})=(\Theta', \bar{\theta}')$.
\end{proof}
\begin{proposition}\label{prop:flopcompatible}
	Let $\phi \colon Y \dashrightarrow Y'$ be a birational map between smooth Calabi--Yau threefolds. 
The induced isomorphisms
	\[
	\phi_* \colon H^3(Y, \ZZ) \cong H^3(Y', \ZZ) \qquad \bar{\phi} \colon \Def(Y) \rightarrow \Def(Y')
	\]
	are compatible with the period map, i.e. $P\circ \bar{\phi}=\PP(\phi_*) \circ P$
\end{proposition}
\begin{proof}
$\phi$ factors as a composite of flops, and we denote $\Phi \colon \cY \dashrightarrow \cY'$ the induced flop of Kuranishi families. Let $\cE$ denote the exceptional locus. 
For every $s \in \Def(Y)$, we have that $\cE \cap \cY_s$ has codimension at least two in $\cY_s$, so by Lemma~\ref{lem:compatibilityH3}, the induced map $\Phi^{0,3}$ factors as
%
\begin{equation}\label{eq55}
H^3(Y, \ZZ) \xrightarrow{PT} H^3(\cY_s, \ZZ) \xrightarrow{\phi_{s*}} H^3(\cY'_{\bar{\phi}(s)}, \ZZ) \xrightarrow{PT} H^3(Y', \ZZ)
\end{equation}
Since $\phi_{s*}$ is an isomorphism of polarized Hodge structures, it sends the line spanned by the holomorphic volume form of $\cY_s$ to the line spanned by the holomorphic volume form of $\cY'_{\bar{\phi}(s)}$, so that we have
\begin{equation}\label{eq3}
	\Phi^{0,3}(\Omega_{s})=\lambda_s\Omega_{\bar{\phi}(s)}
\end{equation}
for some nonzero scalar $\lambda_s$. 
Taking $s=0$ in \eqref{eq55} shows that $\Phi^{0,3}=\phi_*$, and compatibility with the period map follows.
\end{proof}
\subsection{The action of $W_Y^{big}$}
We now want to associate to $E \in \Delta_Y^{big}$ a birational automorphism of $\cY/\Def(Y)$, using the flop of $E$ in the total space of the Kuranishi family constructed in Proposition~\ref{prop:flopofroot}. We start with the case that $E$ is contractible, recall that this means that there exists a contraction $p_E \colon Y \rightarrow \bar{Y}$ of $E$ of relative Picard rank $1$.
\begin{proposition}\label{prop:reflectioncompatible}
 Let $E \in \Delta_Y^{big}$ be contractible. $E$ induces a birational automorphism
		\[
		\begin{tikzcd}
			\cY \ar[r, "{\Phi_E}", dashed] \ar[d, ]& \cY \ar[d]\\
			(\Def(Y), 0) \ar[r, "\bar{\phi}_E"] &(\Def(Y), 0)
		\end{tikzcd}
		\]
		between Kuranishi families. The restriction
		\[
		\Phi_{E, s} \colon \cY_s \dashrightarrow \cY_{\bar{\phi}_E(s)}
		\]
		for $s \in \Def(Y)$ very general is an isomorphism in codimension one. If $E \in \Delta_Y^{sm}$, then $\Phi_{E, s}$ is an isomorphism. \\
		The induced isomorphism $\Phi_E^{0,3} \in \Aut(H^3(Y, \ZZ))$ and $\bar{\phi}_E$ are compatible with the period map.
		In particular, $\bar{\phi}_E$  is an isomorphism, and $\Phi_E^{0,3}$ is a Hodge isometry.
	\end{proposition}
	\begin{proof}
			Let
		\[
		\begin{tikzcd}
			\cY \ar[rr, dashed, "\Phi_E"] \ar[dr]&& \cY^+ \ar[dl]\\
			&S&
		\end{tikzcd}
		\]
		be the flop constructed in Proposition~\ref{prop:flopofroot}. Since $\cY^+$ is again a deformation of $Y$, there is a unique morphism of germs
		\[
		\bar{\phi}_E \colon (\Def(Y), 0) \rightarrow (\Def(Y), 0)
		\]
		such that $\bar{\phi}_E^* \cY \cong \cY^+$. This gives a diagram 
		\[
		\begin{tikzcd}
			\cY \ar[d] \ar[r, dashed, "\Phi_E"]& \cY \ar[d]\\
			(\Def(Y), 0) \ar[r, "\bar{\phi}_E"] &(\Def(Y), 0)
		\end{tikzcd}
		\]
		By construction, the restriction of $\Phi_E$ to central fibers extends to the identity map. Let $\cE$ be the exceptional locus of $\Phi_E$. 
		Since $E$ does not deform, $\cE \cap \cY_s$ has codimension at least two for $s$ very general, so that $\Phi_{E, s}$ is an isomorphism in codimension one.
		If $E \in \Delta_Y^{sm}$, then no curve $Z$ with class proportional to the fiber class $\ell_E$ of $E$ deforms generically, so that $\cE \cap \cY_s=\emptyset$ for $s$ very general, and $\Phi_{E, s}$ is an isomorphism. By Lemma~\ref{lem:compatibilityH3}, $\Phi_E^{0,3}$ equals the composite
		\[
			H^3(Y) \xrightarrow{PT} H^3(\cY_s) \xrightarrow{\Phi_{E, s,*}} H^3(\cY_{\bar{\phi}_E(s)}) \xrightarrow{PT} H^3(Y) 
		\]
		for $s$ very general. 
		Since $\Phi_{E, s,*}$ is an isomorphism in codimension one, it is a Hodge isometry by Proposition~\ref{prop:PsAutH3}. As in Proposition~\ref{prop:flopcompatible}, we find that
		\begin{equation}\label{eq:continuity}
			\Phi_E^{0,3}(\Omega_{s})=\lambda_s\Omega_{\bar{\phi}_E(s)}
		\end{equation}
		for $s \in S^0$, where $S^0 \subset \Def(Y)$ is the dense subset over which $E$ doesn't deform. This shows that 
		\begin{equation}\label{eq:compatibility}
		P\circ \bar{\phi}_E|_{S^0}=\PP(\Phi_E^{0,3}) \circ P|_{S^0}
		\end{equation}
		By continuity, \eqref{eq:compatibility} must then hold for all $s \in \Def(Y)$, and therefore, \eqref{eq:continuity} holds for all $s \in \Def(Y)$ as well. Since $P$ is an immersion by the infinitesimal Torelli theorem, the analytic inverse function theorem shows that $\bar{\phi}_E$ is invertible. \\
		To check that $\Phi_E^{0,3}$ is an isometry, note that the polarization $Q$ on $H^3(Y)$ is just given by the cup product pairing, and every morphism in the definition of $\Phi_E^{0,3}$ preserves this pairing. Next, setting $s=0$ in \eqref{eq:continuity} shows that $\Phi_E^{0,3}$ preserves $H^{3,0}$. Since $\Phi_E^{0,3}$ is defined over $\ZZ$, it also preserves $H^{0, 3}=\overline{H^{3,0}}$. 
		By Griffiths transversality (see for example \cite[Proposition 10.12]{Voisin1}), the derivative of the period map is identified with the inclusion 
		\[
		H^{2,1} \hookrightarrow H^{2,1} \oplus H^{1,2} \oplus H^{0,3}.
		\]
		Since $\Phi_E^{0,3}$ is compatible with the period map, it preserves the subspace $F^2=H^{3,0} \oplus H^{2,1}$. Since $\Phi_E^{0,3}$ preserves the polarization $Q$ and 
		\[
		H^{2,1}=\left\{ x \in F^2 \mid Q(x, H^{0,3})=0 \right\},
		\]
		$\Phi_E^{0,3}$ must also preserve $H^{2,1}$, and therefore also $H^{1,2}=\overline{H^{2,1}}$.
		So $\Phi_E^{0,3}$ is a Hodge isometry.
		\end{proof}
 \begin{proposition}
 	Let $E \in \Delta_Y^{big}$ be arbitrary. Then Proposition~\ref{prop:reflectioncompatible} continues to hold. 
 \end{proposition}
\begin{proof}
Choose an SQM $\psi \colon Y \dashrightarrow Z$ such that $E_Z=\psi_*E$ is contractible. $\psi$ deforms to an SQM $\Psi \colon \cY \dashrightarrow \cZ$ between Kuranishi families, and we define $\Phi_E$ and $\bar{\phi}_E$ as the composite 
	\[
	\begin{tikzcd}
		\cY \ar[d] \ar[r, dashed, "\Psi"]&	\cZ \ar[d] \ar[r, dashed, "\Phi_{E_Z}"]& \cZ \ar[d] \ar[r, dashed, "\Psi^{-1}"]& \cY \ar[d]\\
		(\Def(Y), 0) \ar[r, "\bar{\psi}"]&	(\Def(Z), 0) \ar[r, "\bar{\phi}_{E_Z}"] &(\Def(Z), 0)\ar[r, "\bar{\psi}^{-1}"]&	(\Def(Y), 0)
	\end{tikzcd}
	\]
	We claim that the definition of $\Phi_E$ and $\bar{\phi}_E$ is independent of the choice of SQM. Indeed, the induced automorphism $\Phi^{0,2}_E$ of $H^2(Y, \ZZ)$ is 
	\[
	\psi^{-1}_* \circ \sigma_{E_Z} \circ \psi_*=\sigma_E
	\]
	 and hence independent of the SQM. By Lemma~\ref{lem:h2enough}, $\Phi_E$, $\bar{\phi}_E$, and $\Phi^{0,3}_E$ are independent of the SQM as well.
\end{proof}
\begin{remark}
	If $E$ is ruled, or quasi-ruled, the automorphism $\bar{\phi}_E$ of $\Def(Y)$ has appeared in the work of Szendr\H{o}i~\cite{Szendroi}. 
\end{remark}
Proposition~\ref{prop:reflectioncompatible} constructs for each generator $\sigma_E \in W_Y^{big}$ a birational automorphism of $\cY/\Def(Y)$. We now extend this to an action of the whole group $W_Y^{big}$.
\begin{proposition}\label{prop:Weylgroupacts}
	The assignments
	\[
	E\longmapsto \Phi_E,\qquad
	E\longmapsto \bar{\phi}_E,\qquad
	E \longmapsto \Phi_E^{0,3}
	\]
	respect the Coxeter relations, and therefore extend to actions of $W_Y^{{big}}$ on the Kuranishi
	family $\cY\rightarrow\Def(Y)$, on $\Def(Y)$, and $H^3(Y, \ZZ)$. 
	The action on $H^3(Y,\mathbb Z)$ is by Hodge isometries and is compatible
	with the period map. For $w \in W_Y^{big}$, write 
	\[
	\begin{tikzcd}
		\cY \ar[r, "{\Phi_w}", dashed] \ar[d, ]& \cY \ar[d]\\
		(\Def(Y), 0) \ar[r, "\bar{\phi}_w"] &(\Def(Y), 0)
	\end{tikzcd}
	\]
	 for the action on the Kuranishi family. The restriction
	 \[
	 \Phi_{w, s} \colon \cY_s \dashrightarrow \cY_{\bar{\phi}_w(s)}
	 \]
	 for $s \in \Def(Y)$ very general is an isomorphism in codimension one. If $w \in W_Y^{sm}$, then $\Phi_{w, s}$ is an isomorphism. 
\end{proposition}
\begin{proof}
Given $w \in W_Y^{big}$, write $w=\sigma_{E_1} \dots \sigma_{E_n}$ as a product of simple reflections, and define $\Phi_w$ and $\bar{\phi}_w$ as the corresponding composites. To see that this is a well-defined action, it suffices to check that Coxeter relations act trivially. 
If $w=(\sigma_E\sigma_{E'})^{m(E, E')}=1$, then we have that 
$\Phi_w^{0,2}=\id_{H^2(Y)}$. By construction, $\Phi_w$ extends to the identity on central fibers, so Lemma~\ref{lem:h2enough} then shows that $\Phi_w=\bar{\phi}_w=\id$. It follows that $\Phi_w^{0,3}=\id$ as well.
\end{proof}
\subsection{The action of $W_Y^{big}\rtimes \PsAut(Y) $.}
The image of an SQM-conic fibration surface under a pseudo-automorphism is again an SQM-conic fibration surface, so that $\PsAut(Y)$ naturally acts on $\Delta_Y^{big}$. This action preserves the relations in $W_Y^{big}$, so that $\PsAut(Y)$ acts on $W_Y^{big}$, we write this action as 
\[
\theta(w)=\theta_* \circ w \circ \theta_*^{-1}.
\]
The group operation in $W_Y^{big}\rtimes \PsAut(Y) $ is given by 
\[
(w_1, \theta_1)(w_2, \theta_2)=(w_1\theta_1(w_2), \theta_1\theta_2).
\]
The image of a (quasi)-ruled surface under an automorphism is again (quasi)-ruled, so that $\Aut(Y)$ acts on $\Delta_Y^{sm}$, and also on $W_Y^{sm}$, so that 
\[
W_Y^{sm} \rtimes \Aut(Y) \leq W_Y^{big} \rtimes \PsAut(Y) .
\]
\begin{theorem}\label{thm:biggroupacts}
	The actions of $W_Y^{big}$ and $\PsAut(Y)$ constructed above combine to compatible actions of
	$W_Y^{big}\rtimes\PsAut(Y)$
	on the Kuranishi
	family $\cY\rightarrow\Def(Y)$, on $\Def(Y)$, and $H^3(Y, \ZZ)$. An element $(w, \theta)$ acts on the Kuranishi family as 
				\[
	\begin{tikzcd}
		\cY \ar[d] \ar[r, "\Theta", dashed]&\cY \ar[d] \ar[r, "{\Phi_w}", dashed] &\cY \ar[d, ]\\
		(\Def(Y), 0) \ar[r, "\bar{\theta}"] &(\Def(Y), 0) \ar[r, "\bar{\phi}_w"] &(\Def(Y), 0) 
	\end{tikzcd}
	\]
	The action on $H^3(Y,\mathbb Z)$ is by Hodge isometries and compatible
	with the period map.
	The restriction
	\[
	 \cY_s \dashrightarrow \cY_{\bar{\phi}_w \circ \bar{\theta}(s)}
	\]
	for $s \in \Def(Y)$ very general is an isomorphism in codimension one. If $(w, \theta) \in W_Y^{sm} \rtimes \Aut(Y)$, then it is an isomorphism. 
%
\end{theorem}
	\begin{proof}
		By Propositions~\ref{prop:PsAutH3}, ~\ref{prop:PsAutbirational}, and~\ref{prop:flopcompatible}, the groups $\Aut(Y)$ and $\PsAut(Y)$ act on the Kuranishi family, on $\Def(Y)$, and on
		$H^3(Y,\mathbb Z)$ as claimed.
		By Proposition~\ref{prop:Weylgroupacts}, the groups $W_Y^{sm}$ and $W_Y^{big}$ act as claimed.  It remains to check that the semidirect product acts. It is enough to show that for $(w, \theta) \in W_Y^{big} \rtimes \PsAut(Y)$, the induced maps of total spaces satisfy
		\begin{equation}\label{eq:conjugation}
		\Phi_{\theta(w)}
		=
		\Theta\circ\Phi_w\circ\Theta^{-1}.
		\end{equation}
		To see this, we note that the induced maps on $H^2(Y,\mathbb Z)$ agree, since
		\[
		(\Theta\circ\Phi_w\circ\Theta^{-1})^{0,2}
		=
		\theta_*\circ w\circ\theta_*^{-1}
		=
		\theta(w)
		=
		\Phi_{\theta(w)}^{0,2}.
		\]
		Moreover, both sides of~\eqref{eq:conjugation} restrict on the central fiber to maps
		extending to the identity on $Y$. Lemma~\ref{lem:h2enough} then proves~\eqref{eq:conjugation}.
	\end{proof}
We now arrive at the main result of this section, which shows that even though  $W_Y^{big} \rtimes \PsAut(Y) $ is often an infinite group, it always acts on the deformation space via a finite group.
\begin{proposition}\label{prop:actsfinitegroup}
	Let $Y$ be a smooth Calabi--Yau threefold. Then 
	\begin{enumerate}
		\item $W_Y^{big} \rtimes \PsAut(Y)$ acts on $\Def(Y)$ via a finite group.
		\item  $W_Y^{sm} \rtimes \Aut(Y)$ acts on $\Def(Y)$ via a finite group.
	\end{enumerate}
In particular, a finite index subgroup lifts to pseudo-automorphisms (resp. automorphisms) of a very general deformation of $Y$.
\end{proposition}
\begin{proof}
Since $W_Y^{sm} \rtimes \Aut(Y)$ is a subgroup of $W_Y^{big} \rtimes \PsAut(Y)$, the first statement clearly implies the second.	In view of the compatibility with the period map established in Proposition~\ref{thm:biggroupacts}, it is enough to check that $W_Y^{big} \rtimes \PsAut(Y)$ acts via a finite group on $H^3(Y, \CC)$. 
	Since $H^1(Y, \CC)=0$, we must have $H^5(Y, \CC)=0$, so that cup product with a  $(1,1)$-class $\omega$ annihilates $H^3(Y, \CC)$, so that $H^3(Y, \CC)$ is primitive. It follows that each of the Hodge summands $H^{p,q}(Y)$ of $H^3(Y, \CC)$ carries a canonical positive definite Hermitian form (see \cite[Theorem 6.32]{Voisin1}). Note that the form on the middle primitive cohomology does not depend on a choice of K\"ahler class. Since $W_Y^{big} \rtimes \PsAut(Y)$ acts via Hodge isometries, it follows that the image of $W_Y^{big} \rtimes \PsAut(Y)$ in $\GL(H^3(Y, \CC))$ is contained in the compact group corresponding to the product of the unitary groups acting on the $H^{p,q}$. However, the action also preserves the discrete group $H^3(Y, \ZZ)_{\text{tf}}$ (where $M_{\text{tf}}=M/\text{torsion})$, so the image is finite. 
\end{proof}
\begin{remark}
In fact, a similar argument shows that $\Aut(Y)$ acts via a finite group on $\Def(Y)$ for any smooth Calabi--Yau manifold $Y$ of dimension $n \geq 3$ (recall that this means that $K_Y=0$ and $h^i(\cO_Y)=0$ for $0 < i <n$): indeed, by type, we see that the holomorphic volume form is contained in the orthogonal complement $S$ of the image of $H^{2}(Y, \CC) \times H^{n-2}(Y, \CC)$ in $H^n(Y, \CC)$ (since $H^2(Y, \CC)=H^{1,1}(Y)$). Since $S$ is clearly contained in the primitive cohomology, it carries a canonical positive definite Hermitian form. Now argue as in Proposition~\ref{prop:actsfinitegroup}, using the infinitesimal Torelli theorem.
\end{remark}
Proposition~\ref{prop:actsfinitegroup} shows that the actions of $W_Y^{big}$ on $H^2(Y, \RR)$ and $H^3(Y, \RR)$ are very different. While the former action is always faithful (see Theorem~\ref{thm:Weylgroupactionisfaithful}), the latter acts via a finite group. Note that $W_Y^{big}$ is often infinite (see Section~\ref{sec:examples}).
As an application, we give an affirmative answer to a question raised in \cite[Remarks 4.7]{WilsonCone}: for each of the countably many classes $E_i \in H^2(Y, \ZZ)$ which can represent a SQM-conic fibration surface on some small deformation of $Y$, let $\Gamma_i \subset \Def(Y)$ be the submanifold over which $E_i$ deforms.
\begin{corollary}\label{cor:finitearrangement}
The arrangement of the $\Gamma_i$ in $\Def(Y)$ is locally finite.
\end{corollary}
\begin{proof}
$\Gamma_i \subset \Def(Y)$ coincides with the fixed locus of $\Phi^{0,3}_{E_i}$ acting on $\Def(Y)$. Since $W_Y^{big}$ acts via a finite group, the number of $\Gamma_i$ must be finite. 
\end{proof}
We remark that for $E \in \Delta_Y^{sm}$, the automorphism $\Phi^{0,3}_{E}$ of $H^3(Y)$ is realized by monodromy. Indeed, denoting $\pi \colon Y \rightarrow \bar{Y}$ the contraction of $E$, with $C=\pi(E)$, we have a decomposition 
	\[
H^3(Y, \QQ) \cong \pi^*H^3(\bar{Y}, \QQ) \oplus H^1(C, \QQ)
\]
where $H^1(C) \rightarrow H^3(Y)$ is the composite of pullback $(\pi|_E)^*  \colon H^1(C) \rightarrow H^1(E)$ and the Gysin morphism $H^1(E) \rightarrow H^3(Y)$.
By \cite{Szendroi}, $\Phi^{0,3}_{E}$ acts by $(-1)$ on the summand $H^1(C)$, and by $1$ on the summand $H^3(\bar{Y})$, and one checks that this agrees with the monodromy action computed in Section~\ref{sec:monodromy}.
\section{The Morrison cone conjecture under deformation}\label{sec:MCC}
\subsection{The work of Looijenga}
We recall a few definitions and theorems from \cite{LooijengaCones}. Let $V$ be a real finite dimensional vector space equipped with a rational structure $V(\QQ) \subset V$, and $\cC$ an open nondegenerate convex cone in $V$. We also assume that $\Gamma$ is a subgroup of $\GL(V)$ which preserves $\cC$ and some lattice in $V(\QQ)$. 
We define $\cC^+=\text{Conv}(\bar{\cC} \cap V(\QQ))$.
\begin{proposition}[\cite{LooijengaCones} Proposition-Definition 4.1]\label{prop:Looijenga4.1}
	The following are equivalent:
	\begin{enumerate}[(i)]
		\item There exists a polyhedral cone $\Pi$ in $\cC^+$ with $\Gamma \cdot \Pi=\cC^+$.
		\item There exists a polyhedral cone $\Pi$ in $\cC^+$ with $\Gamma \cdot \Pi \supset \cC$.
	\end{enumerate}
	Moreover, in case (ii) we necessarily have $\Gamma \cdot \Pi=\cC^+$. If one of these equivalent conditions is fulfilled, we say that $(V(\QQ),\cC,\Gamma)$ is a polyhedral triple or simply, that $(\cC^+,\Gamma)$ is of polyhedral type.
\end{proposition}
An important consequence of this proposition is that if $(\cC^+, \Gamma)$ is of polyhedral type, the action of $\Gamma$ admits a rational polyhedral fundamental domain (henceforth abbreviated to RPFD).
\begin{theorem}\cite[Application 4.14]{LooijengaCones}\label{thm:RPFDcriterion}
	Suppose that $(V(\QQ),\cC,\Gamma)$ is a polyhedral triple. Then the action of $\Gamma$ on $\cC^+$ admits an RPFD.
\end{theorem}

The RPFD in Theorem~\ref{thm:RPFDcriterion} can be constructed explicitly as follows: given any integral class $\xi \in (\cC^*)^\circ$, Looijenga constructs a rational polyhedral cone
\[
\Pi_\xi=\{x \in \cC^+ \mid \xi (\gamma x) \geq \xi(x) \; \text{for all} \; \gamma \in \Gamma\}
\]
If the stabilizer $\Gamma_\xi$ is trivial, $\Pi_\xi$ is an RPFD for the action of $\Gamma$ on $\cC^+$.
Since the action of $\Gamma$ on $\cC$ is properly discontinuous, this condition is satisfied by a generic rational point of $(\cC^*)^\circ$.
We will also use the following result:
\begin{theorem} \cite[Proposition 4.6]{LooijengaCones} \label{thm:RPFDfiniteindex}
	Let $(\cC^+, \Gamma)$ be a pair of polyhedral type, and let $\Gamma'$ be a subgroup of $\GL(V(\QQ))$ which stabilizes $\cC$. Then $(\cC^+, \Gamma')$ is of polyhedral type if and only if $\Gamma'$ is commensurable with $\Gamma$.
\end{theorem}
We recall here that two subgroups of some group are said to be commensurable if their intersection is of finite index in each of them.
\begin{remark}\label{rem:subgroup}
	Keeping the assumptions of Theorem~\ref{thm:RPFDfiniteindex} we see that if either $\Gamma' \subset \Gamma$ or $\Gamma \subset \Gamma'$ is a subgroup of finite index, $(\cC^+,\Gamma')$ is of polyhedral type as well. If on the other hand $\Gamma'$ contains $\Gamma$ as a subgroup, not necessarily of finite index, but we happen to know that $\Gamma'$ preserves a lattice in $V(\QQ)$, it is immediate from the definitions that $(\cC^+,\Gamma')$ is of polyhedral type. 
	Theorem~\ref{thm:RPFDfiniteindex} then has the nontrivial implication that the index of $\Gamma$ in $\Gamma'$ must in fact have been finite.
\end{remark}
\subsection{The cone conjecture for Calabi--Yau threefolds}
Let $Y$ be a smooth Calabi--Yau threefold. We take $V=H^2(Y, \RR)\cong \Pic(Y)_\RR$ with its natural rational structure $H^2(Y, \QQ)$, and let $\cC=A(Y)$ be the ample cone, so that $\cC^+=\Nef^+(Y)$.
Unlike for $K3$ surfaces, the group $\Aut(Y)$ generally does not act faithfully on $H^2(Y)$, although the kernel $K$ is at worst a finite group: 
\begin{lemma}\label{lem:kernelfinite}
	Let $Y$ be a smooth Calabi--Yau threefold.
	\begin{itemize} 
		\item The kernel of the action of $\Aut(Y)$ on $H^2(Y)$ is a finite group.
		\item The kernel of the action of $\PsAut(Y)$ on $H^2(Y)$ is a finite group.
	\end{itemize}
\end{lemma}
\begin{proof}
	The first part is standard (see for example \cite[Proposition 2.4]{Oguiso}). If an element of $\PsAut(Y)$ acts trivially on $H^2(Y)$, it must fix an ample divisor and is therefore biregular. The second part follows.
\end{proof}
For a representation $\rho \colon G \rightarrow \GL(V)$ we write $\bar{G}=G/\ker(\rho)$.
Taking $\Gamma=\oAut(Y)$, Theorem~\ref{thm:RPFDcriterion} says that the nef cone conjecture is equivalent to the statement that $(\Nef^+(Y), \oAut(Y))$ is of polyhedral type.\\
In a similar vein, one can take $\cC$ to be the interior of the movable cone of $Y$, and $\Gamma=\oPsAut(Y)$. Then the movable cone conjecture is equivalent to the statement that $(\Mov^+(Y), \oPsAut(Y))$ is of polyhedral type. \\
For later use, we note the following easy lemma:
\begin{lemma}\label{lem:kerneltrivial} Let $Y$ be a smooth Calabi--Yau threefold. Then 
	\begin{itemize}
	\item $\overline{W_Y^{sm} \rtimes \Aut(Y)}=W_Y^{sm} \rtimes \overline{\Aut(Y)} $.
	\item $\overline{W_Y^{big} \rtimes \PsAut(Y)}=W_Y^{big} \rtimes \overline{\PsAut(Y)} $.
	\end{itemize}
\end{lemma}
\begin{proof}
Suppose $(w,g)$ acts trivially on $H^2(Y)$. Then $w$ preserves the ample cone $A(Y)$ of $Y$. Since $A(Y) \subset \Bbig(Y) \cap \Nef(Y)$, we must have $w=1$ by Theorem~\ref{thm:main2}, and therefore $g$ acts trivially on $H^2(Y)$. The proof of the second item is analogous.
\end{proof}

\begin{corollary}\label{cor:relationMCC}
	Let $Y$ be a smooth Calabi--Yau threefold.
	\begin{itemize}
		\item 
		Suppose that the action of $\overline{\Aut}(Y)$ on $\Nef^e(Y)$ admits an RPFD. Then the action of $\overline{\Aut}(Y)$ on $\Nef^+(Y)$ admits an RPFD. 
		\item Suppose that the action of $\overline{\PsAut}(Y)$ on $\Mov^e(Y)$ admits an RPFD. Then the action of $\overline{\PsAut}(Y)$ on $\Mov^+(Y)$ admits an RPFD. 
	\end{itemize}
\end{corollary}
\begin{proof}
	If the action of $\Gamma=\overline{\Aut}(Y)$ on $\Nef^e(Y)$ admits an RPFD $\Pi$, then clearly we have $\Gamma \cdot \Pi \supset A(Y)$, where $A(Y)$ denotes the ample cone. 
	It follows from Proposition~\ref{prop:Looijenga4.1} that $(\Nef^+(Y), \Gamma)$ is of polyhedral type, and therefore admits an RPFD by Theorem~\ref{thm:RPFDcriterion}. 
	The proof of the second part is identical.
\end{proof}
\subsection{The cone conjecture under deformation}
Throughout this section, $S$ denotes an algebraic variety.
Let $\pi \colon \cY \rightarrow S$ be a family of Calabi--Yau threefolds, and let $\eta \in S$ be the generic point. Let $K=k(\eta)$ be the function field of $S$, and $\bar{K}$ its algebraic closure. We denote the generic fiber by $Y_\eta$ and let $Y_{\bar{\eta}}=Y_\eta \times_K \bar{K}$ be the geometric generic fiber, and $p \colon Y_{\bar{\eta}} \rightarrow  Y_\eta$ the induced morphism. 
The following result is standard.
\begin{lemma}\label{lem:AutomorphismDefinedOverFiniteExtension}
	Let $Y$ be a projective variety defined over a field $K$. Write $\bar{Y}=Y \times_K \bar{K}$ for the base change to the algebraic closure. 
	\begin{enumerate}
		\item Suppose that $\Aut(\bar{Y})$ is finitely generated. Then there exists a finite extension $E/K$ such that $\Aut(Y \times_K E) \cong \Aut(\bar{Y})$.
		\item The statement of (1) holds with $\Aut$ replaced by $\PsAut$.
	\end{enumerate}
\end{lemma}
\begin{proof}
	For the first claim, choose generators $\theta_1, \dots, \theta_r$ for $\Aut(\bar{Y})$, and for each $i$, let
		\[
		\Gamma_i \subset \overline{Y} \times_{\overline{K}} \overline{Y}
		\]
		be the graph of $\theta_i$.  
		Each subscheme $\Gamma_i$ is defined over a finite extension $E_i/K$, so passing to a common finite extension $E/K$, we may assume that each $\Gamma_i$ has a model $\Gamma_{i,E}$ over $E$. 
The two projections
\[
\Gamma_{i,E}\longrightarrow Y_E
\]
become isomorphisms after base change to $\bar{Y}$, since $\Gamma_i$ is the graph of $\theta_i$. By faithfully flat descent (see~\cite{stacks-project}), they are already isomorphisms over $E$, and hence $\Gamma_{i,E}$ is the graph of an automorphism of $Y_E$.\\
		For the second claim, choose generators $\theta_1, \dots, \theta_r$ for $\PsAut(\bar{Y})$. By definition, there are open subsets $\bar{U}_i$, $\bar{V}_i$ whose complements have codimension at least two, and isomorphisms 
		\[
	{\theta_i}|_{\bar{U}_i} \colon	\bar{U}_i \xrightarrow{\sim} \bar{V}_i. 
		\]
		Arguing as above, after passing to a finite extension $E/K$, we may assume that the open sets and the isomorphisms are defined over $E$, we therefore obtain open sets
	$U_{i,E},V_{i,E} \subset Y_E$
	whose complements have codimension at least two, together with
	isomorphisms
$
	U_{i,E} \xrightarrow{\sim} V_{i,E}.
$
	It follows that each $\theta_i$ descends to a pseudo-automorphism of $Y_E$. In either case, the natural map 
	\[
	\Aut(Y_E) \rightarrow \Aut(\bar{Y}), \qquad \text{respectively}\qquad \PsAut(Y_E) \rightarrow \PsAut(\bar{Y})
	\]
	is injective, and we showed that the image contains the chosen generators, so it is an isomorphism.
\end{proof} 
We next relate the automorphism group of $Y$ to the automorphism group of a very general deformation $Y^{\mathrm{gen}}$. 
\begin{proposition}\label{prop:genericvsgeneral}
			Let $\cY \rightarrow S$ be a family of Calabi--Yau threefolds.
	\begin{enumerate} 
		\item Suppose that $\Aut(Y^{\mathrm{gen}})$ is finitely generated. After a base change by a generically finite morphism $S' \rightarrow S$, with $S'$ smooth, the automorphism groups $\Aut(Y^{\mathrm{gen}})$ of a very general fiber of $\pi$ and the automorphism group $\Aut(\cY'_{\eta})$ are identified, where $\cY'_\eta$ is the generic fiber of $\cY':=\cY \times_{S} S'$.
	In particular, any automorphism of $Y^{\mathrm{gen}}$ lifts to a birational automorphism of $\cY'$. 
	\item The statement of (1) holds with $\Aut$ replaced by $\PsAut$.
	\end{enumerate}
\end{proposition}
\begin{proof}
 $\bar{K}$ is an algebraically closed field of the same cardinality as $\CC$, so there is an isomorphism $\alpha \colon \bar{K} \cong \CC$. By \cite[Lemma 2.1]{Vial}, there is an isomorphism $\tilde{\alpha} \colon Y_{\bar{\eta}} \cong Y^{\mathrm{gen}}$ over $\alpha$. Conjugating with $\alpha$ then defines an isomorphism $\Aut(Y^{\mathrm{gen}}) \cong \Aut(Y_{\bar{\eta}})$. Since $\Aut(Y^{\mathrm{gen}})$ is finitely generated, Lemma~\ref{lem:AutomorphismDefinedOverFiniteExtension} ensures that there is a finite extension $E/K$ such that $\Aut(Y_{\bar{\eta}}) \cong \Aut(Y_\eta \times_K E)$. Let $\bar{S}$ be the normalization of $S$ in $E$, then $p \colon \bar{S} \rightarrow S$ is a finite surjective map, and the pullback of $\cY$ to $\bar{S}$ has generic fiber $Y_\eta \times_K E$. Let $q \colon S' \rightarrow \bar{S}$ be a resolution of singularities. The pullback $\cY'=q^*\cY$ is the required family.
The proof for $\PsAut$ is completely analogous. 
\end{proof}
\begin{proposition}\label{prop:genericautomorphismgroupsubgroup}
Let $\pi \colon \cY \rightarrow S$ be a very general deformation family of $Y$ over a smooth germ $(0 \in S)$.
There are injective group homomorphisms
	\begin{align}
	 \Aut(\cY_\eta) &\rightarrow W_Y^{sm} \rtimes \Aut(Y) \\
	 \PsAut(\cY_\eta) &\rightarrow W_Y^{big} \rtimes \PsAut(Y).
	\end{align}
\end{proposition}
\begin{proof}
	An element $\Theta \in \Aut(\cY_\eta)$ induces a birational automorphism $\Theta \colon \cY \dashrightarrow \cY$ over $S$. We consider the restriction of $\Theta$ to the inverse image of a small open ball $U \subset S$ (in the analytic topology) centered at $0$, and as usual denote 
	\[
	\Theta^{0} \colon H^2(Y) \rightarrow H^2(Y)
	\]
	the isomorphism induced by parallel transport, followed by strict transform, followed by parallel transport. 
	Let $D \in A(Y)$. Since $\Theta$ restricts to an automorphism on a very general fiber, we have 
	\[
	\Theta^{0}(D) \in A(Y^{\mathrm{gen}}) \subset \Nef(Y^{\mathrm{gen}}) \cap \Bbig(Y^{\mathrm{gen}}).
	\]
By Theorem~\ref{thm:main2}, the group $W_Y^{sm}$ acts on $\Nef(Y^{\mathrm{gen}}) \cap \Bbig(Y^{\mathrm{gen}})$ with fundamental domain $\Nef(Y) \cap \Bbig(Y)$. By choosing $D$ sufficiently general, we may assume that 
	$\Theta^{0}(D)$ lies in the interior of a $W_Y^{sm}$-chamber.
Theorem~\ref{thm:main2} and Theorem~\ref{thm:MMP} show that there is a unique element $w(\Theta) \in W_Y^{sm}$ such that 
	\[
	w(\Theta) \circ \Theta^0(D) \in A(Y).
	\]
Note that $w(\Theta)$ is independent of the sufficiently general $D$. 
	By Theorem~\ref{thm:MMP}, there is a composite of relative flops 
	\[
	\Phi_\Theta \colon \cY \dashrightarrow \cY_\Theta
	\]
	such that the restriction of $\Phi_\Theta$ to central fibers extends to an isomorphism
	$\bar{\phi}_\Theta \colon Y \cong Y_\Theta$, and
	\[
	\Phi^0_\Theta \in W_Y^{sm}, \qquad \Phi_\Theta^0\circ\Theta^0(D) \in A(Y).
	\]
	By uniqueness, we must have $\Phi_\Theta^0=w(\Theta)$. 
%
	The composite $\Psi_\Theta=\Phi_\Theta\circ\Theta$:
\[
	\begin{tikzcd}
	\cY \ar[dr] \ar[r, dashed, "\Theta"] &\cY \ar[d]\ar[r, dashed, "\Phi_\Theta"] &\cY_\Theta \ar[dl]\\
	&U&
	\end{tikzcd}
\]
	is a birational map between Calabi--Yau fiber spaces sending a relatively ample class to a relatively ample class, so is an isomorphism. 
	It follows that $\Theta=\Phi_\Theta^{-1} \circ \Psi_\Theta$
	is defined at the generic point of the central fiber and we denote 
	$\theta \colon Y \dashrightarrow Y$
	the restriction to central fibers. 
	We have that:
	\begin{itemize}
		\item The restriction of the isomorphism $\Psi_\Theta$ to central fibers is an isomorphism $\psi_\Theta \colon Y \rightarrow Y_\Theta$.
		\item The restriction of the birational map $\Phi_\Theta$ to central fibers
		extends to an isomorphism $\bar{\phi}_\Theta \colon Y \rightarrow Y_\Theta$ by Theorem~\ref{thm:MMP}.
	\end{itemize}
	It follows that $\theta=\bar{\phi}_\Theta^{-1} \circ \psi_\Theta$ extends to an automorphism  
	$\bar{\theta} \colon Y \rightarrow Y$.
Its action on $H^2(Y)$ is given by
	\begin{equation}\label{eq:inducedisooncentral1}
	\bar{\theta}_*=\bar{\phi}_{\Theta,*}^{-1} \circ \psi_{\Theta,*}=(\bar{\phi}_{\Theta,*}^{-1} \circ\Phi_\Theta^0) \circ  \Theta^0=w(\Theta) \circ \Theta^0
	\end{equation}
where the last equality follows from Proposition~\ref{prop:geometricrealization}.
	We define 
	\[
	i(\Theta)=(w(\Theta)^{-1}, \bar{\theta}).
	\]
	We show that $i$ is a group homomorphism. Let $\Theta_1, \Theta_2 \in \Aut(\cY_\eta)$. 
We have that 
\begin{align*}
	i(\Theta_1) i(\Theta_2)
	&=(w(\Theta_1)^{-1},\overline{\theta_1}) \cdot (w(\Theta_2)^{-1},\overline{\theta_2})\\
	&=(W^{-1},  \bar{\theta}_1 \circ  \bar{\theta}_2),
	\end{align*}
	where 
	\[
	W=\bar{\theta}_{1,*} \circ w(\Theta_2)\circ \bar{\theta}_{1, *}^{-1} \circ w(\Theta_1)=w(\Theta_1) \circ (\Theta_1)^0 \circ w(\Theta_2) \circ ((\Theta_1)^0)^{-1} 
	\]
	using Equality~\eqref{eq:inducedisooncentral1}. We compute that 
	\[
	W \circ (\Theta_1 \circ \Theta_2)^0(D) \in A(Y)
	\]
	so that $W=w(\Theta_1 \circ \Theta_2)$. It follows that 
	\[
	i(\Theta_1 \circ \Theta_2)=	i(\Theta_1) i(\Theta_2)
	\]
	as required.
We now show that $i$ is injective. 
Suppose that  
\[
i(\Theta)=(1, \id_Y).
\]
Since $w(\Theta)=1$, we may take $\Phi_\Theta=\id_\cY$ in the above proof, so it follows that $\Theta \in \Aut(\cY/U)$. Since $\bar{\theta}=\id_Y$, $\Theta$ is an automorphism of $\cY/U$ restricting to the identity on central fibers. Since $H^0(Y, T_Y)=0$, it follows that $\Theta$ is the identity (see the argument in Lemma~\ref{lem:h2enough}) and $i$ is injective.\\
The proof of the second part is very similar:
An element $\Theta \in \PsAut(\cY_\eta)$ induces a birational automorphism $\Theta \colon \cY \dashrightarrow \cY$ over $S$. Since $\cY \rightarrow S$ is a smooth Calabi--Yau fiber space, $\Theta$ is an SQM over $S$, so $\Theta \in \PsAut(\cY/S)$. We again consider the restriction of $\Theta$ to the inverse image of a small open ball $U \subset S$ (in the analytic topology) centered at $0$, and denote 
\[
\Theta^0 \colon H^2(Y) \rightarrow H^2(Y)
\]
the isomorphism induced by parallel transport, followed by strict transform, followed by parallel transport. 
Let $D \in A(Y)$. Since $\Theta$ restricts to a pseudo-automorphism on a very general fiber, we have 
\[
\Theta^0(D) \in \Mov^\circ(Y^{\mathrm{gen}}) \subset \overline{\Mov}(Y^{\mathrm{gen}}) \cap \Bbig(Y^{\mathrm{gen}}).
\]
By Theorem~\ref{thm:main2}, the group $W_Y^{big}$ acts on $\overline{\Mov}(Y^{\mathrm{gen}}) \cap \Bbig(Y^{\mathrm{gen}})$ with fundamental domain $\overline{\Mov}(Y) \cap \Bbig(Y)$. By choosing $D$ sufficiently general, we may assume that $\Theta^0(D)$ lies in the interior of a $W_Y^{big}$ chamber. 
Theorem~\ref{thm:main2} and Theorem~\ref{thm:MMP} show that there is a unique element $w(\Theta) \in W_Y^{big}$ such that 
\[
w(\Theta)\Theta^0(D) \in \Mov^\circ(Y).
\]
Note that $w(\Theta)$ is independent of the sufficiently general $D$. 
By Theorem~\ref{thm:MMP}, there is a composite of relative flops 
\[
\Phi_\Theta \colon \cY \dashrightarrow \cY_\Theta
\]
such that the restriction of $\Phi_\Theta$ to central fibers extends to an isomorphism
$\bar{\phi}_\Theta \colon Y \cong Y_\Theta$, and
\[
\Phi^0_\Theta \in W_Y^{big}, \qquad \Phi_\Theta^0\circ\Theta^0(D) \in \Mov^\circ(Y).
\]
By uniqueness, we must have $\Phi_\Theta^0=w(\Theta)$. 
Using the decomposition of the movable cone into nef cones~\eqref{MovingConeDecomposition} and using that $D$ is sufficiently general, there is an SQM $\alpha \colon Y_\Theta \dashrightarrow Y'$ such that 
\[
\alpha_* \circ w(\Theta)\circ \Theta^0(D)\in A(Y').
\]
The SQM $\alpha$ deforms to a SQM $\tilde{\alpha} \colon \cY_\Theta \dashrightarrow \cY'$, and as before, the composite $\Psi_\Theta=\tilde{\alpha} \circ\Phi_\Theta\circ\Theta$:
\[
\begin{tikzcd}
	\cY \ar[dr] \ar[r, dashed, "\Theta"] &\cY \ar[d]\ar[r, dashed, "\Phi_\Theta"] &\cY_\Theta \ar[dl] \ar[r, dashed, "\tilde{\alpha}"] &\cY' \ar[dll]\\
	&U&&
\end{tikzcd}
\] 
sends a relatively ample class to a relatively ample class and is therefore an isomorphism.	It follows that
\[
\Theta=\Phi_\Theta^{-1} \circ\tilde{\alpha}^{-1} \circ \Psi_\Theta
\]
is defined at the generic point of the central fiber and we denote 
$
\theta \colon Y \dashrightarrow Y
$
the restriction to central fibers. 
We have that:
\begin{itemize}
	\item The restriction of the isomorphism $\Psi_\Theta$ to central fibers is an isomorphism $\psi_\Theta \colon Y \rightarrow Y'$.
	\item The restriction of the SQM $\tilde{\alpha}$ to central fibers
	is the SQM $\alpha \colon Y_\Theta \dashrightarrow Y'$. 
	\item The restriction of the birational map $\Phi_\Theta$ to central fibers
	extends to an isomorphism $\bar{\phi}_\Theta \colon Y \rightarrow Y_\Theta$ by Theorem~\ref{thm:MMP}.
\end{itemize}
This shows that $\theta=\bar{\phi}_\Theta^{-1} \circ \alpha^{-1} \circ \psi_\Theta$ extends to an SQM 
$\bar{\theta} \colon Y \dashrightarrow Y$, 
and its action on $H^2(Y)$ is given as before by
\begin{equation}\label{eq:inducedisooncentral}
	\bar{\theta}_*=w(\Theta) \circ \Theta^0.
\end{equation}
We define 
\[
i(\Theta)=(w(\Theta)^{-1}, \bar{\theta}),
\]
and the same argument as before shows that $i$ is a group homomorphism. For injectivity, suppose that 
\[
i(\Theta)=(1, \id_Y).
\]
Since $w(\Theta)=1$, we may take $\Phi_\Theta=\id_\cY$ in the above proof, so it follows that $\Theta \in \PsAut(\cY/U)$ restricts to an SQM on every fiber. Since $\bar{\theta}=\id_Y$,  $\Theta^0$ acts trivially on $H^2(Y)$. In particular, $\Theta^0$ preserves a relatively ample divisor, so that $\Theta \in \Aut(\cY/U)$. The argument from before now shows that $\Theta=\id_\cY$.  
	\end{proof}

\begin{theorem}\label{thm:MCClocal}
Let	$\pi \colon \cY \rightarrow \Def(Y)$ be the Kuranishi family of a Calabi--Yau threefold $Y$. The nef cone conjecture (resp. movable cone conjecture) holds for $Y$ if and only if it holds for $Y^{\mathrm{gen}}$.
\end{theorem}
\begin{proof}
We first consider the nef case.
We denote $\tilde{\pi} \colon \tilde{\cY} \rightarrow S$ an algebraization of $\pi$. Since $\Def(Y)$ is smooth, we may assume that $S$ is a smooth variety. It is clear that the nef cone conjecture holds for a very general fiber of $\pi$ if and only if it holds for a very general fiber of $\tilde{\pi}$, so we will denote both by $Y^{\mathrm{gen}}$.\\
	Suppose first that the nef cone conjecture holds for $Y^{\mathrm{gen}}$. This means that $\oAut(Y^{\mathrm{gen}})$ acts on $\Nef^+(Y^{\mathrm{gen}})$ with RPFD, so that
	\[
	(\oAut(Y^{\mathrm{gen}}), \Nef^+(Y^{\mathrm{gen}}))
	\]
	 is of polyhedral type.
	In particular, this implies that $\overline{\Aut}(Y^{\mathrm{gen}})$ is finitely generated (see \cite[Corollary 4.15]{Looijenga}). By Lemma~\ref{lem:kernelfinite}, $\Aut(Y^{\mathrm{gen}})$ is also finitely generated. After a base change by a generically finite morphism, Proposition~\ref{prop:genericvsgeneral} shows that we may assume $\Aut(Y^{\mathrm{gen}})=\Aut(\tilde{\cY}_\eta)$.
	By Proposition~\ref{prop:genericautomorphismgroupsubgroup}, we have 
	\[
	\Aut(Y^{\mathrm{gen}}) =\Aut(\tilde{\cY}_\eta) \subset W_Y^{sm} \rtimes \Aut(Y)
	\]
and hence, using Lemma~\ref{lem:kerneltrivial},
\[
i \colon \overline{\Aut}(Y^{\mathrm{gen}}) \hookrightarrow \overline{W_Y^{sm} \rtimes \Aut(Y)}=W_Y^{sm} \rtimes \overline{\Aut}(Y)
\]
$\Theta \in {\Aut}(Y^{\mathrm{gen}})$ acts on $H^2(Y, \ZZ)$ via the identification given by parallel transport, which is given by $\Theta^0 \in \Aut(H^2(Y, \ZZ))$.
It follows that $i$ is an inclusion of subgroups of $\Aut(H^2(Y))$: writing $i(\Theta)=(w(\Theta)^{-1}, \bar{\theta})$ as in the proof of Proposition~\ref{prop:genericautomorphismgroupsubgroup}, Equality~\eqref{eq:inducedisooncentral1} shows that 
\[
i(\Theta)_*=w(\Theta)^{-1} \circ \bar{\theta}_*=\Theta^0
\]
Moreover, the group $W_Y^{sm} \rtimes \overline{\Aut}(Y)$ preserves $A(Y^{\mathrm{gen}})$. Indeed, this is true for the action of $W_Y^{sm}$ by Corollary~\ref{thm:weylpreservescones}. For the action of $\Aut(Y)$, recall that $\theta \in \Aut(Y)$ deforms uniquely to an automorphism
\[
\begin{tikzcd}
	\cY \ar[r, "{\Theta}"] \ar[d, ]& \cY \ar[d]\\
	\Def(Y) \ar[r, "\bar{\theta}"] &\Def(Y)
\end{tikzcd}
\]
between universal families, and therefore induces an isomorphism $\Theta_s \colon Y_s \rightarrow Y_{\bar{\theta}(s)}$ for all $s \in \Def(Y)$. If $s$ is very general, so is ${\bar{\theta}}(s)$, so that $\Theta^0$ preserves $\Nef(Y^{\mathrm{gen}})$, and its interior $A(Y^{\mathrm{gen}})$. \\
By Remark~\ref{rem:subgroup}, it now follows that the pair 
\[
(W_Y^{sm} \rtimes \overline{\Aut}(Y), \Nef^+(Y^{\mathrm{gen}}))
\]
 is of polyhedral type as well. 
 By Theorem~\ref{thm:RPFDcriterion}, the cone 
 \[
 \Pi_\xi=\{x \in \Nef^+(Y^{\mathrm{gen}}) \mid \xi \cdot (\gamma(x)) \geq \xi \cdot x \; \text{for all} \; \gamma \in W_Y^{sm} \rtimes \overline{\Aut}(Y)\}
 \]
 is an RPFD for the action of $W_Y^{sm} \rtimes \overline{\Aut}(Y)$ on $\Nef^+(Y^{\mathrm{gen}})$, for $\xi$ a generic rational point in $(\Nef(Y^{\mathrm{gen}})^*)^\circ$. 
 We want to show that 
$
 \Pi_\xi \subset \Nef(Y).
 $
 	Fix a rational class $H \in A(Y) \subset A(Y^{\mathrm{gen}})$. By the Hard Lefschetz theorem, multiplication by $H$ induces an isomorphism
 	\[
 	H^2(Y,\mathbb R)\longrightarrow H^4(Y,\mathbb R)
 	\qquad
 	H'\longmapsto H\cdot H'.
 	\]
 	Consequently, as $H'$ varies in the ample cone $A(Y)$, the classes $H\cdot H'$ form an open subset $U \subset H^4(Y)$. We have that  $H \cdot H' \cdot D \geq 0$ for 
 	\[
 	D \in \Nef(Y^{\mathrm{gen}}) \subset \overline{\Eff}(Y^{\mathrm{gen}}),
 	\]
 	so we must have that 
 	$U \subset  \Nef(Y^{\mathrm{gen}})^*$,
 	where we have identified $H^4(Y)$ with the dual of $H^2(Y)$. So $\xi=H\cdot H'$
 	is a generic rational point in $(\Nef(Y^{\mathrm{gen}})^*)^\circ$ if $H'$ is generic.
 Applying the inequality in the definition of $\Pi_{\xi}$ to $\gamma=\sigma_E$ gives
 \[
 (x \cdot \ell_E)H \cdot H' \cdot E \geq 0.
 \]
 Since $H \cdot H' \cdot E>0$ ($E$ is effective on $Y$), this implies that $x \cdot \ell_E \geq 0$ for every $E \in \Delta_Y^{sm}$, so that 
 \[
 \Pi_\xi \subset \bar{\cC}^{sm} \cap \Nef(Y^{\mathrm{gen}}) = \Nef(Y),
 \]
 where we have used \eqref{eq:nef}. 
 Since $\Pi_\xi$ is rational polyhedral, it is spanned by finitely many rational points in $\Nef(Y)$, so we in fact have $\Pi_\xi \subset \Nef^+(Y)$.
 Let now $x \in A(Y)$. There exists $(w,\theta) \in  W_Y^{sm} \rtimes \overline{\Aut}(Y)$ such that 
 \begin{equation}\label{eq11}
 w(\theta(x)) \in \Pi_\xi \subset \Nef(Y).
 \end{equation}
 We have that 
 \[
 \theta(x) \in A(Y) \subset \Nef(Y) \cap \Bbig(Y^{\mathrm{gen}}).
 \]
Since the action of $W_Y^{sm}$ preserves $\Bbig(Y^{\mathrm{gen}})$, Equality~\eqref{eq11} shows that 
\[
w(\theta(x)) \in \Nef(Y) \cap \Bbig(Y^\mathrm{gen})=\Nef(Y) \cap \Bbig(Y).
\]
It follows that $\theta(x) \in A(Y)$ and $w(\theta(x)) \in \Nef(Y) \cap \Bbig(Y)$, so by Theorem~\ref{thm:main2}, we must have $w=1$. This shows that 
	\[
	\overline{\Aut}(Y) \cdot \Pi_\xi \supset A(Y).
	\]
	It follows that $(\oAut(Y), \Nef^+(Y))$ is of polyhedral type, so the action of $\oAut(Y)$ on $\Nef^+(Y)$ admits an RPFD, i.e. the nef cone conjecture holds for $Y$. \\
	
Conversely, suppose that the nef cone conjecture holds for $Y$, 
so that $\oAut(Y)$ acts on $\Nef^+(Y)$ with RPFD $\Pi$. Let $x \in A(Y^{\mathrm{gen}}) \subset \Nef(Y^{\mathrm{gen}}) \cap \Bbig(Y^{\mathrm{gen}})$. By Theorem~\ref{thm:MMP}, there exists $w \in W_Y^{sm}$ such that $w(x) \in \Nef(Y) \cap \Bbig(Y) \subset \Nef^+(Y)$. By assumption, there is then $\theta \in \oAut(Y)$ such that $y:=\theta(w(x)) \in \Pi$. So 
\[
%
(W_Y^{sm} \rtimes \overline{\Aut}(Y)) \cdot \Pi \supset A(Y^{\mathrm{gen}}),
\]
and it follows that 
\[
(W_Y^{sm} \rtimes \overline{\Aut}(Y), \Nef^+(Y^{\mathrm{gen}}))
\] is of polyhedral type. 
Let $K$ be the kernel of the action of $W_Y^{sm} \rtimes \Aut(Y)$ on $\Def(Y)$. By Proposition~\ref{prop:actsfinitegroup}, $K$ has finite index, and therefore so has $\bar{K} \subset  W_Y^{sm} \rtimes \overline{\Aut}(Y)$. Using Theorem~\ref{thm:RPFDfiniteindex}, we see that 
\[
(\bar{K}, \Nef^+(Y^{\mathrm{gen}}))
\] 
is of polyhedral type as well. By Theorem~\ref{thm:biggroupacts}, $W_Y^{sm} \rtimes \Aut(Y)$ acts on the universal family $\cY \rightarrow \Def(Y)$ via birational automorphisms which restrict to isomorphisms on a very general fiber. We therefore have a restriction map 
$r \colon K \rightarrow \Aut(Y^{\mathrm{gen}})$. We claim that $r$ is injective. To see this, let $(w, \theta) \in K$ and let
				\[
\begin{tikzcd}
	\cY \ar[d] \ar[r, "\Theta"]&\cY \ar[d] \ar[r, "{\Phi_w}", dashed] &\cY \ar[d, ]\\
	(\Def(Y), 0) \ar[r, "\bar{\theta}"] &(\Def(Y), 0) \ar[r, "\bar{\phi}_w"] &(\Def(Y), 0) 
\end{tikzcd}
\]
be the induced birational automorphism of the Kuranishi family. By assumption, we have that 
\[
\bar{\phi}_w \circ \bar{\theta}=\id_{\Def(Y)}.
\]
Via parallel transport, the action of $r(w, \theta)$ on $H^2(Y)$ is given by 
\begin{equation}\label{eq:kernelcompatible}
r(w, \theta)_*=	w \circ \theta_*.
\end{equation}
If $r(w, \theta)=1$, we have $w^{-1} = \theta_*$.
Since $\theta$ is an automorphism, it preserves $\Nef(Y)$, so the same must be true for $w$. By Theorem~\ref{thm:main2}, this forces $w=1$. It now follows that $\Theta$ restricts to the identity on a general fiber, so we must have $\Theta=\id_\cY$ and therefore also $\theta=\id_Y$, so $r$ is injective. By \eqref{eq:kernelcompatible}, the induced inclusion $\bar{K} \subset \overline{\Aut}(Y^{\mathrm{gen}})$ is compatible with the actions on $H^2(Y)$, so Remark~\ref{rem:subgroup} shows that
\[
(\oAut(Y^{\mathrm{gen}}), \Nef^+(Y^{\mathrm{gen}}))
\]
is of polyhedral type, so the nef cone conjecture holds for $Y^{\mathrm{gen}}$.\\
The proof for the movable cone conjecture proceeds along exactly the same lines, replacing $\Nef$ by $\overline{\Mov}$, $A(Y)$ by $\text{Mov}(Y)^\circ$, $W_Y^{sm}$ by $W_Y^{big}$, and ``automorphism" by ``pseudo-automorphism". We only spell out three points where the argument differs in a nontrivial way: 
\begin{itemize}
	\item In the construction of $\xi$ one uses that 
	\[
	\overline{\Mov}(Y^{\mathrm{gen}}) \subset \overline{\Eff(Y^{\mathrm{gen}})}.
	\]
\item We apply the inequality in the definition of $\Pi_\xi$ to every $E \in \Delta_Y^{big}$, which gives 
\[
\Pi_\xi \subset \bar{\cC}^{big} \cap \overline{\Mov}(Y^{\mathrm{gen}})=\overline{\Mov}(Y).
\]
\item 
	In the proof of the converse direction, we use the map 
	\[
	r \colon K \rightarrow \PsAut(Y^{\mathrm{gen}}),
	\]
	and injectivity is proved using the movable fundamental domain. 
	\end{itemize}
\end{proof}
We are now ready to prove Theorem~\ref{thm:mainA}.
\begin{proof}[Proof of Theorem~\ref{thm:mainA}]
Let $\pi \colon \cY \rightarrow \Def(Y)$ be the Kuranishi family of $Y$. After shrinking $\Def(Y)$, by openness of versality, we may assume that $\pi$ is a Kuranishi family for each of its fibers. Given $p, q \in \Def(Y)$, we may then choose a fiber $Y_t$ which is a very general deformation of both $Y_p$ and $Y_q$, so if the nef (resp. movable) cone conjecture holds for $Y_p$, it also holds for $Y_q$, by applying Theorem~\ref{thm:MCClocal} twice. \\
If $\pi \colon \cY \rightarrow S$ is a family of smooth Calabi--Yau threefolds over a connected base $S$, then $\pi$ is locally pulled back from the Kuranishi family, so the validity of the cone conjecture is locally constant on $S$. Since $S$ is connected, the validity of the cone conjecture is invariant under deformation. 
\end{proof}

\section{Examples and applications}\label{sec:examples}
We conclude this paper with a series of examples of the Wilson Weyl group and use Theorem~\ref{thm:mainA} to prove a new case of the cone conjecture.
\subsection{Examples of the Wilson Weyl group}
\begin{example}
Let $\pi  \colon Y \rightarrow S$ be a Calabi--Yau threefold with an elliptic fibration, let $C$ be a component of the discriminant locus $\Delta$, and suppose that $g(C^\nu)>0$, where $C^\nu \rightarrow C$ denotes the normalization.
Let $D=\pi^{-1}(C)_{red}$. Suppose that $D$ is reducible with components $E_1, \dots, E_m$, and denote 
\[
p \colon D \rightarrow C, \qquad p_i \colon E_i \rightarrow C
\]
the induced morphisms. The general fiber $F$ of $p$ is a singular fiber of Kodaira type 
\[
I_n (n>1), I_n^*, IV^*, IV, III^*, III, II^*,
\]
corresponding to the affine Dynkin diagrams 
\[
\tilde{A}_{n-1}, \tilde{D}_{n+4}, \tilde{E}_6,  \tilde{A}_2, \tilde{E}_7, \tilde{A}_1, \tilde{E}_8
\]
respectively. 
We recall that every irreducible component of such a singular fiber is a $\PP^1$. 
Let $C^0 \subset C^{reg}$ be an open subset disjoint from other components of $\Delta$, and over which $p \colon D \rightarrow C$ is a locally trivial fiber bundle (in the analytic topology). Then $\pi_1(C^0)$ acts on the components of $F$, and the orbits of this action biject with the irreducible components of $D$. 
Since $g(C^\nu)>0$, the proof of Proposition~\ref{prop:everysurfacecontracts} shows that there is an SQM $\alpha \colon Y \dashrightarrow Z$ such that the strict transform $E_Z$ of $E_i$ can be contracted via a primitive contraction $q \colon Z \rightarrow \bar{Z}$ along its ruling to a curve $C_Z$, with $C_Z$ smooth by \cite{WilsonCone}. Since the SQM only flops curves contained in fibers of $p_i \colon E_i \rightarrow C$, there is also an induced morphism $E_Z \rightarrow C$, with Stein factorization 
\[
E_Z \xrightarrow{q|_{E_Z}}  C_Z \rightarrow C.
\]
In particular, $g(C_Z)>0$. 
 It follows that $E_i$ is an SQM-conic fibration surface over a curve of positive genus, so that $[E_i] \in \Delta_Y^{big}$. Computing the intersection numbers $E_i \cdot \ell_{E_j}$, we see that the subgroup $G \subset W_Y^{big}$ generated by the reflections $\sigma_{E_i}$ is isomorphic to the Weyl group of the folded affine Dynkin diagram obtained from the diagram of $F$ according to the monodromy action.
These folded Weyl groups have appeared in the work of Szendr\H oi~\cite[Section 5]{SzendroiEnhanced}. See also \cite{Markman} for similar results in the setting of Hyperk\"ahler manifolds.
We summarize all of this as follows:
	\begin{proposition}
		Let $ \pi \colon Y \rightarrow S$ be a Calabi--Yau threefold with an elliptic fibration. Let $\Delta=C_1+\dots +C_n$ be the decomposition of the discriminant locus of $\pi$ into irreducible components.
		Let $X$ be the subset of those components satisfying:
		\begin{itemize}
			\item $g(C_i^\nu)>0$.
			\item $\pi^{-1}(C_i)$ is reducible.
		\end{itemize}
		Then $W^{big}_Y$ contains a subgroup isomorphic to $\prod_{C_i \in X} W_{C_i}$ where $W_{C_i}$ is the Weyl group corresponding to the (folded) Dynkin diagram given by the irreducible components of $\pi^{-1}(C_i)$. 
	\end{proposition}
	\begin{proof}
		We have already verified this for the case of a single curve in the discussion above. We check that generators of Weyl groups coming from two different curves $C_i$ and $C_j$ commute. This is easy to see: suppose that $E_i \in \Delta_Y^{big}$ is an irreducible component of $\pi^{-1}(C_i)$ and $E_j \in \Delta_Y^{big}$ is an irreducible component of $\pi^{-1}(C_j)$. Then 
		\[
		E_i \cdot \ell_{E_j}=0,
		\]
		since $E_i$ is disjoint from a general fiber of $E_j \rightarrow C_j$.
		Similarly, we have $E_j \cdot \ell_{E_i}=0$. It follows that the corresponding elements $\sigma_{E_i}$ and $\sigma_{E_j}$ of $W_Y^{big}$ commute. Moreover, $W_{C_i}$ intersects the group generated by all other $W_{C_j}$ trivially: $W_{C_i}$ acts faithfully on the vector space 
		\[
		V_i=\Span\{[E_i] \mid E_i \subset \pi^{-1}(C_i)\},
		\]
		whereas $W_{C_j}$ with $j \neq i$ fixes $V_i$ pointwise. It follows that 
		\[
		\prod_{C_i \in X} W_{C_i}
		\]
		is a subgroup of $W_Y^{big}$.
	\end{proof}
\end{example}
The same arguments show the following:
\begin{example}
	Suppose that $\pi \colon Y \rightarrow \bar{Y}$ is a crepant resolution of a Calabi--Yau threefold with canonical singularities. Assume also that the singular locus of $\bar{Y}$ is a smooth curve $C$ of positive genus along which the transverse singularities are cDV. Then the irreducible components $E_i$ of the exceptional divisor give elements $E_i \in \Delta_Y^{big}$. 
	Computing intersection numbers, the group $G \subset W_Y^{big}$ generated by the $\sigma_{E_i}$ is isomorphic to the Weyl group of the corresponding ADE singularity, possibly folded by the corresponding monodromy action on the ADE diagram.
\end{example}
\begin{example}[Example of a Calabi--Yau threefold with $W_Y^{sm}$ infinite]
 Let $B$ and $B'$ be rational elliptic surfaces and let $Y=B \times_{\PP^1} B'$. If the discriminant loci of $B$ and $B'$ are disjoint, $Y$ is a smooth elliptic Calabi--Yau threefold. These threefolds were studied by Schoen~\cite{Schoen}. If either $B$ or $B'$ has a reducible fiber $F$, then $W_Y^{sm}$ contains the Weyl group of the affine Dynkin diagram corresponding to $F$ and is infinite.
\end{example}

\subsection{The results of Cantat--Oguiso}\label{sec:CO}
Let $Y$ be a generic $(2, \dots, 2)$-divisor in $(\PP^1)^n$. 
Cantat--Oguiso~ \cite{CantatOguiso} have studied the birational geometry of $Y$, and have verified the movable cone conjecture for $Y$ by relating $\Nef(Y)$ and $\Mov^e(Y)$ via an explicit Coxeter group. 
In the case $n=4$ (so that $Y$ is a Calabi--Yau threefold), we use Theorem~\ref{thm:mainA} to give a new proof of the movable cone conjecture for an arbitrary smooth $(2,2,2,2)$-divisor in $(\PP^1)^4$ (not necessarily very general). Moreover, we show that the Coxeter group appearing in their work is nothing other than the big Weyl group of a certain degeneration of $Y$. \\
We briefly review the setup.
Let $Y$ be a smooth $(2,2,2,2)$-divisor in $V=(\PP^1)^4$. 
Let $\pi_i$, $1 \leq i\leq 4$ be the projections onto $(\PP^1)^3$.
For any $1 \leq i \leq 4$, we may write the equation of $Y$ as 
\[
F_{i,1}x_0^2+F_{i,2}x_0x_1+F_{i,3}x_1^2=0
\]
where $F_{i,j}$ are homogeneous quadratic in the remaining coordinates. 
Away from the locus 
\[
B_i=\{F_{i,1}=F_{i,2}=F_{i,3}=0\} \subset (\PP^1)^3,
\]
$\pi_i$ exhibits $Y$ as a double cover.
If $x \in B_i$, then $\pi_i^{-1}(x)$ is a $\PP^1$. We let $\ell_i$ be the class of this $\PP^1$-fiber. 
If $Y$ is very general, then $B_i$ consists of $48$ points, so that $Y$ contains $48$ smooth rational curves of class $\ell_i$.
The birational involution $\iota_i \colon Y \dashrightarrow Y$ induced by the double cover $\pi_i$ is the simultaneous flop of these $48$ curves. 
If $Y$ is smooth, by the Lefschetz hyperplane theorem, $i^* \colon \Pic(V) \rightarrow \Pic(Y)$ is an isomorphism, and by \cite[Theorem 3.1]{CantatOguiso}, we have $i^*(A(V))=A(Y)$.
In particular this shows that $\Nef(Y)$ is rational polyhedral and is bounded by the four hyperplanes $\ell_i^\perp$. 
\begin{theorem}[Cantat--Oguiso~\cite{CantatOguiso}]
	Let $Y^{\mathrm{gen}}$ be a generic $(2,2,2,2)$-divisor in $(\PP^1)^4$.
\begin{enumerate}
	\item 
$
\Bir(Y^{\mathrm{gen}})=\langle i_1, \dots i_4 \rangle \cong \ZZ_2 \star \ZZ_2 \star \ZZ_2 \star \ZZ_2.
$
\item $\Nef(Y^{\mathrm{gen}})$ is a fundamental domain for the action of $\PsAut(Y^{\mathrm{gen}})$ on $\Mov^e(Y^{\mathrm{gen}})$. 
\end{enumerate}
In particular, the movable cone conjecture holds for $Y^{\mathrm{gen}}$.
\end{theorem}
We now explain how to recover the cone-conjecture statement and give a geometric interpretation of the Coxeter group and its fundamental chamber on the big movable cone as a corollary of Theorem~\ref{thm:mainA}. To do this, we consider a special degeneration of $Y^{\mathrm{gen}}$. 
Note first that if $F_{i,1}, F_{i,2}, F_{i,3}$ are very general subject to the condition that they are linearly dependent, then the locus $B_i$ is a curve of genus $g=25$, and $E_i:=\pi_i^{-1}(B_i)$ is a surface ruled over $B_i$, so that $E_i \in \Delta_Y^{big}$, with $\ell_i=\ell_{E_i}$ the conic fiber class. By Wilson's results \cite{WilsonCone}, $E_i$ does not deform to $Y^{\mathrm{gen}}$, but $2g-2=48$ fibers do. These are precisely the $48$ isolated $\PP^1$'s on $Y^{\mathrm{gen}}$. 
We  claim that there exists a smooth degeneration $Y_0$ on which all the $E_i$, $1 \leq i \leq 4$ are present simultaneously. 
To see this, let $V \subset H^0(\PP^1, \cO(2))$ be the subspace spanned by the monomials $x_0^2$ and $x_1^2$, so that $L:=V^{\boxtimes 4} \subset H^0((\PP^1)^4, \cO(2,2,2,2))$ is spanned by monomials of the form $x_i^2y_j^2z_k^2w_l^2$, where $i,j,k,l$ are each $0$ or $1$. $L$ is base-point-free, so the general member $Y_0$ is smooth by Bertini's theorem. If we now write the equation of $Y_0$ as
\[
F_{i,1}x_0^2+F_{i,2}x_0x_1+F_{i,3}x_1^2=0
\]
then $F_{i,2}=0$ by construction, so that $F_{i,1}, F_{i,2}, F_{i,3}$ are linearly dependent and $E_i$ is present on $Y_0$. 
This shows that the four hyperplanes $\ell_i^\perp$ also bound $\overline{\Mov}(Y_0)$, so that 
\[
\Nef(Y_0) \subset \overline{\Mov}(Y_0) \subset \bigcap_{i=1}^4\{D \mid D \cdot \ell_i \geq 0\} =\Nef(Y_0)
\]
 In particular, $\overline{\Mov}(Y_0)=\Nef(Y_0)$ is rational polyhedral, and the movable cone conjecture holds for $Y_0$.
Using Theorem~\ref{thm:mainA}, we obtain:
\begin{theorem}
	The movable cone conjecture holds for any smooth $(2,2,2,2)$-divisor in $(\PP^1)^4$. 
	\end{theorem}
$W_{Y_0}^{big}$ is generated by the reflections $\sigma_{E_i}$ for $1 \leq i \leq 4$. 
It is not hard to show using the projection formula on $\pi_i$ that $E_i \cdot \ell_j=2$ for $i \neq j$, so by definition (see \eqref{eq:defnofm}), there is no relation between $\sigma_{E_i}$ and $\sigma_{E_j}$ and 
\[
W_{Y_0}^{big}=\ZZ_2 \star \ZZ_2 \star \ZZ_2 \star \ZZ_2
\]
Theorem~\ref{thm:main2} now shows that $W_{Y_0}^{big}$ acts on $\overline{\Mov}(Y^{\mathrm{gen}}) \cap \Bbig(Y^{\mathrm{gen}})$ with fundamental domain $\overline{\Mov}(Y_0) \cap \Bbig(Y_0)$. 
Since $\overline{\Mov}(Y_0)=\Nef(Y_0)=\Nef(Y^{\mathrm{gen}})$, we obtain that $W_{Y_0}^{big}$ acts on
$\overline{\Mov}(Y^{\mathrm{gen}}) \cap \Bbig(Y^{\mathrm{gen}})$ with fundamental domain $\Nef(Y^{\mathrm{gen}}) \cap \Bbig(Y^{\mathrm{gen}})$. 
This explains the appearance of the Coxeter group in the result of Cantat--Oguiso in dimension 3.
\subsection{Double covers of klt log Calabi--Yau pairs}
Recall that a pair $(X, \Delta)$ with $X$ a smooth projective variety and $\Delta$ an effective $\QQ$-divisor is a klt log Calabi--Yau pair if $(X, \Delta)$ is klt and $K_X+\Delta \sim_\QQ 0$. We allow $\Delta=\emptyset$, in which case $X$ is a Calabi--Yau variety.
The Morrison cone conjecture was generalized to klt log Calabi--Yau pairs by Totaro~\cite{Totaro}:
\begin{conjecture}\label{conj:MCClogCY}[Morrison--Totaro cone conjecture]
	Let $(X, \Delta)$ be a klt log Calabi--Yau pair. 
	\begin{itemize}
		\item The action of $\Aut(X, \Delta)$ on $\Nef^+(X)$ admits a rational polyhedral fundamental domain.
		\item  The action of $\PsAut(X, \Delta)$ on ${\Mov}^+(X)$ admits a rational polyhedral fundamental domain.
	\end{itemize}
\end{conjecture}
Here, a (pseudo)-automorphism of $(X, \Delta)$ is a (pseudo)-automorphism of $X$ that preserves the support of $\Delta$. 
The conjecture was proven by Totaro in dimension 2 \cite{Totaro}. A preprint of Xu~\cite{Xu} proves the cone conjecture in a wide range of cases, of which we will use the following:
\begin{theorem}\cite[Theorem 2, Corollary 2]{Xu}\label{thm:Xu}
Let $(X, \Delta)$ be a klt log Calabi--Yau pair of dimension $3$, and let $f \colon (X, \Delta_X) \rightarrow (Y, \Delta_Y)$ be a crepant birational morphism between $\QQ$-factorial klt log Calabi--Yau pairs. 
Suppose further that the support of $\Delta_X$ contains the exceptional divisors of $f$. Then the movable cone conjecture holds for $(Y, \Delta_Y)$ if and only if it holds for $(X, \Delta_X)$.
\end{theorem}
The cone conjecture is stated in \cite{Xu} in terms of the cone $\Nef^e$ (respectively $\Mov^e$). As noted in Corollary~\ref{cor:relationMCC}, this implies the cone conjecture for $\Nef^+$ (respectively $\Mov^+$). 
We will use Theorem~\ref{thm:Xu} to prove:
\begin{proposition}\label{prop:doublecover}
Let $X$ be a smooth projective threefold, with $h^1(X, \cO_X)=h^2(X, \cO_X)=0$. Let $D \in |-2K_X|$ be a smooth divisor, and let $\Delta=\tfrac{1}{2}D$, so that $(X, \Delta)$ is a klt log Calabi--Yau pair. 
Let $p \colon Y \rightarrow X$ be the double cover of $X$ branched along $D$.
Then $Y$ is a smooth Calabi--Yau threefold. Suppose further that
	\begin{itemize}
		\item the nef (resp. movable) cone conjecture holds for $(X, \Delta)$. 
		\item $\rho(X)=\rho(Y)$.
	\end{itemize}
	Then the nef (resp. movable) cone conjecture holds for $Y$, and therefore for any smooth Calabi--Yau threefold deformation-equivalent to $Y$. 
\end{proposition}
\begin{proof}
	By the double cover formula, we have that 
	\[
	K_Y=p^*(K_X+\frac{1}{2}D) \sim 0,
	\]
	so that $Y$ is $K$-trivial. Since $D$ is smooth, $Y$ is smooth as well. Moreover, since 
	\[
	p_*\cO_Y=\cO_X \oplus\cO_X(K_X),
	\]
	Serre duality and the vanishing of $h^1(X, \cO_X)$ and $h^2(X, \cO_X)$ imply that 
	\[
	h^1(Y, \cO_Y)=h^2(Y, \cO_Y)=0,
	\]
	so that $Y$ is a smooth Calabi--Yau threefold. 
	Since $p$ is finite, \cite[Example 1.7.4]{Fulton} shows that $p_* \circ p^* \colon H^2(X, \RR) \rightarrow H^2(X, \RR)$ is multiplication by $2$, and since $\rho(X)=\rho(Y)$, it follows that $p^* \colon H^2(X, \RR) \rightarrow H^2(Y, \RR)$ is an isomorphism. 
	Since $p$ is surjective, $L \in \Pic(X)$ is nef if and only if $p^*(L)$ is nef (\cite[1.4.4]{Lazarsfeld}), so it follows that $p^*(\Nef(X))=\Nef(Y)$. 
 $p^*$ is a linear map defined over $\QQ$, so we obtain 
	\[ 
	p^*(\Nef^+(X))=\Nef^+(Y)
	\]
	Let 
	\[
	{G}=\{{g} \in \Aut(Y) \mid p \circ {g}=\bar{g} \circ p \; \text{for some} \; \bar{g} \in \Aut(X, \Delta)\}.
	\]
	Since every automorphism of $X$ preserving $\Delta$ lifts to an automorphism of $Y$, we obtain an exact sequence 
	\[
	1 \rightarrow \langle \iota \rangle \rightarrow {G} \rightarrow \Aut(X, \Delta) \rightarrow 1,
	\]
	where $\iota \in {G}$ is the covering involution. 
	Since $p^*$ is surjective, the action of $\iota$ on $H^2(Y, \RR)$ is trivial, so that conjugation by $p^*$ identifies $\overline{{G}}$ with $\overline{\Aut(X, \Delta)}$.
The action of $\Aut(X, \Delta)$ on $\Nef^+(X)$ admits an RPFD by assumption, hence so does the action of $\overline{G}$ on $\Nef^+(Y)$. By Remark~\ref{rem:subgroup}, the action of the larger group $\overline{\Aut}(Y)$ on $\Nef^+(Y)$ also has an RPFD. \\
The movable case is very similar, only the statement that 
\[
p^*(\overline{\Mov}(X))=\overline{\Mov}(Y)
\]
needs a separate argument. Suppose first that $L$ is movable, so the stable base locus of $L$ has codimension at least two. Since $p$ is surjective, we have an inclusion 
\[
H^0(X, mL) \hookrightarrow H^0(Y, mp^*L)
\]
for every $m>0$, and therefore 
\[
\textrm{Bs}|mp^*L| \subset p^{-1}\textrm{Bs}|mL|.
\]
Since $p$ is finite, we see that if $L$ is movable, then so is $p^*L$.  
Conversely, suppose that $p^*L$ is movable. Let $B \subset X$ be a prime divisor. By assumption, there is a section $s \in H^0(Y, mp^*L)$ which does not vanish entirely along any component of $p^{-1}(B)$, and neither does the section $\iota^*s$. The section
\[
s \otimes \iota^*s \in H^0(Y, 2mp^*L), 
\]
is $\iota$-invariant and therefore the pull back of a section $t \in H^0(X, 2mL)$, and $t$ does not vanish entirely along $B$ by construction. It follows that $L$ is movable, and we obtain the equality by taking closures.\\
 The final statement follows from Theorem~\ref{thm:mainA}.
\end{proof}
\begin{remark}
	The condition that $\rho(X)=\rho(Y)$ in Proposition~\ref{prop:doublecover} is quite restrictive, and must usually be checked on a case-by-case basis.
\end{remark}
Let $Z$ be a smooth projective K3 surface, equipped with an antisymplectic involution $\iota$, and let $E$ be an elliptic curve, with $i \colon E \rightarrow E$ the involution.
Then 
\[
\bar{Y}=Z \times E /(\iota, i)
\]
is a Calabi--Yau threefold with canonical singularities, and the crepant resolution ${Y} \rightarrow \bar{Y}$ is a smooth Calabi--Yau threefold. Calabi--Yau threefolds of this type were studied by Borcea and Voisin~\cite{Borcea}. Under some additional assumptions, we are able to prove the cone conjecture for $Y$:
\begin{theorem}\label{thm:BorceaVoisin}
	Let $Z$ be a $K3$ surface which is a double cover of $\PP^2$ branched along a very general sextic $B$, let $E$ be an arbitrary elliptic curve, and let $Y$ be the associated Borcea--Voisin Calabi--Yau threefold. Then the movable cone conjecture holds for $Y$, and for any smooth Calabi--Yau threefold deformation-equivalent to $Y$.
\end{theorem}
\begin{proof}
	The involution $\iota \colon Z \rightarrow Z$ is just the involution induced by the double cover, so that the fixed locus of $\iota$ is the sextic $B$. It follows that the singular locus of $\bar{Y}$ is the union of the curves $C_i=B \times p_i$, for $p_i$ a $2$-torsion point on $E$. The threefold $\bar{Y}$ has uniform $A_1$-singularities along $C_i$, so the resolution $Y \rightarrow \bar{Y}$ is just the blowup of the $C_i$. 
	Since $B$ is very general, $\rho(Z)=1$. We compute that $\rho(\bar{Y})=2$, and therefore $\rho(Y)=6$. \\
	By construction, $\bar{Y}$ is a double cover of $\PP^2 \times \PP^1$, branched along the union $D$ of five divisors $D_0=B \times \PP^1$ and $D_1, \dots, D_4$ given by 
	$\PP^2 \times p$ where $p$ ranges over the ramification points of $E \rightarrow \PP^1$. Let $\Gamma$ be the $1$-skeleton of $D$, which is the disjoint union of $4$ curves of genus $10$, corresponding to the images of $B \times p$. Let $X$ be the blowup of $\PP^2 \times \PP^1$ along $\Gamma$, and denote $\tilde{D}$ the strict transform of $D$. Note that $\tilde{D}$ is the disjoint union of five smooth divisors $\tilde{D}_0, \dots, \tilde{D}_4$. Note also that $\rho(X)=6$. We have a diagram 
	\[
	\begin{tikzcd}
		{Y} \ar[r] \ar[d]&\bar{Y} \ar[d] \\
		X \ar[r] & \PP^2 \times \PP^1
	\end{tikzcd}
	\]
	which expresses ${Y}$ as a double cover of $X$ branched along $\tilde{D}$. Let $\Delta=\tfrac{1}{2}\tilde{D}$. Since $\tilde{D} \in |-2K_X|$ is a smooth divisor, the pair $(X, \Delta)$ is klt log Calabi--Yau.
	 For
	$0<\epsilon\ll 1$, the pair
	$(X,\Delta+\epsilon \tilde D)$
	is klt, and
	\[
	K_X+\Delta+\epsilon \tilde D
	\sim_{\mathbb Q}\epsilon\tilde D.
	\]
	We run an MMP for this pair. Every negative extremal ray is
	generated by a rational curve contained in one of the disjoint
	components $\tilde D_i$. For $i\geq 1$, such a curve is
	numerically proportional to a line
	$\ell_i\subset \tilde D_i\simeq \mathbb P^2$. Any rational curve
	in $\tilde D_0\simeq B\times\mathbb P^1$ is contained in a fiber
	of $\tilde D_0\to B$, since $g(B)>0$, denote the class of such a
	fiber by $\ell_0$. Thus the MMP contracts the divisors
	$\tilde D_0,\ldots,\tilde D_4$, in some order.\\
	Let $f\colon X\longrightarrow \bar X$
	be the resulting morphism. Then $\bar X$ is $\QQ$-factorial
	and klt, and
	\[
	K_X+\Delta
	\sim_\QQ f^*K_{\bar X},
	\]
	so that $f\colon (X,\Delta)\rightarrow(\bar X,0)$ is
	crepant. By construction, $\bar{X}$ is a klt Calabi--Yau threefold with 
	\[
	\rho(\bar X)=\rho(X)-5=1,
	\] 
	so the cone conjecture is trivially true for $\bar{X}$. 
Using Theorem~\ref{thm:Xu}, the movable cone conjecture holds for $(X, \Delta)$. We have $\rho(X)=\rho(Y)$, and $h^1(X, \cO_X)=h^2(X, \cO_X)=0$, so Proposition~\ref{prop:doublecover} shows that the movable cone conjecture holds for $Y$, and therefore for every smooth deformation of $Y$.
\end{proof}
Note that by \cite[Theorem 4.4.2]{CoxKatz} $h^{2,1}({Y})=60$. Therefore, ${Y}$ has a large deformation space. The Borcea--Voisin threefolds in the above double cover construction only vary in a $20$-dimensional family. Theorem~\ref{thm:BorceaVoisin} therefore proves the movable cone conjecture for a large family of examples.
	\bibliography{MovingCone4}
	\bibliographystyle{alpha}
\end{document}